\documentclass[10pt,bezier,amstex]{article}  
\usepackage{color}              
\usepackage{graphicx}
\usepackage[]{amsmath}
\usepackage{amssymb}
\usepackage{hyperref,bbm}
\definecolor{MyDarkBlue}{rgb}{0,0.08,0.50}
\definecolor{BrickRed}{rgb}{0.65,0.08,0}

\hypersetup{
colorlinks=true,       
    linkcolor=MyDarkBlue,          
    citecolor=BrickRed,        
    filecolor=red,      
    urlcolor=cyan           
}

\usepackage{enumerate}
\usepackage[mathscr]{eucal}
\usepackage{graphicx}
\usepackage{latexsym}
\usepackage{amssymb}
\usepackage{psfrag}
\usepackage{epsfig}
\usepackage{tikz}
\usetikzlibrary{automata,positioning}

\usepackage{mhchem}

\include{psbox}

\usepackage{amsfonts}

\usepackage{graphicx}

\usepackage{amsfonts}

\usepackage{color}

\newtheorem{Theorem}{Theorem}
\newtheorem{Assumption}{Assumption}
\newtheorem{Lemma}{Lemma}[section]
\newtheorem{lemma}{Lemma}[section]
\newtheorem{Proposition}[Lemma]{Proposition}

\newtheorem{Corollary}[Lemma]{Corollary}

\newtheorem{remark}[Lemma]{Remark}

\newcommand{\skp}{\vspace{\baselineskip}}
\newcommand{\noi}{\noindent}

\newcommand{\s}{\quad}
\newcommand{\ep}{\varepsilon}

\newcommand{\R}{\mathbb{R}}

\newcommand{\non}{\nonumber}

\newcommand{\beq}{\begin{eqnarray*}}
\newcommand{\eeq}{\end{eqnarray*}}
\newcommand{\beqn}{\begin{eqnarray}}
\newcommand{\eeqn}{\end{eqnarray}}

\newcommand{\bt}{\begin{Theorem}}
\newcommand{\et}{\end{Theorem}}

\newcommand{\bas}{\begin{Assumption}}
\newcommand{\eas}{\end{Assumption}}
\newcommand{\lt}{\left}
\newcommand{\rt}{\right}

\newcommand{\be}{\begin{equation}}
\newcommand{\ee}{\end{equation}}

\newcommand{\eps}{\epsilon}

\newcommand{\clc}{\mathcal{C}}

\newcommand{\Lip}{\mbox{Lip}}
\newcommand{\clw}{\mathcal{W}}

\newcommand{\clp}{\mathcal{P}}

\newcommand{\cle}{\mathcal{E}}
\newcommand{\clf}{\mathcal{F}}
\newcommand{\clh}{\mathcal{H}}
\newcommand{\N}{\mathbb{N}}

\newcommand{\vs}{\varsigma}
\newcommand{\qed}{\hfill $\square$}
\newcommand{\ti}{\tilde}

\newcommand{\e}{{\varepsilon}}
\newcommand{\gr}{\nabla}

\newcommand{\nn}{\nonumber}

\numberwithin{equation}{section}

\definecolor{darkgreen}{rgb}{0,.4,0}
\definecolor{darkagenta}{rgb}{.5,0,.5}
\definecolor{darkred}{rgb}{1,0,0}
\definecolor{darkblue}{rgb}{0,0,.4}

\begin{document}

\tikzset{every node/.style={auto}}
 \tikzset{every state/.style={rectangle, minimum size=0pt, draw=none, font=\normalsize}}
 \tikzset{bend angle=20}

	\author{Abhishek Pal Majumder
	\thanks{University of Copenhagen
	}
      }

\title{Quantitative evaluation of an active Chemotaxis model in  Discrete time.}
\maketitle
\begin{abstract}
\noi
A system of $N$ particles in a chemical medium in $\mathbb{R}^{d}$ is studied in a discrete time setting. Underlying interacting particle system in continuous time can be expressed as  
\begin{eqnarray}
dX_{i}(t) &=&[-(I-A)X_{i}(t) + \bigtriangledown h(t,X_{i}(t))]dt + dW_{i}(t), \,\, X_{i}(0)=x_{i}\in \mathbb{R}^{d}\,\,\forall i=1,\ldots,N\nonumber\\
\frac{\partial}{\partial t} h(t,x)&=&-\alpha h(t,x) + D\bigtriangleup h(t,x) +\frac{\beta}{n} \sum_{i=1}^{N} g(X_{i}(t),x),\quad h(0,\cdot) = h(\cdot).\label{main}
\end{eqnarray}   
where  $X_{i}(t)$ is the location of the $i$th particle at time $t$ and $h(t,x)$ is the function measuring the concentration of the medium at location $x$ with $h(0,x) = h(x)$. In this article we describe a general discrete time non-linear formulation of the model (\ref{main}) and a strongly coupled particle system approximating it. Similar models have been studied before (Budhiraja et al.(2010)) under a restrictive compactness assumption on the domain of particles. 
In current work the particles take values in $\R^{d}$ and consequently the stability analysis is particularly challenging. We provide sufficient conditions for the existence of a unique fixed point for the dynamical system governing the large $N$ asymptotics of the particle empirical measure. We also provide uniform in time convergence rates for the particle empirical measure to the corresponding limit measure under suitable conditions on the model.

\noi {\bf AMS 2010 subject classifications:} Primary 60J05, 60K35, 60F10.\\ \ 

\noi {\bf Keywords:} Weakly interacting particle system, propagation of chaos,  nonlinear Markov chains, Wasserstein distance, McKean-Vlasov equations, exponential concentration estimates, transportation inequalities, metric entropy, stochastic difference equations, long time behavior, uniform concentration estimates.
\end{abstract}

\section{Introduction}
 There have been a surge of significant research activities aimed towards understanding the dynamics of  collective behavior of a multi-agent system in the time limit. Motivations for such problems come from various examples of self organizing systems such as consensus formation in opinion dynamics \cite{gomez2012bounded}, active chemotaxis \cite{budhiraja2011discrete}, self organized networks \cite{latane1997self}, large communication systems \cite{graham2009interacting}, multi target tracking \cite{caron2011particle}, swarm robotics \cite{pace2013swarm} (additional applications can be found in \cite{schweitzer2007brownian}) etc. One of the basic challenges is to understand how a large group of autonomous agents with decentralized local interactions that gives rise to a coherent behavior. 

In this paper we consider a reduced model motivated by both \cite{budhiraja2011discrete},\cite{budhiraja2014long} for a system of interacting agents in a stochastic diffusing environment, variations of which have been proposed (see \cite{budhiraja2011discrete},\cite{pace2013swarm} and references therein). Consider for each $i=1,\ldots,N$ $X_{i}(0)=x_{i}\in \mathbb{R}^{d}$
\begin{eqnarray}
dX_{i}(t) &=&\bigg[-(I-A)X_{i}(t)+ \gr h(t,X_{i}(t))+\frac{1}{N}\sum_{j=1, j\neq i}^{N}K\big(X_{i}(t),X_{j}(t)\big)\bigg]dt + dW_{i}(t), \quad\quad \label{main33}\\
\frac{\partial}{\partial t} h(t,x)&=&-\alpha h(t,x) + D\bigtriangleup h(t,x) +\frac{\beta}{N} \sum_{i=1}^{N} g(X_{i}(t),x),\quad h(0,\cdot) = h(\cdot).  \nonumber
\end{eqnarray}

Here $W_{i}, i=1,...,N$ are independent Brownian motions that drive the state process $X_{i}$ of the $N$ interacting particles. The interaction between the particles arises directly from the evolution equation (\ref{main33}) and indirectly through the underlying potential field $h$ which changes continuously according to a diffusion equation and through the aggregated input of the $N$ particles. One example of such an interaction is in Chemotaxis where cells preferentially move towards a higher chemical concentration and themselves release chemicals into the medium, in response to the local information on the environment, thus modifying the potential field dynamically over time. In this context, $h(t,x)$ represents the concentration of a chemical at time $t$ and location $x$. Diffusion of the chemical in the medium is captured by the Laplacian in (\ref{main33}) and the constant $\alpha>0$ models the rate of decay or dissipation of the chemical. The first equation in (\ref{main33}) describes the motion of a particle in terms of diffusion process with a drift consisting of three terms. The first term models a restoring force towards the origin where origin represents the natural rest state of the particles. The second term is the gradient of the chemical concentration and captures the fact that particles tend to move particularly towards regions of higher chemical concentration. Finally the third term captures the interaction(e.g attraction or repulsion) between the particles. Contribution of the agents to the chemical concentration field is given through the last term in the second equation. The function $g$ captures the agent response rules and can be used to model a wide range of phenomenon \cite{schweitzer2007brownian}.  

In \cite{budhiraja2011discrete} the authors considered a discrete time model which captures some of the key features of the dynamics in (\ref{main33}) and studied several long time properties of the system. One aspect that greatly simplified the analysis of \cite{budhiraja2011discrete} is that the state space of the particles is taken to be a compact set in $\R^{d}$. However this requirement is restrictive and may be unnatural for the time scales at which the particle evolution is being modeled. In \cite{pace2013swarm} authors had considered a  number of variations of (\ref{main33}). The theoretical properties obtained in this work on the long time behavior of the particle system can be also  applied for such systems with some minor modifications.

\vspace{0.2cm}
We now give a general description of the $N$- particle system that gives a discrete time approximation of the mechanism outlined above. The space of real valued bounded measurable functions on $S$ is denoted as $BM(S)$. Borel $\sigma$ field on a metric space will be denoted as $\mathcal{B}(S)$. $\mathcal{C}_{b}(S)$ denotes the space of all bounded and continuous functions $f: S \to \mathbb{R}$. For a measurable space S, $\mathcal{P}(S)$ denotes the space of all probability measures on $S$. For $k\in \mathbb{N},$ let $\mathcal{P}_{k}(\mathbb{R}^{d})$ be the space of $\mu \in \mathcal{P}(\mathbb{R}^{d})$ such that $$\|\mu\|_{k}:= \lt(\int|x|^{k}d\mu(x)\rt)^{\frac{1}{k}} <\infty.$$ 
Consider a system of $N$ interacting particles that evolve in $\mathbb{R}^{d}$ governed by a random dynamic chemical field according to the following discrete time stochastic evolution equation given on some probability space $(\Omega,\mathbb{F},P)$. Suppose that the chemical field at time instant $n$ is 
 given by a nonnegative $C^{1}$(i.e continuously differentiable)  real function on $\R^{d}$ satisfying $\int_{\R^{d}}\eta(x) dx=1$. Then, given that particle state at time instant $n$ is $x$ and the empirical measure of the particle states at time $n$ is $\mu,$  the particle state $X^{+}$ at time $(n+1)$ is given as
\begin{eqnarray}\label{main34}
X^{+} = Ax + \delta f(\gr \eta (x),\mu,x,\epsilon) + B(\epsilon),\label{tr}
\end{eqnarray}
where $A$ is a $d\times d$ matrix,  $\delta $ is a small parameter, $\epsilon$ is a $\mathbb{R}^{m}$ valued random variable with probability law $\theta$ and $f :\mathbb{R}^{d}\times\mathcal{P}(\R^d)\times\mathbb{R}^{d}\times\mathbb{R}^{m} \longrightarrow \mathbb{R}^{d}$ is a measurable function. 
Here we consider a somewhat more general form of dependence of the particle evolution on the concentration profile than the additive form that appears in (\ref{main33}). Additional assumptions on $A,\theta, f$ will be introduced shortly. Nonlinearity (modeled by $f$ and $B$) of the system can be very general and as described below. Denote by $X_{n}^{i}\equiv X_{n}^{i,N}$ (a $\mathbb{R}^{d}$ valued random variable) the state of the $i$-th particle $(i=1,\ldots,N)$ and by $\eta_{n}^{N}$ the chemical concentration field  at time instant $n$. Let $\mu_{n}^{N}:= \frac{1}{N}\sum_{i=1}^N \delta_{X_{n}^{i}}$ be the empirical measure of the particle values at time instant $n$. The stochastic evaluation equation for the $N$-particle system is given as 
\begin{eqnarray}
X_{n+1}^{i} &=& AX_{n}^{i} + \delta f(\gr \eta_{n}^{N}(X_{n}^{i}),\mu_{n}^{N}, X_{n}^{i},\epsilon_{n+1}^{i})+B(\epsilon_{n+1}^{i}),\text{\hspace{10 mm}} \quad i=1,\ldots,N, \quad n\in\mathbb{N}_{0}.\label{nlps2}
\end{eqnarray}
 In (\ref{nlps2}) $\{\epsilon_{n}^{i},i=1,...,N, \quad n\geq 1 \}$ is an i.i.d array of $\mathbb{R}^{m}$ valued random variables with common probability law $\theta$. Here $\{X_{0}^{i},i=1,...,N\}$ are assumed to be exchangeable with common distribution $\mu_{0}$ where $\mu_{0} \in \clp_{1}(\R^{d})$. Note that in the notation we have suppressed the dependence of the sequence $\{X_{n}^{i}\}$ on $N$. 

\vspace{0.1cm}

We now describe the evolution of the chemical field approximating the second equation in \eqref{main33} and its interaction with the particle system.  A transition probability kernel on $S$ is a map $P: S\times \mathcal{B}(S) \to [0,1]$ such that $P(x,\cdot)\in \mathcal{P}(S)\quad \forall x\in S$ and for each $A\in \mathcal{B}(S),$ $P(\cdot,A) \in BM(S)$. Given the concentration profile at time $n$ is a $C^{1}$ probability density function $\eta$ on $\R^{d}$ and the empirical measure of the state of $N$-particles at time instant $n$ is $\mu$, the concentration probability density $\eta^{+}$ at time $(n+1)$ is given by the relation
\beqn
\eta^{+}(y)= \int_{\R^{d}}\eta(x)R_{\mu}^{ \alpha}(x,y)l(dx)\label{evoconc}
\eeqn
where $l$ denotes the Lebesgue measure on $\R^{d},$ and $R^{\alpha}_{\mu}(x,y)$ is the Radon-Nikodym derivative of the transition probability kernel with respect to the Lebesgue measure $l(dy)$ on $\R^{d}.$ The kernel $R^{\alpha}_{\mu}$ is given as follows. We considered the same model as introduced in \cite{budhiraja2011discrete}. Let $P$ and $P'$ betwo transition probability kernels  on $\R^{d}.$ For $\mu \in \mathcal{P}(\R^{d})$ and $ \alpha\in(0,1)$ define the transition probability kernel $R_{\mu}^{\alpha}$ on $\R^{d}$ as 
$$R_{\mu}^{\alpha}(x,C) := (1-\alpha)P(x,C) + \alpha \mu P'(C),\quad\quad x\in \R^{d},C\in \mathcal{B}(\R^{d}).$$ 
Here $P$ represents the background diffusion of the chemical concentration while $\delta_{x}P'$ captures the contribution to the field by a particle with location $x$. So the kernel $P'$ gives a spike at origin which can be approximated by a smooth density function as $P(x,dy)=\frac{1}{\sqrt{2\pi}\lambda}e^{-\frac{(x-y)^2}{2\lambda^{2}}}dy$ with very small $\lambda>0$. The parameter $\alpha$ gives a convenient way for combining  the contribution from the background diffusion and the individual particles. For each $x\in\R^{d},$ both $P(x,\cdot)$ and $P'(x,\cdot)$ are assumed to be absolutely continuous with respect to Lebesgue measure and throughout this article we will denote the corresponding Radon-Nikodym derivatives with the same notations $P(x,\cdot)$ and $P'(x,\cdot)$ respectively. Additional properties of $P$ and $P'$ will be specified shortly. The evolution equation for the chemical field is then given as 
\beqn
\eta_{n+1}^{N}(y)= \int_{\R^{d}}\eta_{n}^{N}(x)R_{\mu_{n}^{N}}^{\alpha}(x,y)l(dx)\label{evoconcpart}.
\eeqn

In contrast to the model studied in \cite{budhiraja2014long}, the situation here is somewhat more involved. Note that $\{X_{n}(N)\}_{n\ge 0}:=(X_{n}^{1,N},X_{n}^{2,N},\ldots,X_{n}^{N,N})_{n\ge 0}$   is not a Markov process and in order to get a Markovian state descriptor one needs to consider $\{X_{n}(N),\eta_{n}^{N}\}_{n\ge 0}$ which is a discrete time Markov chain with values in $(\R^{d})^{N}\times \clp_{}(\R^{d})$.

\vspace{0.2cm}

We will show that as $N\to\infty$ $(\mu^{N}_{n},\eta^{N}_{n})_{n\in\mathbb{N}_{0}}$ converges to a deterministic nonlinear dynamical system $(\mu_{n},\eta_{n})_{n\in\mathbb{N}_{0}}$ with methods followed in \cite{budhiraja2011discrete}. We established further sharp quantitative bounds (with techniques used in \cite{fournier2013rate} and \cite{budhiraja2014long}) for weakly interacting particle system jointly with the stochastic field potential to the nonlinear system of interest. 
For both polynomial and exponential concentration bound it requires further constraints on the tail of the transition kernels $P,P'$ used in modeling the  diffusive environment. One major motivation of cthe current article is giving a sharp uniform in time quantitative estimate for the particle system $(\mu_{n}^{N},\eta_{n}^{N})$ to the non-linear system of interest $(\mu_{n},\eta_{n})$ so that any functional of the form $\big<\phi_{1},\mu_{n}^{}\big>+\big<\phi_{2},\eta_{n}\big>$ can be approximated by $\frac{1}{N}\sum_{i=1}^{N}\phi_{1}(X^{i}_{n})+\big<\phi_{2},\eta_{n}^{N}\big>$ with desired precision. Previous work on concentration bounds for similar particle system in discrete time includes \cite{del2011concentration} but that involves a Dobrushin type
 stability condition which is not very effective if the particles are assumed to come from a non-compact domain. A very recent work \cite{budhiraja2016uniform} addresses several quantitative bounds for Chemotaxis model motivated by Patlak-keller-segel type non-linear equations.

\vspace{0.2cm}

The following notations will be used in this article. $\mathbb{R}^{d}$ will denote the $d$ dimensional Euclidean space with the usual Euclidean norm $|\cdot|$. The set of natural numbers (resp. whole numbers) is denoted by $\mathbb{N}$ (resp. $\mathbb{N}_{0}$). Cardinality of a finite set $S$ is denoted by $|S|$.  For $x\in \mathbb{R}^{d}$, $\delta_{x}$ is the Dirac delta measure on $\mathbb{R}^{d}$ that puts a unit mass at location $x$. The supremum norm of a function $f:S\to \mathbb{R}$ is $\|f\|_{\infty}=\sup_{x\in S}|f(x)|$. When $S$ is a metric space, the Lipschitz seminorm of $f$ is defined by $\|f\|_{1}=\sup_{x\not= y}\frac{|f(x)-f(y)|}{d(x,y)}$ where $d$ is the metric on the space $S$. For a bounded Lipschitz function $f$ on $S$ we define $\|f\|_{BL}:=\|f\|_{1}+\|f\|_{\infty}$.
$\mbox{Lip}_{1}(S)$ (resp. $BL_{1}(S)$ )  denotes the class of Lipschitz (resp. bounded Lipschitz) functions $f :S\to \mathbb{R}$ with $\|f\|_{1}$ (resp. $\|f\|_{BL}$)
bounded by 1. Occasionally we will suppress $S$ from the notation and write $\mbox{Lip}_{1}$ and  $BL_{1}$ when clear from the context.  For a Polish space $S$, $\mathcal{P}(S)$ is equipped with the topology of weak convergence. A convenient metric metrizing this topology on $\mathcal{P}(S)$ is given as $\beta(\mu,\gamma) = \sup \{|\int fd\mu - \int fd\gamma|:\|f\|_{BL_{1}}\leq 1 \}$ for $ \mu,\gamma\in  \mathcal{P}(S)$. For a signed measure $\gamma$ on $\mathbb{R}^{d}$, we define $\langle f,\gamma\rangle:= \int f d\gamma$ whenever the integral makes sense. The space $\mathcal{P}_{1}(\mathbb{R}^{d})$ will be equipped with the Wasserstein-1 distance that is defined as follows:
$$\mathcal{W}_{1}(\mu_{0},\gamma_{0})  := \inf_{X,Y}E|X-Y|,\s\s \mu_0,\nu_0\in \mathcal{P}_{1}(\R^d),$$ where the infimum is taken over all  $\mathbb{R}^{d}$ valued random variables $X,Y$ defined on a common probability space and where the marginals of $X, Y$ are $\mu_{0}$ and $\gamma_{0}$ respectively.
From Kantorovich-Rubenstein duality (cf. \cite{villani2003topics}) one sees the Wasserstein-1 distance has the following characterization
\begin{eqnarray}
\mathcal{W}_{1}(\mu_{0},\gamma_{0}) =\sup_{f\in \mbox{Lip}_{1}(\mathbb{R}^{d})}|\langle f,\mu_{0} - \gamma_{0}\rangle|,\s\s\mu_0,\nu_0\in \mathcal{P}_{1}(\R^d).
\end{eqnarray}
For a signed measure $\mu$ on $(S,\mathcal{B}(S))$, the total variation norm of $\mu$ is defined as $|\mu|_{TV}:=\sup_{||f||_{\infty} \leq 1}\langle f,\mu \rangle$. Probability distribution of a $S$ valued random variable $X$ will be denoted as $\mathcal{L}(X)$. Convergence in distribution of a $S$ valued  sequence $\{X_{n}\}_{n\geq 1}$ to a $S$ valued random variable $X$ will be written as $X_{n}\Rightarrow X$. 

A finite collection $\{Y_{1},Y_{2},\ldots,Y_{N}\}$ of $S$ valued random variables is called exchangeable if $$\mathcal{L}(Y_{1},Y_{2},\ldots,Y_{N}) = \mathcal{L}(Y_{\pi(1)},Y_{\pi(2)},\ldots,Y_{\pi(N)})$$ for every permutation $\pi$ on the $N$ symbols $\{1,2,\ldots,N\}$.  Let $\{Y_{i}^{N},i=1,\ldots,N\}_{N\geq 1}$ be a collection of $S$ valued random variables, such that for every $N$, $\{Y_{1}^{N},Y_{2}^{N},\ldots,Y_{N}^{N}\}$ is exchangeable. Let $\nu_{N} = \mathcal{L}(Y_{1}^{N},Y_{2}^{N},\ldots,Y_{N}^{N})$. The sequence $\{\nu_{N}\}_{N\geq 1}$ is called $\nu$ -chaotic (cf.  \cite{sznitman1991topics})  for a $\nu \in \mathcal{P}(\mathcal{S})$, if for any $k\geq 1$, $f_{1},f_{2},\ldots,f_{k} \in \mathcal{C}_{b}(\mathcal{S}),$ one has
\begin{eqnarray}
\lim_{N\to \infty} \langle f_{1}\otimes f_{2}\otimes \ldots\otimes f_{k}\otimes 1 \ldots \otimes1,\nu_{N} \rangle = \prod_{i=1}^{k}\langle f_{i} ,\nu \rangle. \label{chaotic}
\end{eqnarray}

Denoting the marginal distribution on first $k$ coordinates of $\nu_{N}$ by $\nu_{N}^{k}$, equation (\ref{chaotic}) says that, for every $k\geq 1,$ $\nu_{N}^{k} \rightarrow \nu^{\otimes k}$.
 The gradient of a real differentiable function $f$ on $\R^{d}$ denoted by $\gr f$ is defined as the $d$ dimensional vector field $\gr f := (\frac{\partial f}{\partial x_{1}},\frac{\partial f}{\partial x_{2}},\ldots,\frac{\partial f}{\partial x_{d}})'$. For a function $f:\R^{d}\times \R^{m}\to \R$ $$\gr_{x}f(x,y) :=  \lt(\frac{\partial f}{\partial x_{1}},\frac{\partial f}{\partial x_{2}},\ldots,\frac{\partial f}{\partial x_{d}}\rt)'.$$ The function $\gr_{y}f(x,y)$ is defined similarly. Absolute continuity of a measure $\mu$ with respect to a measure $\nu$ will be denoted by $\mu \ll \nu.$ We will denote the Radon-Nikodym derivative of $\mu$ with respect to $\nu$ by $\frac{d\mu}{d\nu}$. For $f\in BM(\mathcal{S})$ and a transition probability kernel $P$ on $S$, define $Pf \in  BM(\mathcal{S})$ as $Pf(\cdot) = \int_{S} f(y)P(\cdot,dy) $. For any closed subset $B \in S$, and $\mu \in \mathcal{P}(B),$ define $\mu P\in \mathcal{P}(S)$ as $\mu P(A) = \int_{B} P(x,A)\mu(dx)$. For a matrix $B$ the usual operator norm is denoted by $\|B\|$. 

\section{Description of the nonlinear system:}
\label{sec:mod-desc}

We now describe the nonlinear dynamical system obtained on taking the limit $N\to \infty$ of $(\mu_{n}^{N},\eta_{n}^{N})$.
 Given a $C^1$ density function $\rho$ on $\R^{d}$ and $\mu\in\mathcal{P}(\R^d)$, define a transition probability kernel $Q^{\rho,\mu}$ on $\R^{d}$ as  $$Q^{\rho,\mu}(x,C) = \int_{\R^{m}} 1_{\{Ax+\delta f(\gr{\rho}(x),\mu,x,z)+B(z)\in C\}} \theta(dz),\quad  \quad (x,C)\in \R^{d} \times \mathcal{B}(\R^{d}).$$

With an abuse of notation we will also denote by $Q^{\rho,\mu}$ the map from $BM(\R^{d})$ to itself, defined as 
$$Q^{\rho,\mu} \phi(x) = \int_{\R^{d}} \phi(y) Q^{\rho,\mu}(x,dy), \quad \phi \in BM(\R^{d}), x\in \R^{d}.$$

For $\mu,\mu_{1} \in \clp(\R^{d})$, let $\mu Q^{\rho,\mu_{1}}\in \mathcal{P}(\R^{d})$ be defined as
\beqn
\mu Q^{\rho,\mu_1}(C) = \int_{\R^{d}} Q^{\rho,\mu_1}(x,C) \mu(dx), \quad C\in \mathcal{B}(\mathbb{R}^{d}).\label{trans}
\eeqn
Note that $\mu_{} Q^{\rho,\mu_{1}} = \mathcal{L}\big(AX + \delta f(\gr \rho(X),\mu_1,X,\epsilon)+B(\epsilon)\big)$ where $\mathcal{L}(X,\epsilon) = \mu_{} \otimes \theta$.

Define $\clp_{1}^{*}(\R^{d}):=\{\mu \in \clp_{1}(\R^{d}): \mu \ll  l , \frac{d \mu}{d l}$  is continuously differentiable and $\|\gr \frac{d \mu}{d l}\|_{1} <\infty\}.$ For notational simplicity we will identify an element in $\clp_{1}^{*}(\R^{d})$ with its density and denote both by the same symbol. Define the map $\Psi : \clp(\R^{d})\times\clp_{1}^{*}(\R^{d})\to \clp(\R^{d})\times\clp(\R^{d})$  as
\beqn
\Psi(\mu,\eta) &=& (\mu Q^{\eta,\mu},\eta R_{\mu}^{ \alpha}), \s\s (\mu,\eta)\in \clp(\R^{d})\times\clp_{1}^{*}(\R^{d}). \label{vlasovm}
\eeqn
Under suitable  assumptions (which will be introduced in Section \ref{Assumptionsec}) it will follow that for every $(\mu,\eta)\in \clp_{1}(\R^{d})\times\clp_{1}^{*}(\R^{d}),$ $\eta^{+}$ defined by (\ref{evoconc}) is in $\clp_{1}^{*}(\R^{d})$ and $\mu Q^{\eta,\mu}$ defined by (\ref{trans}) is in $\clp_{1}(\R^{d})$. Thus (under those assumptions) $\Psi$ is a map from $ \clp_{1}(\R^{d})\times\clp_{1}^{*}(\R^{d})$ to itself. Using the above notation we see that $\{(X_{n}^{1},...,X_{n}^{N}),\mu_{n}^{N},\eta_{n}^{N}\}_{n\geq 0}$  is a $\mathbb(\mathbb{R}^{d})^{N} \times\mathcal{P}_{1}(\mathbb{R}^{d})\times \clp_{1}^{*}(\R^{d})$ valued discrete time Markov chain  defined recursively as follows. Let $X_{k}(N) \equiv (X_{k}^{1},X_{k}^{2},...,X_{k}^{N})$, and  $\eta_{0}^{N}$ be the initial chemical field which is a random element of $\clp_{1}^{*}(\R^{d})$. Let $\mathcal{F}_{0} = \sigma\{X_{0}(N),\eta_{0}^{N}\}.$ Then, for $k\geq 1$
\begin{eqnarray}\label{total22}
\begin{cases}
P(X_{k}(N) \in C |\mathcal{F}_{k-1}^{N}) = \bigotimes_{i=1}^{N}(\delta_{X_{k-1}^{j}} Q^{\eta_{k-1}^{N},\mu_{k-1}^{N}})(C)\hspace{3mm} \forall C \in \mathcal{B}(\mathbb{R}^{dN}),  \\
\mu_{k}^{N}= \frac{1}{N}\sum_{i=1}^N \delta_{X_{k}^{i}}, \\
\eta_{k}^{N}= \eta_{k-1}^{N} R^{ \alpha}_{\mu_{k-1}^{N}}, \\
\mathcal{F}_{k}^{N} =\sigma\{\eta_{k}^{N},X_{k}(N)\} \vee \mathcal{F}_{k-1}^{N}. 
\end{cases}
\end{eqnarray}
We will call this particle system as $\mathbb{IPS}_{1}$. We next describe a nonlinear dynamical system which is the formal Vlasov-Mckean limit of the above system, as $N\to \infty$. Given $(\mu_{0},\eta_{0}) \in \clp_{1}(\R^{d})\times \clp_{1}^{*}(\R^{d})$ define a sequence $\{(\mu_{n},\eta_{n})\}_{n\ge 0}$ in $\clp_{1}(\R^{d})\times \clp_{1}^{*}(\R^{d})$ as
\beqn
\mu_{n+1}= \mu_{n}Q^{\eta_{n},\mu_{n}},\s\s\ \eta_{n+1}= \eta_{n}R^{\alpha}_{\mu_{n}},\s\s\s n\ge 0.\label{extra}
\eeqn 
Using (\ref{vlasovm}) the above evolution can be represented as 
\beqn
(\mu_{n+1},\eta_{n+1})=\Psi(\mu_{n},\eta_{n}),\s\s\s n\in \mathbb{N}_{0}\label{Vlasovmseq}.
\eeqn 
As in \cite{budhiraja2014long}, the starting point of our investigation on long time asymptotics of the above interacting particle system will be to study the stability properties of (\ref{extra}).
We identify $\eta,\eta' \in \mathcal{P}(\R^d)$ that are equal a.e under the Lebesgue measure on $\R^{d}$. From a computational point of view we are approximating $(\mu_n,\eta_{n})$ by $(\mu_n^N,\eta_{n}^{N})$  uniformly in time parameter $n$, with explicit uniform concentration bounds. Such results are particularly important for developing sampling methods for approximating the steady state distribution of the mean field models such as in (\ref{extra}).

\vspace{0.3cm}

The third equation in \eqref{total22} makes the simulation of $\mathbb{IPS}_{1}$ numerically challenging. In section \ref{Assumptionsec} we will mention another particle system (based on the second particle system in \cite{budhiraja2011discrete}) referred to as $\mathbb{IPS}_{2}$ which also gives an asymptotically consistent approximation of \eqref{extra} and is computationally more tractable. We show in THeorem \ref{pointwise} that under conditions that include a Lipschitz property of $f$ (Assumptions \ref{As2} and \ref{As2.5}), smoothness assumptions on the transition kernels of the background diffusion of the chemical medium (Assumption \ref{As3}) the Wasserstein-1($\mathcal{W}_1$) distance between the occupation measure of the particles along with the chemical medium $(\mu_{n}^{N},\eta_{n}^{N})$  and $(\mu_{n},\eta_{n})$ converges to $0$, for every time instant $n.$ Under  an additional condition on the contractivity of $A$ and $\delta, \alpha$ being sufficiently small we show that the nonlinear system \eqref{Vlasovmseq} has a unique fixed point and starting from an arbitrary initial condition, convergence to the fixed point occurs at a geometric rate.  Using these results we next argue in Theorem \ref{unifintime} that under some integrability conditions  (Assumption \ref{unifas1}-\ref{unifas2}), as $N \to \infty$, the empirical occupation measure of the $N$-particles and density of the chemical medium at time instant $n$, namely $(\mu_n^N,\eta_{n}^{N})$ converges to $(\mu_n,\eta_{n})$ in the $\mathcal{W}_1$ distance, in $L^1$,  {\em uniformly} in $n$. This result in particular shows that the $\mathcal{W}_1$ distance between $(\mu_n^N,\eta_{n}^{N})$ and the unique fixed point $(\mu_\infty,\eta_{\infty})$ of \eqref{Vlasovmseq} converges to zero as $n \to \infty$ and $N\to \infty$ in any order. We next show that for each $N$, there is  unique invariant measure $\Theta^N_{\infty}$ of the $N$-particle dynamics with integrable first moment and this sequence of measures is $\mu_{\infty}$-chaotic, namely as $N \to \infty$, the projection of $\Theta^N_{\infty}$ on the first $k$-coordinates converges to $\mu_{\infty}^{\otimes k}$ for every $k \ge 1$.  This propagation of chaos property all the way to $n=\infty$ crucially relies on the uniform in time convergence of $(\mu_n^N,\eta_{n}^{N})$ to $(\mu_{\infty},\eta_{\infty})$.
Such a result is important since it says that the steady state of a $N$-dimensional fully coupled Markovian system has a simple approximate description in terms of  a product measure when $N$ is large. This result is key in developing particle based numerical schemes for approximating the fixed point of the evolution equation \eqref{Vlasovmseq}. Next we present some uniform in time concentration bounds of $\clw_{1}(\mu^{N}_{n},\mu_n)+\clw_{1}(\eta_{n}^{N},\eta_{n})$. Proof is very similar to that of Theorem 3.8 of \cite{budhiraja2014long} so we only provide a sketch after showing necessary conditions. 

\section{Main Results:}\label{Assumptionsec}
We now introduce our main assumptions on the problem data. Recall that $\{X_{0}^{i},i=1,\ldots N\}$ is assumed to be exchangeable with common distribution $\mu_0.$ We assume further $(\mu_{0},\eta_{0})\in \clp_1(\R^{d})\times\clp^*_1(\R^{d}).$  For a $d\times d$ matrix B we denote its norm by $\|B\|,$ i.e. $\|B\| = \sup_{x \in \mathbb{R}^{d}\setminus\{0\}} \frac{|Bx|}{|x|}$.

\bas\label{As2}
The error distribution $\theta$ is such that $\int A_{1}(z)\theta(dz) := \sigma\in(0,\infty)$ where
\beqn
A_{1}(\epsilon)   &:=& \sup_{\{x_{1}, x_{2}, y_{1},y_{2}\in\R^{d},\mu_{1},\mu_{2}\in\mathcal{P}_{1}(\R^{d}) : \mu_{1}\neq\mu_{2},x_{1} \neq x_{2}, y_{1} \neq y_{2}\}}           \frac{       |f(y_{1},\mu_{1},x_{1},\epsilon)-f(y_{2},\mu_{2},x_{2},\epsilon)|}{|x_{1}- x_{2}|+ |y_{1}-y_{2}|+\mathcal{W}_{1}(\mu_{1},\mu_{2})}.\label{as1}
\eeqn

It follows that $\forall x,y\in\R^{d},\mu\in\mathcal{P}_{1}(\R^{d}),$ 
\beqn\label{cru}
 |f(y,\mu,x,\epsilon)| \leq (|y|+\|\mu\|_{1}+|x|)A_{1}(\epsilon) +A_{2}(\epsilon)
\eeqn
where $A_{2}(\epsilon):= f(0,0,\epsilon)$.
\eas
Recall the function $B:\R^{m}\to\R^{d}$ introduced in \eqref{tr}.
\bas\label{As2.5}
The error distribution $\theta$ is such that 
$$\int_{\R^{m}} \Big(A_{2}(z)+|B(z)|\Big)\theta(dz)<\infty.$$
\eas

\bas\label{As6}
$\eta_{0}^{N}$ (the density function) is a Lipschitz function on $\R^{d}$ and $\eta_{0}^{N} \in\clp_{1}^{*}(\R^{d})$ .
\eas

Assumptions \ref{As3} and \ref{As7} on the kernels $P$ and $P'$ hold quite generally. In particular, they are satisfied for Gaussian kernels.
\bas\label{As3}
There exist $l_{P}^{\gr}\in (0,1]$ and $l_{P'}^{\gr} \in (0,\infty)$ such that for all $x,y,x',y'\in \R^{d}$
\beqn
 |\gr_{y} P(x,y)-  \gr_{y} P(x',y')|&\leq& l_{P}^{\gr}(|y-y'| +|x-x'|)\quad \label{lip1}\\
|\gr_{y} P'(x,y)-  \gr_{y} P'(x',y')|&\leq& l_{P'}^{\gr}(|y-y'| +|x-x'|).\quad \label{lip2}
\eeqn
Furthermore 
\beqn
\sup_{x\in\R^{d}}\{|\gr_{y}P(x,0)|\vee|\gr_{y}P'(x,0)|\}<\infty.\label{lingrfin}
\eeqn
 Using the Lipschitz property in (\ref{lip1}) and the growth condition \eqref{lingrfin} one has the linear growth property for some $M_{P}^{\gr}\in(0,\infty)$ 
\beqn
sup_{x\in\R^{d}} |\gr_{y} P(x,y)|\le  M_{P}^{\gr}(1+|y|).\label{lingrfin}
\eeqn
A similar inequality holds for $P'$ from (\ref{lip2}) with $M_{P'}^{\gr}\in(0,\infty)$.
\eas
Denote $(1-\alpha) l_{P}^{\gr}+\alpha  l_{P'}^{\gr}$ by $ l_{PP'}^{\gr,\alpha}$.

\bas\label{As7}
For every $f\in \Lip_{1}(\R^{d}),$  $Pf$ and $P'f$ are also Lipschitz and

$$\sup_{f\in Lip_{1}(\R^{d})}\sup_{x\neq y \in \R^{d}}\frac{Pf(x)-Pf(y)}{|x-y|}:=l(P)<\infty$$
Also $l(P')$ defined as above for $P'$ is finite. 
\eas

\bas\label{As8unifintegr}
Both $P(x,\cdot)$ and $P'(x,\cdot)$ are such that for any compact set $K\subset \R^d,$ the families of probability measures $\{P(x,\cdot): x\in K\}$ and $\{P'(x,\cdot): x\in K\}$ are both uniformly integrable.
\eas

Let $\max\{l(P),l(P')\} =l_{PP'}$. 
\begin{remark}\label{rem1}
Assumption \ref{As7} is satisfied if $P,P'$ are given as follows. For any $f\in \clc_{b}(\R^d),$ let 
\beqn\label{tran1}
Pf(\cdot) := E f(g_{1}(\cdot,\varepsilon_1)),\s\s\s\s P'f(\cdot) := E f(g_{2}(\cdot,\varepsilon_2))
\eeqn
where $\e_{1},\e_{2}$ are $\R^{m}$ valued random variables and $\varepsilon_1,\varepsilon_2$ and $g_1,g_{2}: \R^{d}\times \R^{m}\to \R^{d}$ are maps with following properties:
\beqn\label{tran2}
E(G_1(\varepsilon_{1}))\le l(P)\s\s\s\text{and}\s\s\s E(G_2(\varepsilon_{2}))\le l(P'),
\eeqn
where
\beqn
G_1(y) :=\sup_{x_1\neq x_2}\frac{g_1(x_1,y)-g_1(x_2,y)}{|x_{1}-x_2|}\s\s\text{and}\s\s G_2(y) :=\sup_{x_1\neq x_2}\frac{g_2(x_1,y)-g_2(x_2,y)}{|x_{1}-x_2|}.
\eeqn
\end{remark}

Simulation of the system  is numerically intractable due to the step that involves the updating of $\eta_{n-1}^{N}$ to $\eta_{n}^{N}.$  This requires computing the integral in \eqref{evoconc} which, since $R_{\mu}^{\alpha}$ is a mixture of two transition kernels, over time leads to an explosion of terms in the mixture that need to be updated. An approach (proposed in \cite{budhiraja2011discrete}) that addresses this difficulty is, without directly updating $\eta_{n-1}^{N}$, to use the empirical distribution of the observations drawn independently from $\eta_{n-1}^{N}.$

Denote $\bar{X}_{0}(N)$ by $(\bar{X}_{0}^{1,N},\ldots,\bar{X}_{0}^{N,N})$ a sample of size $N$ from $\mu_0.$ Let $M\in\mathbb{N}$. The new particle scheme will be described as a family $(\bar{X}_{k}(N), \bar{\mu}^{N}_{k},\bar{\eta}^{M}_{k})_{k\in \mathbb{N}_{0}}$ of $(\R^{d})^{N}\times\clp(\R^{d}) \times\clp^*(\R^d)$ valued random elements on some probability space defined recursively as follows. Set $\bar{X}_{0}(N)=(\bar{X}_{0}^{1,N},\ldots,\bar{X}_{0}^{N,N}), \bar{\eta}^{M}_{0}=\eta_{0},\bar{\mathcal{F}}^{M,N}_{0}=\sigma(\bar{X}^N(0))$. For $k\ge 1$
\begin{eqnarray}\label{total23}
\begin{cases}
\bar{\mu}_{k}^{N}= \frac{1}{N}\sum_{i=1}^N \delta_{\bar{X}_{k}^{i}},
  \\
P(\bar{X}_{k}(N) \in C |\mathcal{F}_{k-1}^{M,N}) = \bigotimes_{i=1}^{N}(\delta_{\bar{X}_{k-1}^{j}} Q^{\bar{\eta}_{k-1}^{M},\bar{\mu}_{k-1}^{N}})(C)\hspace{3mm} \forall C \in \mathcal{B}(\mathbb{R}^{d})^{N},\\
\bar{\eta}_{k}^{M}= (1-\alpha)(S^{M}(\bar{\eta}_{k-1}^{M})P) +\alpha \bar{\mu}^{N}_{k-1}P',\\
\bar{\mathcal{F}}_{k}^{M,N} =\sigma\{\bar{\eta}_{k}^{M},\bar{X}_{k}(N)\} \vee \bar{\mathcal{F}}_{k-1}^{M,N}
\end{cases}
\end{eqnarray} 
where $S^M(\bar{\eta}_{k-1}^{M})$ is the random measure defined as $\frac{1}{M}\sum_{i=1}^M \delta_{Y^{i,M}_{k-1}}$ where $\{Y^{i,M}_{k-1}\}_{i=1,\ldots,M}$ conditionally on $\bar{\mathcal{F}}_{k-1}^{M,N},$ are $M$ i.i.d distributed according to $\bar{\eta}_{k-1}^{M}.$ We will call this particle system as $\mathbb{IPS}_{2}$. We remark that our notation is not accurate since both the quantities $\bar{\mu}_{k}^{N},\bar{\eta}_{k}^{M}$ depend on $M,N.$ The  superscripts only describe the number of particles/samples used in the procedure to combine them.  Note that like $\mathbb{IPS}_{1},$ here $(\bar{X}_{k}(N),\bar{\eta}_{k}^{M})_{k\ge 0}$ is not a Markov chain on $(\R^{d})^N\times \clp_{1}^{*}(\R^d)$ anymore. Rather $(\bar{X}_{}^{N}(k),\bar{\eta}_{k}^{M},S^{M}(\bar{\eta}_{k}^{M}))_{k\ge 0}$ is a discrete time Markov chain on $(\R^{d})^{N}\times \clp_{1}^{*}(\R^d)\times \clp_{1}^{}(\R^d)$. 

For any random variable $Z$ we denote $E\big[Z\big|\mathcal{F}_{k}^{M,N}\big]$ by $E_{k}^{M,N} \big[Z\big]$. The following result shows that the particle systems in \eqref{total22} and \eqref{total23} approximate the dynamical system in \eqref{extra} as $N$ (respectively $\min\{M,N\}$ for $\mathbb{IPS}_{2}$) becomes large for a fixed time instant.

\vspace{0.2cm}

\begin{Proposition}\label{pointwise}
Suppose Assumptions \ref{As2},\ref{As2.5},\ref{As3} and \ref{As7} hold.  
\begin{enumerate}[(a)]
\item Consider the particle system $\mathbb{IPS}_1$ in (\ref{nlps2},\ref{evoconcpart}). Suppose the sampling of the exchangeable datapoints $X_{0}(N) \equiv (X_{0}^{1},X_{0}^{2},\ldots,X_{0}^{N})$ is exchangeable and $\{\mathcal{L}(X_{0}(N))\}_{N\in\mathbb{N}}$ is $\mu_0$- chaotic. Suppose $E\clw_1(\eta_{0}^{N},\eta_0)\to 0$  as $N\to \infty$. Then, as $N\to \infty$
\beqn
E\left[\clw_1(\mu_{n}^{N},\mu_{n})+ \clw_1(\eta_{n}^{N},\eta_{n})\right]\to 0
\eeqn
for all $n\ge 0$ where $\mu_{n},\eta_{n}$ are as in (\ref{extra}).

\item Consider the second particle system $\mathbb{IPS}_2.$ Suppose that in addition Assumption \ref{As8unifintegr} holds. Suppose the sampling of the exchangeable datapoints $\bar{X}_{0}(N) \equiv (\bar{X}_{0}^{1},\bar{X}_{0}^{2},\ldots,\bar{X}_{0}^{N})$ is exchangeable and $\{\mathcal{L}(\bar{X}_{0}(N))\}_{N\in\mathbb{N}}$ is $\mu_0$- chaotic. Then as $\min\{N,M\}\to \infty$ 
\beqn\label{ptth1}
E\left[\clw_1(\bar{\mu}_{n}^{N},\mu_{n})+ \clw_1(\bar{\eta}_{n}^{M},\eta_{n})\right]\to 0
\eeqn
for all $n\ge 0$. 
\end{enumerate}
\end{Proposition}

As a consequence of Proposition \ref{pointwise}, we have a finite time propagation of chaos result of the following form. Let $\nu_{n}^{N}=\mathcal{L}(X_{n}^{1,N},X_{n}^{2,N},\ldots,X_{n}^{N,N}).$ 
\begin{Corollary}
Under Assumptions as in Proposition \ref{pointwise} the family $\{\nu_{n}^{N}\}_{N\ge 1}$ is $\mu_{n}$ chaotic for every $n\ge 1$.
\end{Corollary}

As noted in introduction, the primary goal is studying long time properties of (\ref{nlps2}) and the non-linear dynamical system (\ref{extra}). Following proposition identifies the range of values of the modeling parameters  that leads to stability of the system.

\begin{Proposition} \label{Thmstable}
Suppose Assumptions \eqref{As2}-\eqref{As7} hold. Then there exist $\omega_{0}, \alpha_{0}, \delta_{0}\in(0,1) $  such that for all $\|A\| < \omega_{0}, \alpha \in(0, \alpha_{0})$, and $ \delta \in(0,\delta_{0})$. The map $\Psi$ defined in (\ref{vlasovm}) has a unique fixed point $(\mu_{\infty},\eta_{\infty})$ in $\clp_{1}(\R^{d})\times \clp_{1}^{*}(\R^{d}).$
\end{Proposition}

Now we will give more stringrent conditions under which a non-asymptotic bound on convergence rates of the particle system to the deterministic nonlinear dynamics and  their  consequences for the steady state behavior can be established. 

\bas\label{unifas1}
For some $\tau>0,$ 
$$\mu_{0} \in \clp_{1+\tau}(\R^d),\s\s \int A_{1}(z)^{1+\tau}\theta(dz):=\sigma_1(\tau)<\infty\s\s\int \Big(A_{2}(z)+|B(z)|\Big)^{1+\tau}\theta(dz):=\sigma_2(\tau)<\infty.$$
\eas
 We need to impose the following condition on $P,P'$ for uniform in time convergence.
\bas\label{unifas2}
For some $\left <|x|^{1+\tau},\eta_{0}\right> <\infty.$  There exist $m^{}_{\tau}(P)$ and $m_{\tau}(P')$ in $\R^{+}$ such that following holds for all $x\in\R^d$
\beqn
\int_{\R^{d}} |y|^{1+\tau}P(x,dy)\le m^{}_{\tau}(P) \lt(1+|x|^{1+\tau}\rt),\s \text{and}\s \int_{\R^{d}} |y|^{1+\tau}P'(x,dy)\le m^{}_{\tau}(P') \lt(1+|x|^{1+\tau}\rt).\non
\eeqn

\eas

Now we state a generalization of the Proposition \ref{pointwise}, which gives the convergence rate of $$E\left\{\clw_1(\bar{\mu}_{n}^{N},\mu_{n})+ \clw_1(\bar{\eta}_{n}^{M},\eta_{n})\right\} \to 0$$ uniformly over all $n\ge 0$ in a nonasymptotic manner.

Recall $l_{P}^{\gr}, l_{P'}^{\gr}$ introduced in Assumption \ref{As6}. For $\alpha\in(0,1),$ let $l_{PP'}^{\gr,\alpha}=(1-\alpha)l_{P}^{\gr}+\alpha l_{P'}^{\gr}$. With the notations of Assumption \ref{As2} we define
$$a_{0}:= \frac{1-\|A\|}{\sigma(2+l^{\gr,\alpha}_{PP'})}.$$

For $(\mu_{n},\eta_{n}),(\mu'_{n},\eta'_{n})\in \clp_{1}(\R^{d})\times\clp_{1}^{*}(\R^{d})$ define the following distance on $\clp_{1}(\R^{d})\times\clp_{1}^{*}(\R^{d})$ $$\mathcal{W}_{1}((\mu_{n},\eta_{n}),(\mu'_{n},\eta'_{n})):= \mathcal{W}_{1}(\mu_{n},\mu'_{n})+\mathcal{W}_{1}(\eta_{n},\eta'_{n}).$$

\begin{Theorem}\label{unifintime}
Consider the particle system $\mathbb{IPS}_{2}$. Suppose Assumptions \eqref{As2}-\eqref{As7} and Assumptions (\ref{unifas1}),(\ref{unifas2}) hold for some $\tau>0$. Let $N_{1}:=\min\{M,N\}.$ Also assume $\delta \in (0,a_{0}),\s (1-\alpha)m_{\tau}(P)<  1$ and
$$\max\Big\{\Big(\|A\|+\delta\sigma(2 +l_{PP'}^{\gr,\alpha}) + \alpha l(P')\Big ),(1-\alpha)l(P) \Big\}+\delta \sigma \max\big\{\alpha l_{P'}^{\gr},(1-\alpha) l_{P}^{\gr}\big\}< 1,\s\text{}\s .$$
Then there exists $\theta <1,$ and $a\in(0,\infty)$ such that for each $n\ge 0,$ the upperbound $b(N_1,\tau,d)$ of $$E\clw_1\lt((\bar{\mu}_{n}^{N},\bar{\eta}_{n}^{M}),(\mu_{n},\eta_{n}) \rt)-a\theta^{n}E\clw_1\lt((\bar{\mu}_{0}^{N},\bar{\eta}_{0}^{M}),(\mu_{0},\eta_{0})\rt)$$ can be expressed as
\beqn
b(N_1,\tau,d)= C
\begin{cases}
N_{1}^{- \max\{\frac{1}{2},\frac{\tau}{1+\tau}\}}\s\s\s\s\s\s\text{if}\s\s  d=1,\tau \neq 1,\s \\
N_{1}^{-\frac{1}{2}}\log N_{1}\s\s\s\s\s\s\s\,\text{if}\s\s  d=1,\tau = 1,\s \\
N_{1}^{-\frac{1}{2}}\log N_{1}+N_{1}^{-\frac{\tau}{1+\tau}}\s\s\s\text{if}\s\s d=2,\tau\neq 1,\s\s\s\s\s\s\\
N_{1}^{-\frac{1}{2}}(\log N_{1})^{2}\s\s\s\s\s\s\,\text{if}\s\s d=2,\tau = 1,\s \\
N_{1}^{-\max\{\frac{1}{d},\frac{\tau}{1+\tau}\}}\s\s\s\s\s\,\,\,\,\,\,\text{if}\s\s d>2,\tau \neq \frac{1}{d-1},\\
N_{1}^{-\frac{1}{d}}\log N_{1}\s\s\s\s\s\s\s\,\text{if}\s\s d>2,\tau = \frac{1}{d-1},
\end{cases}.
\eeqn
where the value of the constant $C$ will vary for each of the cases.

\end{Theorem}
\begin{remark}
For the first particle system (\ref{nlps2}-\ref{evoconcpart}) similar results hold where the explicit bounds are given in terms of number of particles $N$ instead of $N_1.$ For $\mathbb{IPS}_{2}$ if the initial sampling scheme of $\bar{X}_{0}(N) \equiv (\bar{X}_{0}^{1},\bar{X}_{0}^{2},...,\bar{X}_{0}^{N})$ is $\mu_0$ -chaotic then using the fact $E\mathcal{W}_1(\bar{\mu}_{0}^{N},\mu_{0})\to 0$ as $N\to \infty,$ it follows from the conclusion of the Theorem \ref{unifintime}
$$\sup_{n \ge 0} E\mathcal{W}_{1}\lt((\bar{\mu}_{n}^{N},\bar{\eta}_{n}^{M}),(\mu_{n},\eta_{n})\rt)\to 0$$
as $\min{\{N,M\}}\to \infty.$ For the first particle system in (\ref{nlps2}-\ref{evoconcpart}), if $E\clw_1(\eta_{0}^{N},\eta_{0})\to 0$ as $N\to\infty,$ and $X_{0}(N) \equiv (X_{0}^{1},X_{0}^{2},...,X_{0}^{N})$ is $\mu_0$ -chaotic then following
$$\sup_{n \ge 0} E\mathcal{W}_{1}\lt((\mu_{n}^{N},\eta_{n}^{N}),(\mu_{n},\eta_{n})\rt)\to 0$$
holds for $N\to \infty$.
\end{remark}
 One consequence of above theorem and Proposition \ref{Thmstable} will be the following interchange of limit results which is analogous to Corollary 3.5 from \cite{budhiraja2014long}.

\begin{Corollary}\label{cor2}
Under conditions of the Theorem \ref{unifintime}
\beqn
\limsup_{\min{\{N,M\}}\to \infty} \limsup_{n\to \infty}E\mathcal{W}_{1}((\bar{\mu}_{n}^{N},\bar{\eta}_{n}^{M}),(\mu_{\infty},\eta_{\infty}))&=& \limsup_{n\to \infty}\limsup_{\min{\{N,M\}}\to \infty}E\mathcal{W}_{1}((\bar{\mu}_{n}^{N},\bar{\eta}_{n}^{M}),(\mu_{\infty},\eta_{\infty})) \non\\&=& 0.\s\label{interchange}
\eeqn
\end{Corollary}

Suppose Assumptions of Theorem \ref{unifintime} hold and let $(\mu_{\infty},\eta_{\infty})$ be the fixed point of the map $\Psi$ of (\ref{Vlasovmseq}). We are interested in establishing a propagation of chaos result for $n = \infty.$  Recall for $\mathbb{IPS}_{2},$ $S^M(\bar{\eta}_{n}^{M})$ is the random measure defined as $\frac{1}{M}\sum_{i=1}^M \delta_{Y^{i,M}_{n}}$ where $\{Y^{i,M}_{n}\}_{i=1,\ldots,M}$ conditionally on $\mathcal{F}_{n}^{M,N},$ are $M$ i.i.d distributed $\R^{d}$ valued random variables according to $\bar{\eta}_{k-1}^{M}.$ Denote $Y_{n}(M):=(Y^{1,M}_{n},\ldots,Y^{M,M}_{n})$.

\begin{Theorem}\label{steadystate}
Consider the second particle system $\mathbb{IPS}_{2}$. Suppose Assumptions \ref{As2},\ref{As2.5},\ref{As3},\ref{As7} hold with conditions $$\delta \in (0,a_0),\s\s\s\s \sum_{i=0}^{\infty}(1-\alpha)^{i}\int_{\R^d}|y|P'P^{i}(0,dy)<\infty.$$ Then for every $N, M \ge 1,$ the Markov process $\big(\bar{X}^{N}(n), \bar{\eta}_{n}^{M},S^{M}(\bar{\eta}_{n}^{M})\big)_{n\ge 0}$ on $(\R^{d})^{N}\times \clp_{1}^{*}(\R^d)\times \clp_{}^{}(\R^d)$ has a unique invariant measure $\Theta_{\infty}^{N,M}$ if following holds
$$\max\Big\{\Big(\|A\|+\delta\sigma(2 +l_{PP'}^{\gr,\alpha}) + \alpha l(P')\Big ),(1-\alpha)l(P) \Big\}+\delta \sigma \max\big\{\alpha l_{P'}^{\gr},(1-\alpha) l_{P}^{\gr}\big\} < 1.$$ Let $\Theta_{\infty}^{1,N,M}$ be the marginal distribution on $(\R^{d})^{N}$ of the first co-ordinate of $\Theta_{\infty}^{N,M}$. Suppose additionally Assumption \ref{As3},\ref{As6} and Assumption \ref{unifas1},\ref{unifas2} hold with further condition for some $\tau>0$
$$ (1-\alpha)l^{}_{\tau}(P)<  1.$$
 Then $\Theta_\infty^{1,N,M}$ is $\mu_{\infty}$- chaotic, where $\mu_{\infty}$ is defined in Proposition \ref{Thmstable}.
\end{Theorem}
\begin{remark}
For first particle system $(\mathbb{IPS}_{1})$ similar steady state result holds for the discrete time Markov chain $\big(\bar{X}^{N}(n), \bar{\eta}_{n}^{N}\big)_{n\ge 0}$ on $(\R^{d})^{N}\times \clp_{1}^{*}(\R^d).$
\end{remark}

\subsection{Concentration Bounds:}

In order to obtain uniform in time concentration bounds of $\mathcal{W}_{1}\big((\mu_{n}^{N},\eta_{n}^{N}),(\mu_{n},\eta_{n})\big)$ we proceed according to those in Theorem 3.7 and Theorem 3.8 of \cite{budhiraja2014long} respectively. Here we establish two different types of concentration bounds. The first one is with initial non iid (i.e initial samples are $\mu_{0}$ chaotic) assumption and the second one is without that. 
\bas\label{exponas1}
\begin{description}
		\item{(i)} For some $K \in (1, \infty)$, $A_{1}(x) \le K$ for $\theta$ a.e. $x \in \R^m$. 
	\item{(ii)} There exists $\alpha \in (0, \infty)$ such that 
 $\int e^{\alpha |x|} \mu_0(dx) < \infty$ and  there exists $\alpha(\delta)\in(0,\alpha)$ such that $$\int_{\R^{m}} e^{\alpha(\delta) \big(A_{2}(z) + \frac{|B(z)|}{\delta}\big)} \theta(dz) < \infty.$$

		\end{description}
\eas

\bas\label{exponas2}
 Suppose there exists functions $h_{1}(\cdot)$,$h_{2}(\cdot)$, $h'_{1}(\cdot)$,$h'_{2}(\cdot),h_{3}(\cdot),h'_{3}(\cdot)$ ( $h_{2},h'_{2},h_{3},h'_{3}$ are nondecreasing with $h_{2}(0)=0$, $h'_{2}(0)=0;$), and constants $l_{h_{1}}\in(0,1],\,\,  l_{h'_{1}}\in(0,\infty)$ such that $h_{1}(x),$ and $h'_{1}(x)$ are respectively $l_{h_{1}}$ and $l_{h'_{1}}$ Lipschitz.
There exists $\alpha\in(0,\infty)$ such that following hold for all $\alpha_{1}\in(0,\alpha)$
\beqn
\int e^{\alpha_{1}|y|}P(x,dy)\le e^{h_{2}(\alpha_{1})}\big(e^{\alpha_{1}|h_{1}(x)|}+e^{h_{3}(\alpha_{1})}\big),\s \int e^{\alpha_{1}|y|}P'(x,dy) \le  e^{h'_{2}(\alpha_{1})}(e^{\alpha_{1}|h'_{1}(x)|}+e^{h'_{3}(\alpha_{1})}).\label{as10e1}
\eeqn
\eas

\begin{remark}
\begin{enumerate}[(a)]
\item For Gaussian transtion kernel $P(x,dy)=\frac{1}{\sqrt{2\pi}\lambda}e^{-\frac{(x-y)^2}{2\lambda^{2}}}dy,$ one has $$\int e^{\alpha_{1}|y|}P(x,dy)= e^{\frac{\alpha_{1}^{2}\lambda^{2}}{2}}\Big[e^{-\alpha x}\Phi\big(\frac{x}{\lambda}-\alpha\lambda\big)+e^{\alpha x}\Phi\big(\alpha\lambda+\frac{x}{\lambda}\big)\Big],$$
where $\Phi(\cdot)$ is the cumulative distribution function of Normal distribution. 
So \eqref{as10e1} holds with $h_{1}(x)=x,\s h_{3}(\cdot)=0,\s h_{2}(\alpha_{1})=\frac{\lambda^{2}\alpha_{1}^{2}}{2}.$
\item For Bi-exponential kernel $P(x,dy)=\frac{1}{2\lambda}e^{-\frac{|x-y|}{\lambda}}dy$ one has $$\int e^{\alpha_{1}|y|}P(x,dy)= e^{\alpha_{1} x}\Big[\frac{1}{1-\alpha_{1}^{2}\lambda^{2}}\Big].$$
So \eqref{as10e1} holds under condition $\alpha_{1}<\frac{1}{\lambda_{1}}$ with $h_{1}(x)=x,\s h_{3}(\cdot)=0,\s h_{2}(\alpha_{1})=\log\Big[\frac{1}{1-\alpha_{1}^{2}\lambda^{2}}\Big].$
Note that any kernel with tail lighter than exponential (like Gaussian) will satisfy \eqref{as10e1} for all $\alpha_{1},$ where for kernels with exponential like tail will have a specific restriction on $\alpha_{1}.$
\item We worked here only for $l_{h_{1}}=1$ as the upper bound. It only influences in the choice of $\alpha_{1}$ for which 
\beqn
\sup_{n\ge 0}\sup_{M,N \ge 1}E\lt<e^{\alpha_{1}|x|},\bar{\eta}^{M}_{n}\rt><\infty.\label{expofinite}
\eeqn
 For $l_{h_{1}}=1$ one has a definite upper bound of $\alpha_{1}.$ More precisely denoting  $\alpha_{1}h_{1}(0)\sum_{j=0}^{i}l^{j}_{h_{1}}+ \sum_{j=0}^{i}h_{2}(\alpha_{1}l^{j}_{h_{1}})$ by $g(i)$  if $g(i)$ is linear in $i$ (happens only for $l_{h_{1}}=1$) then there exists $\alpha^*$ such that \eqref{expofinite} holds for $\alpha_{1}<\alpha^{*}$. On the other hand if $g(\cdot)$ is bounded, then $\sup_{n\ge 0}\sup_{M,N \ge 1}E\lt<e^{\alpha_{1}|x|},\bar{\eta}^{M}_{n}\rt>$ will remain finite for all $\alpha_1>0$. If $g(i)$ is exponential in $i$ (when $l_{h_{1}}>1$) then the upper bound of $\sup_{n\ge 0}\sup_{M,N \ge 1}E\lt<e^{\alpha_{1}|x|},\bar{\eta}^{M}_{n}\rt>$ will diverge.
\end{enumerate}
\end{remark}

With $\tau, \sigma_1(\tau)$ defined above in Assumption \ref{unifas1} let 
\beqn
a(\tau) := \frac{4^{-\tau}-\|A\|^{1+\tau}}{\sigma_{1}(\tau)\big[1+(1+l^{\gr,\alpha}_{PP'})^{1+\tau}\big]}.
\eeqn

\begin{Theorem}\label{conv}\begin{enumerate}[(a)]
\item (Polynomial Concentration)\label{polyconv} 
Let $N_{1}=\min\{M,N\}.$ Suppose Assumptions (\ref{As2}-\ref{As7}) and Assumptions (\ref{unifas1}),(\ref{unifas2}) hold for some $\tau>0$. Suppose that
 $\delta \in(0,a(\tau)^{\frac{1}{1+\tau}}), (1-\alpha)l_{\tau}(P)<1$ and  
\beqn
\max\Big\{\Big(\|A\|+\delta\sigma(2 +l_{PP'}^{\gr,\alpha}) + \alpha l(P')\Big ),(1-\alpha)l(P) \Big\}+\delta \sigma \max\big\{\alpha l_{P'}^{\gr},(1-\alpha) l_{P}^{\gr}\big\}< 1.\s\text{}\s\label{polycond}
\eeqn
 Then there exits $\nu >1,\gamma\in(0,1)$, $N_0 \in \mathbb{N}_0$ and $C_1 \in (0, \infty)$ such that for all $\eps > 0,$ and for all $n\ge 0,$
$$ P(\mathcal{W}_{1}((\mu_{n}^{N},\eta_{n}^{M}),(\mu_{n},\eta_{n})>\varepsilon) \leq P(\mathcal{W}_{1}((\mu_{0}^{N},\eta_{0}^{M}),(\mu_{0},\eta_{0}))>\gamma\nu^n\varepsilon) + C_1 \varepsilon^{-(1+\alpha)} N_{1}^{-\frac{\tau}{d+2}}, $$
for all  $N_{1} >N_{0}
\left(\max\left\{1,\log^+ \varepsilon\right\}\right)^{\frac{d+2}{d}}$.
\item \label{expoconv}(Exponential Concentration)Let $N_{1}=\min\{M,N\}.$ Suppose that Assumptions  \ref{exponas1}  and \ref{exponas2} hold with \eqref{polycond}.  Suppose
	$\delta \in \Big[0,\frac{1- \|A\|}{(2+l^{\gr,\alpha}_{PP'})K}\Big)$ and $\alpha_1 \in \lt[0, \min\{\alpha^{*},\frac{\alpha(\delta)}{\delta}\}\rt)$ where
 $$\alpha^{*}|h_{1}(0)|+h_{2}(\alpha^*)=-\log(1-\alpha).$$ Then there exists $N_0 \in \mathbb{N}, \nu>1, \gamma\in(0,1)$ and $C_2\in (0, \infty)$ such that for all
	$\e > 0$ 
	$$P[\mathcal{W}_{1}((\mu_{n}^{N},\eta_{n}^{M}),(\mu_{n},\eta_{n})>\varepsilon]\leq P[\mathcal{W}_{1}((\mu_{0}^{N},\eta_{0}^{M}),(\mu_{0},\eta_{0}))>\gamma\nu^n\varepsilon]
	+ e^{-C_1 \e N_{1}^{1/d+2}},$$
	for all  $n \ge 0$, $N_{1} \ge N_0 \max\{ (\frac{1}{\e} \log^+ \frac{1}{\e})^{d+2}, \e^{(d+2)/(d-1)}\}$, if $d > 1$; and
	$$P[\mathcal{W}_{1}((\mu_{n}^{N},\eta_{n}^{M}),(\mu_{n},\eta_{n})>\varepsilon]\leq P[\mathcal{W}_{1}((\mu_{0}^{N},\eta_{0}^{M}),(\mu_{0},\eta_{0}))>\gamma\nu^n\varepsilon]
	+ e^{-C_1 (\e \wedge 1) N_{1}^{1/d+2}},$$
	for all $n\ge 0$, $N_{1} \ge N_0 \max\{ (\frac{1}{\e} \log^+ \frac{1}{\e})^{d+2}, 1\}$, if $d = 1$.  
\end{enumerate}
\end{Theorem}

\begin{remark}
\begin{enumerate}[(a)]
\item Similar concentration bounds hold for the first particle system $\mathbb{IPS}_{1}.$
\item Here the nonlinearity in the kernel of the nonlinear Markov process has a linear structure (linear combination of $P$ and $\mu P'$) which is handled through $\mathcal{W}_1$ distance. It can be further generalized for any nonlinear Markov process where the nonlinearity in the kernel depends on the higher order moments (of $p$th order) of the law of the chain, then working with $\mathcal{W}_{p}$ distance would yield similar results.
\end{enumerate}
\end{remark}
Note that the bounds in Theorems \ref{conv} are not dimensions independent while the initial sampling assumptions are not restrictive. It will be interesting to see if one can get sharper bounds under stronger conditions than above theorems. The following result shows that such bounds can be obtained 
in cases where initial locations of $N$ particles are i.i.d and under a more stringent condition on other parameters.

\begin{Theorem}\label{conviid}
Consider the first particle system $\mathbb{IPS}_{1}$ with initial condition $\eta_{0}^{N}\equiv\eta_{0}$. Suppose that $\{X_{0}^{i,N}\}_{i=1,\ldots,N}$ are i.i.d. with common distribution  $\mu_{0}$ for each $N$. Let 
\beqn
C_{1}&:=& \delta K\max\{1,(1-\alpha)l^{\gr}_{P}\alpha l(P')\}\frac{ \max\{\|A\|+\delta K(1+l^{\gr,\alpha}_{PP'}),\alpha l^{\gr}_{P'},(1-\alpha)l(P)\}}{\big|\|A\|+\delta K(1+l^{\gr,\alpha}_{PP'}) - \max\{\alpha l^{\gr}_{P'},(1-\alpha)l(P)\} \big|},\,\,\;\;\;\;\;\;\label{cou1}\\
 \chi_1&:=& \delta K \max\{\|A\|+\delta K(1+l^{\gr,\alpha}_{PP'}),\alpha l^{\gr}_{P'},(1-\alpha)l(P)\}+C_{1}.\label{cou2}
\eeqn
 Suppose that Assumptions  \ref{As2},\ref{As3},\ref{As7} and \ref{exponas1} hold with conditions $\chi_{1}\in (0,1)$, $\delta \in \Big[0,\frac{1- \|A\|}{(2+l^{\gr,\frac{\alpha}{\delta}}_{PP'})K}\Big)$ and $\alpha_{1}<\frac{\alpha(\delta)}{\delta}$. Then 
there exist $a_1, a_2,a'_1, a'_2,a''_1, a''_2 \in (0, \infty)$ and  $N_{0},N_{1},N_{2}$
for all $\e > 0$ 
\begin{eqnarray}\label{totalconc}
\sup_{n\ge 0}P[\mathcal{W}_{1}(\mu_{n}^{N},\mu_{n})>\varepsilon] \le 
\begin{cases}
 a_{1}e^{-N a_2 (\e^2 \wedge \e)}1_{\{d=1\}} \s\s\s\s\s\s\, N\ge N_{1}\max\{\frac{1}{\e},\frac{1}{\e^{2}}\},\\
a'_{1}e^{-N a'_2 \big(\big(\frac{\varepsilon}{\log(2+\frac{1}{\varepsilon})}\big)^2 \wedge \e\big)}1_{\{d=2\}}\s\s\; N\ge N_{2}\max\{\frac{1}{\e},\Big(\frac{\log(2+\frac{1}{\varepsilon})}{\e}\Big)^{2}\},\\
a''_{1}e^{-N a''_2 (\e^d \wedge \e)} 1_{\{d>2\}} \s\s\s\s\s\s N\ge N_{3}\max\{\frac{1}{\e},\frac{1}{\e^{d}}\}.
\end{cases}
\end{eqnarray}

\end{Theorem}
\begin{remark}
\begin{enumerate}[(a)]
	\item  If Assumption \ref{exponas1} is strengthened to $\int e^{\alpha(\delta)\lt(A_1^{2}(z) + \frac{|B(z)|)^2}{\delta^{2}}\rt)} \theta(dz) < \infty$ for some $\alpha(\delta) > 0$ then one can 
	strengthen the conclusion of Theorem \ref{conviid} as follows:
	For $\delta,\alpha$ sufficiently small
	there exist $N_0, a_1, a_2 \in (0, \infty)$ and a nonincreasing function $\vs_2: (0, \infty) \to (0,\infty)$ such that $\vs_2(t) \downarrow 0$ as $t \uparrow \infty$ and
	for all $\e > 0$ and $N \ge N_0 \vs_2(\e)$
	\[\sup_{n\ge 0}P[\mathcal{W}_{1}(\mu_{n}^{N},\mu_{n})>\varepsilon] \le 
	    a_1e^{-N a_2 \e^2 }. 
	\] 

\item Here stability condition \eqref{polycond} which is a crucial assumption for Lemma \ref{lem1} is not  used. Such is the power of the coupling that we used in Theorem \ref{conviid}.
\end{enumerate}
\end{remark}

\section{Discussion and Conclusion}

This article decribes a modified version of discrete time particle approximation scheme described in \cite{budhiraja2011discrete} which incorporates the evolution of particles in a non-compact domain. A similar form of stability condition is obtained under which the nonlinear system has a unique fixed point. Our contribution is computing the quantitative nonasymptotic bounds on these approximation schemes and how these relate to the conditions on the tail and smothness of the transition kernels $P,P'$ that were used to model the diffussive environment. As an additional result we obtained the propagation of chaos result of the particle scheme at time $n=\infty.$ There are few questions and remarks that should be addressed in future.
\begin{enumerate}[(a)]
\item Theorem \ref{conviid} is developed exclusvely for $\mathbb{IPS}_{1}$. For $\mathbb{IPS}_{2}$ we would have an extra term $\mathcal{W}_{1}\lt(S^{M}(\bar{\eta}^{M}_{n-1}),\bar{\eta}_{n-1}^{M}\rt)$ in the expression of $\mathcal{W}_{1}(\mu_{n}^{N},\mu_{n})$. Now the problem will arise in computing sharper (than \eqref{lempoly2})  bound of $$P[\mathcal{W}_{1}\lt(S^{M}(\bar{\eta}^{M}_{n-1}),\bar{\eta}_{n-1}^{M}\rt)>\varepsilon]=EP\big[\mathcal{W}_{1}\lt(S^{M}(\bar{\eta}^{M}_{n-1}),\bar{\eta}_{n-1}^{M}\rt)>\varepsilon \big|\bar{\mathcal{F}}_{n-1}^{M,N}\big].$$
Concentration bound of the conditional probability can be given in terms of random $\big<e^{\alpha_{1}|x|},\bar{\eta}_{n-1}^{M}\big>$ but getting an explicit relationship of the bound with the conditional exponential moment is unavailable. After taking expectation it is impossible conclude whether the inequality of upper bound still holds or not. Illustratively if the conditional concentration bound of $P\big[\mathcal{W}_{1}\lt(S^{M}(\bar{\eta}^{M}_{n-1}),\bar{\eta}_{n-1}^{M}\rt)>\varepsilon \big|\bar{\mathcal{F}}_{n-1}^{M,N}\big]$ is a concave function of $\big<e^{\alpha_{1}|x|},\bar{\eta}_{n-1}^{M}\big>$ then by  Jensen's inequality reasonable conclusion would hold but to our knowledge such explicit relationship is not present in literature.

\item The concentration bounds established in \cite{fournier2013rate} for $\mathcal{W}_{1}$ distance of empirical distribution of i.i.d observations to the true distribution is sharp however their method can be applied here only for $\mathbb{IPS}_{1}$ as done in Theorem \ref{conviid} using the well known coupling construction that works for all Vlasov McKean type systems. Without using that coupling, we attempted to use the grid based methods of \cite{fournier2013rate} in order to find sharper bounds for  $P[\mathcal{W}_1\big((\bar{\mu}_{n}^{N},\bar{\eta}_{n}^{M}),\Psi(\bar{\mu}_{n-1}^{N},\bar{\eta}_{n-1}^{M})\big) >\varepsilon]$ along the line of Theorem \ref{conv}. We faced similar problem as in the previous remark. Since one can derive a bound for $P[\mathcal{W}_1\big((\bar{\mu}_{n}^{N},\bar{\eta}_{n}^{M}),\Psi(\bar{\mu}_{n-1}^{N},\bar{\eta}_{n-1}^{M})\big) >\varepsilon\big|\bar{\mathcal{F}}_{n-1}^{M,N}]$ keeping  $\big<e^{\alpha_{1}|x|},\bar{\eta}_{n-1}^{M}\big>,\big<e^{\alpha_{1}|x|},\bar{\mu}_{n-1}^{N}\big>$ as constants but we do not know explicit structure how these bounds are functionally depending on $\big<e^{\alpha_{1}|x|},\bar{\eta}_{n-1}^{M}\big>,\big<e^{\alpha_{1}|x|},\bar{\mu}_{n-1}^{N}\big>,$ so that unconditionally we can conclude something useful. These issues will be addressed in future.
\end{enumerate}

\section{Proofs }\label{sec4}

The following two elementary lemmas give a basic moment bound that will be used in the proofs. We denote the function $f(\cdot,\cdot,\cdot,x)+\frac{B(x)}{\delta}$ by $f_{\delta}(\cdot,\cdot,\cdot,x).$
\begin{lemma}\label{l0}
For an interacting particle system illustrated in (\ref{nlps2}) and (\ref{evoconcpart}),
\begin{enumerate} [(a)]
\item Suppose Assumptions \ref{As2}, \ref{As2.5} and \ref{As3} hold. Then, for every $n\ge 1,\s M_{n}=\sup_{N\ge 1}E|X_{n}^{i}| <\infty.$
Moreover if Assumption \ref{As2} holds, then under $\delta\in (0,a_0)$ then $\s \sup_{n\ge 1}M_n <\infty.$ 
\item With the assumptions in part(a) suppose additionally Assumption \ref{unifas1} holds for some $\tau>0$ and suppose  $\delta\in (0,a(\tau)^{\frac{1}{1+\tau}}).$ Then
$$\sup_{N\ge 1}\sup_{n\ge 1}E|X_{n}^{i}|^{1+\tau}<\infty,$$
where in limit $a(\tau)^{\frac{1}{1+\tau}}\to a_0$ as $\tau \to 0^+$.
\end{enumerate}
\end{lemma}

\begin{remark}\label{rem4.1}
 Note that the same bound for $\sup_{n}\sup_{N,M \ge 1}E|\bar{X}_{n+1}^{i}|$ and $\sup_{n}\sup_{N,M \ge 1}E|\bar{X}_{n+1}^{i}|^{1+\tau}$  also hold for $\mathbb{IPS}_{2}$ under same condition on $\delta$. 
\end{remark}

\subsubsection{Proof of Lemma \ref{l0}}
\begin{enumerate} [(a)]
\item We prove the second statement. Proof of the first statement is similar. For each $n\ge 1$ and $i=1,\ldots,N,$ applying Assumption \ref{As2} on  particle system in (\ref{nlps2}) with definitions of $A_{1}(\cdot)$ and $A_{2}(\cdot)$
\beqn\label{l1eq1}
|X_{n+1}^{i}| \le \|A\| |X_{n}^{i}| + \delta A_1(\epsilon_{n+1}^{i})[|\gr \eta_{n}^{N}(X_{n}^{i})|+\|\mu_{n}^{N}\|_{1}+| X_{n}^{i}|]+\delta A_2(\epsilon_{n+1}^{i})+|B(\epsilon_{n+1}^{i})|.
\eeqn
Now by Assumption \ref{As3} using DCT one has 
\beqn
\gr \eta_{n+1}(y) =\int_{\R^{d}} \eta_{n}(x) [\gr_{y} R^{\alpha}_{\mu_{n}}(x,y)]dx\label{gr}
\eeqn
 for every $y$  since from Assumption \ref{As3} $\sup_{x\in\R^{d}}|\gr_{y} R^{\alpha}_{\mu_{n}}(x,y)|\le l^{\gr,\alpha}_{PP'}\, |y|+\sup_{x\in \R^{d}}\big((1-\alpha)|\gr_{y}P(x,0)|+\alpha|\gr_{y}P'(x,0)|\big).$ Applying the same condition followed by the inequality $|\gr \eta_{n+1}(y)| \le$\\  $\int_{\R^{d}} \eta_{n}(x) |\gr_{y} R^{\alpha}_{\mu_{n}}(x,y)|dx,$ one has
\beqn\label{l1eq1.5}
|\gr \eta_{n}(y)|\le  l^{\gr,\alpha}_{PP'}\, |y|+c^{\alpha}_{PP'}.
\eeqn
Also note by exchangeability $E\|\mu_{n}^{N}\|_{1}=E\int|x|\mu_{n}^{N}(dx)=E|X_{n}^{i}|$. Taking expectation in (\ref{l1eq1}) and using (\ref{l1eq1.5}) and independence between $\epsilon_{n+1}^{i}$ and $\{X_{n}^{j}\}_{j=1}^{N},$ one has
\beqn\label{l1eq2}
E|X_{n+1}^{i}|\le \bigg(\|A\|+\delta \sigma\lt(2+  l^{\gr,\alpha}_{PP'}\rt) \bigg)E|X_{n}^{i}|+  \delta [\sigma c^{\gr,\alpha}_{PP'}+\sigma_2(\delta)].
\eeqn
The assumption on $\delta$ implies that $\gamma:=\|A\|+\delta \sigma\lt(2+  l^{\gr,\alpha}_{PP'}\rt) \in(0,1).$ A recursion on (\ref{l1eq2}) will give
$M_n \le \gamma^{n}E|X_{0}^i|+\frac{\delta [\sigma c^{\gr,\alpha}_{PP'}+\sigma_2]}{1-\gamma},$ from which the result follows.

\item By Holder's inequality for any three nonnegative real numbers $a_1,a_2,a_3,a_{4}$
\beqn\label{ho}
(a_1+a_2+a_3+a_{4})^{1+\tau}\le 4^{\tau}(a_{1}^{1+\tau}+a^{1+\tau}_{2}+a^{1+\tau}_{3}+a^{1+\tau}_{4}).
\eeqn

Starting with (\ref{l1eq1}), applying (\ref{ho}), and Assumption \ref{As2}, on (\ref{l1eq1}) we have
\beqn
|X_{n+1}^{i}|^{1+\tau} &\le& 4^{\tau}\bigg[\|A\|^{(1+\tau)} |X_{n}^{i}|^{1+\tau}  + \lt(\delta A_1(\epsilon_{n+1}^{i})[1+l_{PP'}^{\gr,\alpha}] |X_{n}^{i}|\rt)^{1+\tau} +\lt(\delta A_1(\epsilon_{n+1}^{i})\|\mu_{n}^{N}\|_{1}\rt)^{1+\tau}\non\\&&\s+ \delta^{1+\tau}\big[A_{1}(\epsilon_{n+1}^{i}).c^{\alpha}_{PP'}+A_{2}(\epsilon_{n+1}^{i})+\frac{|B(\epsilon_{n+1}^{i})|}{\delta}\big]^{1+\tau}\bigg].\non 
\eeqn
For any convex function $\phi(\cdot),$ applying Jensen's inequality one gets $\phi(\|\mu_{n}^{N}\|_{1})\le \int |\phi(x)| \mu_{n}^{N}(dx)$ $=\frac{1}{N}\sum_{i=1}^{N}|\phi(X_{n}^{i})|.$ Using $\phi(x)=x^{1+\tau},$ after taking expectation one gets following recursive equation for $E|X_{n+1}^{i}|^{1+\tau}$,
\beqn
E|X_{n+1}^{i}|^{1+\tau}\le 4^{\tau}\bigg[\|A\|^{(1+\tau)} +\delta^{1+\tau}\sigma_{1}(\tau)\big[(1+l^{\gr,\alpha}_{PP'})^{1+\tau}+1\big]\bigg]E|X_{n}^{i}|^{1+\tau}+ \delta^{1+\tau}8^{\tau}\bigg[\sigma_{1}(\tau)c^{\tau}_{PP'}+\sigma_{2}(\delta,\tau)\bigg].\nn
\eeqn
Note that for our condition on $\delta,\s\s \kappa_1 := 4^{\tau}\bigg[\|A\|^{(1+\tau)} +\delta^{1+\tau}\sigma_{1}(\tau)\big[(1+l^{\gr,\alpha}_{PP'})^{1+\tau}+1\big]\bigg]<1.$  Thus 
\beqn\label{finmom} 
\sup_{n\ge 1} E|X_{n}^{i}|^{1+\tau}\le \kappa^{n}_{1} E|X_{0}^{i}|^{1+\tau}+\frac{\delta^{1+\tau}8^{\tau}\big[\sigma_{1}(\tau)c^{\tau}_{PP'}+\sigma_{2}(\delta,\tau)\big]}{1-\kappa_1}.
\eeqn
\end{enumerate}
$\square$

\begin{lemma}\label{l1}
Suppose Assumptions \ref{As2},\ref{As2.5},\ref{As3} and \ref{As7} hold.
\begin{enumerate} [(a)]
\item Consider the interacting particle system described in (\ref{nlps2}) and (\ref{evoconcpart}). Then, for every $n\ge 1,$
\beqn
\left<|x|,\eta_{n}\right><\infty,\s\s \sup_{N\ge 1}E\left<|x|,\eta^{N}_{n}\right><\infty.
\eeqn
Moreover if Assumption \ref{As2} holds, then under conditions 
\beqn\label{condition0}
\delta  \in \left(0,a_{0}\right),\s\s\text{and}\s\s\s\sum_{i=0}^{\infty}(1-\alpha)^{i}\int_{\R^d}|y|P'P^{i}(0,dy)<\infty,
\eeqn
one has $\sup_{n\ge 1}\left<|x|,\eta_{n}\right><\infty.$

Additionally assuming $\sup_{N\ge 1}E\lt<|x|,\eta_{0}^{N}\rt><\infty$ one gets
$$\s\s\sup_{n\ge 1} \sup_{N\ge 1}E\left<|x|,\eta^{N}_{n}\right><\infty.$$ 
\item With the assumptions in part(a) suppose additionally Assumption \ref{unifas1},\ref{unifas2} hold for some $\tau>0$ and suppose  $\delta\in (0,a(\tau)^{\frac{1}{1+\tau}}).$ Then with condition $(1-\alpha) m_{\tau}(P)<  1$ one has
$\sup_{n\ge 1}\left<|x|^{1+\tau},\eta_{n}\right><\infty.$  Additionally assuming $\sup_{N\ge 1}E\lt<|x|^{1+\tau},\eta_{0}^{N}\rt><\infty$ one gets $\sup_{n\ge 1} \sup_{N\ge 1}E\left<|x|^{1+\tau},\eta^{N}_{n}\right><\infty,$ where in limit $a(\tau)^{\frac{1}{1+\tau}}\to a_0$ as $\tau \to 0^+$.

\end{enumerate}
\end{lemma}
\begin{remark}
The second condition in (\ref{condition0}) is very general. It doesn't impose any condition on $\alpha \in(0,1).$ The condition holds for all transition kernels $P(x,\cdot), P'(x,\cdot)$ with finite first moment. Only thing one needs to check $$\int_{\R^d}|y|P'P^{i}(0,dy) =g(i)$$ where $g(i)$ is some polynomial in $i$ (For Gaussian it's linear). If $g(\cdot)$ is an exponential function then it will impose a further lower bound condition on $\alpha$.
\end{remark}
\begin{Corollary}\label{cor}
For $\mathbb{IPS}_{2}$ same conclusion about $\bar{\eta}^{M}_{n}$ holds as $\eta_{n}^{N}$ in first particle system specified in Lemma \ref{l1} under same set of conditions on $\delta,\alpha$. Note that $\bar{\eta}_{0}^{M}=\eta_{0},$ so we don't need to assume anything about the initial sampling scheme like $\sup_{M\ge 1}E\lt<|x|,\bar{\eta}_{0}^{M}\rt><\infty$ (or $\sup_{M\ge 1}E\lt<|x|^{1+\tau},\bar{\eta}_{0}^{M}\rt><\infty$) since they automatically hold for $\eta_{0}\in \clp_1^{*}(\R^d)$ (or $\eta_{0}\in \clp_{1+\tau}^{*}(\R^d)$) respectively.
\end{Corollary}

\subsubsection{Proof of Lemma \ref{l1}}
We will start with the second part of part (a) of the lemma. First part will follow similarly. We will show if $\eta_{0}\in \clp^{*}_1(\R^{d})$ then $\eta_n \in \clp_1(\R^{d})$ for all $n\ge 1.$ Note that 
\beqn\label{etal}
\eta_{k+1} = \sum_{i=0}^{k}\left[\alpha(1-\alpha)^{i}\mu_{k-i}P'P^{i}\right]+(1-\alpha)^{k+1}\eta_{0}P^{k+1}.
\eeqn
From Assumption \ref{As7}, it is obvious that $P'P^{i}f$ is $l(P')l(P)^{i}$ Lipschitz if $f$ is a $1$-Lipschitz function. It implies $|P'P^{i}f(x) - P'P^{i}f(0)|\le l(P')l(P)^{i}|x|$ for any $f\in \Lip_1(\R^d).$ Since $|x|$ is $1$-Lipschitz, one has $$P'P^{i}|x|\le l(P')l(P)^{i}|x|+\int_{\R^{d}}|y|P'P^{i}(0,dy).$$ Using this inequality one has from (\ref{etal})
\beqn
&&\left<|x|,\eta_{k+1}\right> = \sum_{i=0}^{k}[\alpha(1-\alpha)^{i}\left<|x|,\mu_{k-i}P'P^{i}\right>]+(1-\alpha)^{k+1}\left<|x|,\eta_{0}P^{k+1}\right>\nonumber\\
&\le&\sum_{i=0}^{k}[\alpha(1-\alpha)^{i}\left<l(P')l(P)^{i}|x|,\mu_{k-i}\right>]+\alpha \sum_{i=0}^{\infty}(1-\alpha)^{i}\int_{\R^d}|y|P'P^{i}(0,dy)+[(1-\alpha)l(P)]^{k+1}\left<|x|,\eta_{0}\right>\nonumber\\
&\le&\alpha l(P')\left\{\sup_{n\in \mathbb{N}}\left<|x|,\mu_{n} \right>\right\}\sum_{i=0}^{k}\left[(1-\alpha)l(P)\right]^{i} + \alpha \sum_{i=0}^{\infty}(1-\alpha)^{i}\int_{\R^d}|y|P'P^{i}(0,dy)\nonumber\\&&+[(1-\alpha)l(P)]^{k+1}\left<|x|,\eta_{0}\right>.\label{lem4.3e1}
\eeqn
By Assumption \ref{As7}, $l(P)\le 1,$ implies $(1-\alpha)l(P)<1$. From similar derivation done in Lemma \ref{l0}, one has $\sup_{n\in \mathbb{N}}\left<|x|,\mu_{n} \right> <\infty$  if $\delta \in (0,a_0).$ The result follows using all the conditions
$$\sup_{k\in\mathbb{N}}\left<|x|,\eta_{k}\right><\infty.$$
For  $E\left<|x|,\eta^{N}_{k}\right>$ note that for any function $f,$
\beqn\label{etal2}
\lt<f,\eta^{N}_{k+1}\rt> = \sum_{i=0}^{k}\left[\alpha(1-\alpha)^{i}\lt<f,\mu^{N}_{k-i}P'P^{i}\rt>\right]+(1-\alpha)^{k+1}\lt<f,\eta^{N}_{0}P^{k+1}\rt>.
\eeqn

From Lemma \ref{l0} $\sup_{n \ge 0}\sup_{N\ge 1}E\left<|x|,\mu^{N}_{n} \right> <\infty$ for $\delta \in (0,a_0).$
Putting $f(x)=|x|,$ then expanding $\left<|x|,\eta^{N}_{n} \right>$ similarly like (\ref{lem4.3e1}) after taking expectation one gets a similar bound and finiteness of $\sup_{n}\sup_{N\ge 1}E\left<|x|,\eta^{N}_{n} \right>$ follows from that.

$\square$

\textit{Proof of Lemma \ref{l1}(b):}
 From (\ref{etal}), 
\beqn\label{etaltau}
\left<\eta_{k+1},|x|^{1+\tau}\right> = \sum_{i=0}^{k}\left[\alpha(1-\alpha)^{i}\left<\mu_{k-i}P'P^{i},|x|^{1+\tau}\right>\right]+(1-\alpha)^{k+1}\left<\eta_{0}P^{k+1},|x|^{1+\tau}\right>.
\eeqn
From Assumption \ref{unifas2} we get the following recursion for $a_{i}:=\left<\mu P'P^{i},|x|^{1+\tau}\right>$ for any measure $\mu\in\clp_{1+\tau}(\R^d)$
\beqn
a_{i}=\left<\mu P'P^{i-1},P|x|^{1+\tau}\right>\le m_{\tau}(P)(1+a_{i-1})\label{cc1}
\eeqn
since $P|x|^{1+\tau}\le m_{\tau}(P)(1+|x|^{1+\tau})$ from Assumption \ref{unifas2}. Using the fact $a_{0}:=\lt<\mu,P'|x|^{1+\tau}\rt>\le m_{\tau}(P')(1+\lt<\mu,|x|^{1+\tau}\rt>),$ we finally have
\beqn
\left<\eta_{k+1},|x|^{1+\tau}\right>&\le&\alpha\sum_{i=0}^{k}(1-\alpha)^{i}\bigg[m_{\tau}(P)\frac{ l^{i}_{\tau}(P)-1}{m_{\tau}(P)-1}+ m_{\tau}(P') l^{i}_{\tau}(P)\big[1+\left<|x|^{1+\tau},\mu_{k-i}\rt>\big]\bigg]\non\\
&& +(1-\alpha)^{k+1}\bigg[m_{\tau}(P)\frac{ l^{k+1}_{\tau}(P)-1}{m_{\tau}(P)-1}+ l^{k+1}_{\tau}(P)\lt<\eta_{0},|x|^{1+\tau}\right>\bigg].
\eeqn
Under condition $\delta \in (0,a(\tau)^{\frac{1}{1+\tau}})$ and $(1-\alpha)m_{\tau}(P)< 1$ one gets $\sup_{n}\left<\eta_{n},|x|^{1+\tau}\right> <\infty.$ 
Similarly the same bound can be derived for $\sup_{n}\sup_{N\ge 1}E\lt<|x|^{1+\tau},\eta_{n}^{N}\rt>$ under the same set of conditions.

$\square$

\subsubsection{Proof of Corollary \ref{cor}}
To prove the Corollary about $\bar{\eta}_{n}^{M},$ define the random operator $S^{M}\circ P$ acting on the probability measure $\mu$ on $\R^d:\s \mu(S^{M}\circ P) = (S^{M}(\mu))P.$ Note the following recursive form of $\bar{\eta}_{n}^{M}$:
\beqn
\bar{\eta}^{M}_{k+1} = \sum_{i=0}^{k}\left[\alpha(1-\alpha)^{i}\bar{\mu}^{N}_{k-i}P'(S^{M}\circ P)^{i}\right]+(1-\alpha)^{k+1}\eta_{0}(S^{M}\circ P)^{k+1}.
\eeqn
Note that for any function $f$ one has 
\beqn
E\left<\mu(S^{M}\circ P),f\right> = E\left<S^{M}(\mu),Pf\right> =\left<\mu,Pf\right>=\left<\mu P,f\right>.\nonumber
\eeqn 
Now by expanding $\mu(S^{M}\circ P)^{k}$ one gets,
\beqn
\mu(S^{M}\circ P)^{k}= \left[\mu(S^{M}\circ P)^{k-1}\right](S^{M}\circ P)=S^{M}(\mu(S^{M}\circ P)^{k-1}) P.\nonumber
\eeqn
Taking expectation one has
\beqn
E\left<\mu(S^{M}\circ P)^{k},f\right>&=&E\left<S^{M}\left(\mu(S^{M}\circ P)^{k-1}\right)P,f\right>=E\left<S^{M}\left(\mu(S^{M}\circ P)^{k-1}\right), Pf\right>\nonumber\\
&=&E\left<\mu(S^{M}\circ P)^{k-1},Pf\right>=E\left<\mu(S^{M}\circ P)^{k-1}P, f\right>.\nonumber
\eeqn
Continuing this calculation $k-1$ times one has $E\left<\mu(S^{M}\circ P)^{k},f\right>=\left<\mu P^{k},f\right>$ which leads to the following expression
\beqn
E\left<\bar{\mu}^{N}_{k-i}P'(S^{M}\circ P)^{i},f\right>&=&EE\left[\left<\bar{\mu}^{N}_{k-i}P'(S^{M}\circ P)^{i},f\right>\bigg| \mathcal{F}^{M,N}_{k-i}\right]\nonumber\\
&=&E\left[\left<\bar{\mu}^{N}_{k-i}P'P^{i},f\right>\right]=E\left[\left<\bar{\mu}^{N}_{k-i},P'P^{i}f\right>\right].\label{coreq}
\eeqn
The corollary is proved by observing (\ref{coreq}). The same bound holds for both $E\left<\bar{\eta}^{M}_{n},f\right>$, $E\left<\eta^{N}_{n},f\right>$ because of the similarity of bounds of $E\lt<f,\mu_{n}^{N}\rt>$, and $E\lt<f,\bar{\mu}_{n}^{N}\rt>$ for $f(x)=|x|,|x|^{1+\tau}, e^{\alpha|x|^{p}}$ which follows from Remark \ref{rem4.1}.

$\square$

\subsection{Proof of Proposition \ref{pointwise}}
We will prove part (b) of the theorem. Part (a) will follow similarly. We will start with the following lemma.

\begin{lemma}\label{lemth1}
\begin{enumerate}[(a)]
\item Under Assumptions \ref{As2},\ref{As2.5},\ref{As3},  for every $\epsilon >0$ and $n\geq 1$, there exists a compact set $K_{\epsilon,n} \in \mathcal{B}(\mathbb{R}^{d})$ such that 
$$\sup_{M,N\geq 1} E\left\{ \int_{K_{\epsilon,n}^{c}} |x| \left(\mu_{n}^{N}(dx)+\mu_{n-1}^{N}Q^{\bar{\eta}_{n-1}^{M},\mu_{n-1}^{N}}(dx)\right)\right\}<\epsilon .$$
\item Suppose Assumptions \ref{As2},\ref{As2.5},\ref{As3},\ref{As7},\ref{As8unifintegr} hold. Then for every $\epsilon >0$ and $k\geq 1$, there exists a compact set $K_{\epsilon,k} \in \mathcal{B}(\mathbb{R}^{d})$ such that
$$\sup_{M,N\ge 1} E \left<|x|.1_{K_{k,\epsilon}}, S^{M}(\bar{\eta}_{k}^{M}) + \bar{\eta}_{k}^{M}\rt> <\eps.$$
This part of the lemma is exclusively for part (b) of the Proposition \ref{pointwise}.
\end{enumerate}
\end{lemma}
\textbf{Proof:} Note that for any non-negative $\phi: \mathbb{R}^{d} \to \mathbb{R}$,  
\begin{eqnarray}
E \int \phi(x) \mu_{n}^{N}(dx) &=& \frac{1}{N} \sum_{k=1}^{N} E \phi(X_{n}^{k}) = E\phi(X_{n}^{1}),\s\s\s\s\s\label{sum} \\ 
E\int \phi(x) \mu_{n-1}^{N}Q^{\bar{\eta}_{n-1}^{M},\mu_{n-1}^{N}}(dx) &=& \frac{1}{N} \sum_{i=1}^{N} E(E(\langle \phi , \delta_{X_{n}^{i}}Q^{\bar{\eta}_{n-1}^{M},\mu_{n-1}^{N}}\rangle \mid \mathcal{F}_{n}))
\non\\&=& \frac{1}{N}\sum_{i=1}^{N} E\phi\left(AX_{n}^{i} + \delta f_{\delta}(X_{n}^{i},\mu_{n}^{N},\gr \eta_{n}^{N}(X_{n}^{i}), \epsilon_{n+1}^{i})\right)\nonumber\\
&=&\frac{1}{N}\sum_{i=1}^{N} E\phi(X_{n+1}^{i})=  E\phi(X_{n+1}^{1}). \label{eq1.2.1}
\end{eqnarray}
To get the desired result from above equalities it suffices to show that 
\begin{eqnarray}
\text{the family $\{X_{n}^{i,N},i=1,...,N; M,N\geq 1\}$ is uniformly integrable for every $n\geq 0$.}\label{e1} 
\end{eqnarray}
We will prove (\ref{e1}) by induction on $n$. 
Once more we suppress $N$ from the super-script. 
Clearly by our assumptions $\{X_{0}^{i},i=1,...,N; N\geq 1\}$ is uniformly integrable. Now suppose that the Statement (\ref{e1}) holds for some $n$. Note that from (\ref{l1eq1}) and \eqref{l1eq1.5}
\beqn
|X_{n+1}^{i}| &\le& \|A\| |X_{n}^{i}| + \delta A_1(\epsilon_{n+1}^{i})[|\gr \eta_{n}^{N}(X_{n}^{i})|+\|\mu_{n}^{N}\|_{1}+| X_{n}^{i}|]+\delta A_2(\epsilon_{n+1}^{i})+|B(\epsilon_{n+1}^{i})|.\nonumber\\
&\le& \|A\| |X_{n}^{i}| + \delta A_1(\epsilon_{n+1}^{i})[\|\mu_{n}^{N}\|_{1}+(1+l_{PP'}^{\gr,\alpha})| X_{n}^{i}|]+\delta A_2(\epsilon_{n+1}^{i})+|B(\epsilon_{n+1}^{i})|+\delta c^{\alpha}_{PP'}A_{1}(\epsilon^{i}_{n+1})\non\\
&\le& \|A\| |X_{n}^{i}| + \delta A_1(\epsilon_{n+1}^{i})[\frac{1}{N}\sum_{i=1}^{N}| X_{n}^{i}|+(1+l_{PP'}^{\gr,\alpha})| X_{n}^{i}|]+\delta  A_2(\epsilon_{n+1}^{i})+|B(\epsilon_{n+1}^{i})|+\delta c^{\alpha}_{PP'}A_{1}(\epsilon^{i}_{n+1})\non
\eeqn
From Assumptions \ref{As2} and \ref{As2.5} the families $\{A_1(\epsilon_{n+1}^{i});i\geq 1\}$, $\{A_{2}(\epsilon_{n+1}^{i});i\geq 1\}$  $\{B_{2}(\epsilon_{n+1}^{i})$ are uniformly integrable. Now by exchangeability, $\frac{1}{N}\sum_{i=1}^{N}| X_{n}^{i}|=E\Big[| X_{n}^{i}| \Big| \sigma\Big(\frac{1}{N}\sum_{i=1}^{N}\delta_{X_{n}^{i}}\Big)\Big].$ If $\{X_{\alpha}:\alpha\in \Gamma_{1}\}$ is uniformly integrable, and $\{\sigma_{\beta},\beta\in\Gamma_{2}\}$ is a collection of $\sigma$- fields where $\Gamma_{1},\Gamma_{2}$ are arbitrary index sets, then $\{E(X_{\alpha}|\sigma_{\beta}),(\alpha,\beta)\in\Gamma_{1}\times\Gamma_{2}\}$ is also a uniformly integrable family. It follows that $\{\frac{1}{N}\sum_{i=1}^{N}|X_{N}^{i}|, N\ge 1\}$ is a uniformly integrable family from induction  hypothesis. Using \eqref{e1} again along with independence between $\{\epsilon_{n+1}^{i},i=1,\ldots,N\}$ and $\{X_{n}^{i}:i=1,\ldots,N;N\geq 1\}$ yield that the family $\{|X_{n+1}^{i}|:i=1,\ldots,N;N\geq 1\}$ is uniformly integrable. The result follows. $\s\square$

\textit{Proof of Lemma \ref{lemth1}(b):}
Note that $S^M(\bar{\eta}_{k}^{M})=\frac{1}{M}\sum_{i=1}^M \delta_{Y^{i,M}_{k}}$ where $\{Y^{i,M}_{k}\}_{i=1}^{M}\bigg| \mathcal{F}^{M,N}_{k}$
are i.i.d from $\bar{\eta}_{k}^{M}.$ So for any non-negative function $\phi$ we have
\beqn
E\langle\phi, S^M(\bar{\eta}_{k}^{M})\rangle &=& E\frac{1}{M}\sum_{i=1}^M \phi(Y^{i,M}_{k})=EE\left[\frac{1}{M}\sum_{i=1}^M \phi(Y^{i,M}_{k}) \big|\mathcal{F}^{M,N}_{k}\right]=E E\left[ \phi(Y^{i,M}_{k})\big|\mathcal{F}^{M,N}_{k}\right]\nonumber\\
&=&E\phi(Y^{i,M}_{k})= E\langle \phi,\bar{\eta}_{k}^{M}\rangle.
\eeqn
We will prove the result if we can show the family 
\beqn\label{lemunifind0}
\{Y^{i,M}_{k},i=1,\ldots,M; M,N\ge 1\}\s\text{ is uniformly integrable for every }\s k\ge 0.
\eeqn
We will prove (\ref{lemunifind0}) through induction on $k$. For $k=0,$ the result follows trivially since $\{Y^{i,M}_{0},i=1,\ldots,M; M\ge 1\}$ are i.i.d from $\eta_0.$ Suppose it holds for $k=n.$
We will show that both, 
\beqn\label{uniflemind}
&&\{S^M(\bar{\eta}_{n}^{M})P: M,N\ge1\}\s\text{ and }\s\{\bar{\mu}_{n}^{N}P':N\ge 1\}\s\text{are uniformly integrable families of}\s \non\\ &&\text{probability measures}. 
\eeqn
Then from the structure $\bar{\eta}_{n+1}^{M}=(1-\alpha)S^M(\bar{\eta}_{n}^{M})P+\alpha \bar{\mu}_{n}^{N}P,$ it is evident that  $\{\bar{\eta}_{n+1}^{M}:M,N\ge 1\}$ is uniform integrable which equivalently implies $\{Y_{n+1}^{i,M}:i=1,\ldots,M; M,N \ge 1\}$ is UI too.
On proving the first assertion in (\ref{uniflemind}), note that due to the exchangeability of $\{Y_{n}^{i,M}:i=1,\ldots,M\},$ one has 
\beqn
S^M(\bar{\eta}_{n}^{M})P=E\left[\delta_{Y_{n}^{1,M}}P\bigg| \sigma \left(\frac{1}{M}\sum_{i=1}^{M}\delta_{Y_{n}^{i,M}}\rt)\rt].\label{uniflem0}
\eeqn
  We know that if $\{Z_{\alpha}, \alpha \in \Gamma_1\}$ is a uniformly integrable family
and $\{\clh_{\beta}, \beta \in \Gamma_2\}$ is a collection of $\sigma$-fields where $\Gamma_1, \Gamma_2$ are arbitrary index sets, then
$\{E (Z_{\alpha} \mid \clh_{\beta}), (\alpha, \beta) \in \Gamma_1\times \Gamma_2\}$ is a uniformly integrable family. So from (\ref{uniflem0}) it suffices to prove that $\{\delta_{Y_{n}^{i,M}}P: i=1,\ldots,M; M,N\ge 1\}$ is uniformly integrable. Define a function $f_k(.)$ such that, $f_k(x)=0$, if $|x|\in[0,\frac{k}{2}]$ and $f_k(x)=|x|$, if $|x| \ge k$ and linear in between range. Then by construction $f_{k}(.)$ is Lipschitz with coefficient 2 and  $x.1_{\{|x|>k\}}\le f_{k}(x)$ for all $x\in\R^d.$ By Assumption \ref{As8unifintegr} we have that $\{P(z,.):z\in K\}$ is uniformly integrable. So taking the compact set $K=\{|x|\le k\}$ assuming $Y_{n}^{i,M}$ has unconditional law $m^{n}_{i}$ for all $i=1,\ldots,M,$ the quantity
\beqn
\int_{|z|>L}\int y.1_{\{K^{c}\}}P(z,dy)m_{i}^{n}(dz)&\le& \int_{|z|>L}\left[f_{k}(y)P(z,dy)\rt]m_{i}^{n}(dz)\non\\
&\le& \int_{|z|>L}\left[|Pf_{k}(0)|+2l(P)|z|\rt]m_{i}^{n}(dz)\s\s\s\label{lemunif}\\
&\le& Pf_{k}(0)\int_{|z|>L}m_{i}^{n}(dz)+ 2l(P)\int_{|z|>L}|z|m_{i}^{n}(dz).\label{lemunif2}
\eeqn
The display in (\ref{lemunif}) follows from Assumption \ref{As7} and using Lipschitz property of $f_{k}.$ After taking supremum in the set $\{i=1,\ldots,M; M,N \ge 1\}$ in both sides of  (\ref{lemunif2}), second part of R.H.S goes to $0,$ as $L\to \infty$ by induction hypothesis. About the first part $Pf_k(0)$ goes to $0$ as $k\to\infty$ by D.C.T since ($\int |y|P(0,dy)<\infty$) and also $\int_{|z|>L}m_{i}^{n}(dz)$ converges to $0$ (as $L$ goes to $\infty$) due to the tightness of $\{m_{i}^{n}: i=1,\ldots,M;M,N\ge 1\}$ which also follows from induction hypothesis. The second assertion that $\{\bar{\mu}_{n}^{N}P':N\ge 1\}$ is uniformly integrable follows similarly through induction.

$\square$

We will proceed to the main proof via induction on $n\in\mathbb{N}$ for the quantity $E\left[\clw_1(\bar{\mu}_{n}^{N},\mu_{n})+ \clw_1(\bar{\eta}_{n}^{N},\eta_{n})\right]$. For $n=0$, we will first show that $E\clw_1(\bar{\mu}_{0}^{N},\mu_{0})\to 0$ as $N\to \infty.$ From \cite{sznitman1991topics} we have $$(\bar{X}_{0}^{1},\bar{X}_{0}^{2},\ldots,\bar{X}_{0}^{N})\text{ is }\mu_0\text{-chaotic} \Leftrightarrow \bar{\mu}_{0}^N\text{ converges weakly to }\mu_0\text{ in probability }\Leftrightarrow\beta(\bar{\mu}_{0}^N,\mu_0)\overset{p}{\to} 0.$$

 From Lemma \ref{lemth1} one can construct $K_{0,\epsilon}$ compact ball containing $0,$ so that $E\left<|x|.1_{K^c_{0,\epsilon}},\bar{\mu}_{0}^{N}\right><\frac{\eps}{2}$ and $\left<|x|.1_{K^c_{0,\epsilon}},\mu_0\right> <\frac{\eps}{2}$ hold. So using the fact for any $f\in \Lip_1(\R^{d})$ with $f(0)=0,$ one has $|f(x)|\le|x|.$ 
\beqn
E\clw_1(\bar{\mu}_{0}^N,\mu_0)&=&E \sup_{f\in \Lip_1(\R^d)}|\left<f,\bar{\mu}_{0}^N -\mu_0\right>|=E \sup_{f\in \Lip_1(\R^d), f(0)=0}|\left<f,\bar{\mu}_{0}^N - \mu_0\right>|\nonumber\\
&\le& E \sup_{f\in \Lip_1(\R^d), f(0)=0}|\left<f 1_{K_{0,\epsilon}},\bar{\mu}_{0}^N - \mu_0\right>|+E\left<|x|1_{K^c_{0,\epsilon}},\bar{\mu}_{0}^{N}\right>+\left<|x|.1_{K^c_{0,\epsilon}},\mu_0\right>\nonumber\\
&\le& \text{diam}(K_{0,\eps})E \beta(\bar{\mu}_{0}^N,\mu_0) +\eps.
\eeqn
In last display we used the fact that $\sup_{x\in K_{0,\eps}}|f(x)|\le \text{diam}(K_{0,\eps})$. Note that $\beta(\bar{\mu}_0^N,\mu_0)$ is bounded by $2$ (so Uniformly Integrable) and $\beta(\bar{\mu}_{0}^N,\mu_0) \overset{p}{\to} 0$ implies $ E \beta(\mu_{0}^N,\mu_0) \to 0$ as $N\to \infty$ proving the assertion (\ref{ptth1}) for $n=0$. Suppose it holds for $n\le k.$
We start with the following triangular inequality
\begin{eqnarray}
\mathcal{W}_{1}(\bar{\mu}_{k+1}^{N},\mu_{k+1})  &\leq& \mathcal{W}_{1}(\bar{\mu}_{k+1}^{N},\bar{\mu}_{k}^{N} Q^{\bar{\eta}_{k}^{N},\bar{\mu}_{k}^{N}})  +\mathcal{W}_{1}(\bar{\mu}_{k}^{N} Q^{ \bar{\eta}_{k}^{N},\bar{\mu}_{k}^{N} }, \bar{\mu}_{k}^{N} Q^{ \eta_{k},\bar{\mu}_{k}^{N} })\non\\&& +\mathcal{W}_{1}(\bar{\mu}_{k}^{N} Q^{ \eta_{k},\bar{\mu}_{k}^{N} },\mu_{k+1}).\quad  \label{eq1.3.1}
\end{eqnarray}

Consider the third term of (\ref{eq1.3.1}). From the general calculations follwed by (\ref{us})-(\ref{A1}), we  have the following estimate,
\beqn\label{eq1.3.5.1}
\mathcal{W}_{1}(\bar{\mu}_{k}^{N} Q^{ \eta_{k},\bar{\mu}_{k}^{N} },\mu_{k}Q^{\eta_{k},\mu_{k}})\le \lt(\|A\|+\delta \sigma(2+ l_{PP'}^{\gr,\alpha})\rt)\mathcal{W}_{1}(\bar{\mu}_{k}^{N},\mu_{k}).
\eeqn 

Now we consider the first term of the right hand side of (\ref{eq1.3.1}). We will use Lemma \ref{lemth1}(a). Fix $\epsilon>0$ and let $K_{\epsilon}$ be a compact set in $\mathbb{R}^{d}$ such that
$$\sup_{N\geq 1} E\left\{ \int_{K_{\epsilon}^{c}} |x|(\bar{\mu}_{k+1}^{N}(dx)+ \bar{\mu}_{k}^{N} Q^{ \bar{\eta}_{k}^{N},\bar{\mu}_{k}^{N}}(dx))\right\}<\epsilon .$$ 
Let $\Lip_{1}^{0}(\mathbb{R}^{d}):= \{f\in \Lip_{1}(\mathbb{R}^{d}): f(0)=0\}$.
Then,
\begin{eqnarray}
E\sup_{\phi \in \mbox{Lip}_{1}(\mathbb{R}^{d})}|\langle\phi, \bar{\mu}_{k+1}^{N}-\bar{\mu}_{k}^{N} Q^{ \bar{\eta}_{k}^{N},\bar{\mu}_{k}^{N}}\rangle| &=& E\sup_{\phi \in \Lip_{1}^{0}(\mathbb{R}^{d})}|\langle\phi,\bar{\mu}_{k+1}^{N}-\bar{\mu}_{k}^{N} Q^{ \bar{\eta}_{k}^{N},\bar{\mu}_{k}^{N}}\rangle|\nonumber\\
&\leq& E\sup_{\phi \in \Lip_{1}^{0}(\mathbb{R}^{d})}|\langle\phi. 1_{K_{\epsilon}}, \bar{\mu}_{k+1}^{N}-\bar{\mu}_{k}^{N} Q^{ \bar{\eta}_{k}^{N},\bar{\mu}_{k}^{N}}\rangle| +\epsilon.\label{1.3.4} 
\end{eqnarray}
We will now apply Lemma \ref{app12} in the Appendix. Note that for any $\phi\in \Lip_{1}^{0}(\mathbb{R}^{d})$, $\sup_{x \in K_{\epsilon}} |\phi(x)|\leq diam(K_{\epsilon}) := m_{\epsilon}.$

Thus with notation as in Lemma \ref{app12} 
\begin{eqnarray}
 \sup_{\phi \in \Lip_{1}^{0}(\mathbb{R}^{d})}|\langle\phi. 1_{K_{\epsilon}}, \bar{\mu}_{k+1}^{N}-\bar{\mu}_{k}^{N} Q^{ \bar{\eta}_{k}^{N},\bar{\mu}_{k}^{N}}\rangle| 
&\leq& \max_{\phi \in \mathcal{F}^{\epsilon}_{m_{\epsilon,1}}(K_{\epsilon})}|\langle\phi , \bar{\mu}_{k+1}^{N}-\bar{\mu}_{k}^{N} Q^{ \bar{\eta}_{k}^{N},\bar{\mu}_{k}^{N}}\rangle| +2\epsilon.\label{1.3.5}
\end{eqnarray}
where we have denoted the restrictions of $\bar{\mu}_{k+1}^{N}$ and $\bar{\mu}_{k}^{N} Q^{ \bar{\eta}_{k}^{N}}$ to $K_{\epsilon}$ by the same symbols.
Using the above inequality in (\ref{1.3.4}), we obtain 
\begin{eqnarray}
E\mathcal{W}_{1}( \bar{\mu}_{k+1}^{N},\bar{\mu}_{k}^{N} Q^{ \bar{\eta}_{k}^{N},\bar{\mu}_{k}^{N}})  &\leq& \sum_{\phi \in \mathcal{F}^{\epsilon}_{m_{\epsilon,1}}(K_{\epsilon})} E|\langle \phi, \bar{\mu}_{k+1}^{N}-\bar{\mu}_{k}^{N} Q^{ \bar{\eta}_{k}^{N},\bar{\mu}_{k}^{N}}  \rangle| +3\epsilon.\label{1.3.6}
\end{eqnarray}
Using Lemma \ref{app22} we see that the first term on the right hand side can be bounded by $\frac{2m_{\epsilon}| \mathcal{F}^{\epsilon}_{m_{\epsilon,1}}(K_{\epsilon})|}{\sqrt{N}}$.

Consider the second term of R.H.S of (\ref{eq1.3.1}). From Assumption \ref{As3} applying DCT one has
\beqn
\gr \bar{\eta}_{k}^{N}(y)&=& (1-\alpha)\int S^{M}(\bar{\eta}_{k-1}^{N})(dx)\gr_{y}P(x,y) +\alpha\int \bar{\mu}_{k}^{N}(dx)\gr_{y}P'(x,y),\label{eq1.3.1.1}\\ \gr \eta_{k}(y) &=& (1-\alpha)\int  \eta_{k-1}(dx)\gr+\alpha \int \mu_k(dx)\gr_{y}P'(x,y). 
\eeqn 
Suppose $\bar{X}_{k}$ is a random variable conditioned on $\mathcal{F}_{k}^{M,N}$ is distributed with law $\bar{\mu}_{k}^{N}$. Then almost surely  $\mathcal{W}_{1}(\bar{\mu}_{k}^{N} Q^{ \bar{\eta}_{k}^{N},\bar{\mu}_{k}^{N}}, \bar{\mu}_{k}^{N}Q^{ \eta_{k},\bar{\mu}_{k}^{N}})$ is 
\beqn
&\le&\sup_{g\in \mbox{Lip}_{1}(\R^{d})} E_{k}^{M,N}\bigg[\bigg|g(A\bar{X}_{k}+\delta f_{\delta}(\gr \bar{\eta}_{k}^{N}(\bar{X}_{k}),\bar{\mu}_{k}^{N},\bar{X}_{k},\epsilon))\nonumber\\&&-g(A\bar{X}_{k}+\delta f_{\delta}(\gr \eta_{k}(\bar{X}_{k}),\bar{\mu}_{k}^{N},\bar{X}_{k},\epsilon))\bigg|\bigg]\le \delta \sigma E_{k}^{M,N}\left[\left|\gr \bar{\eta}_{k}^{N}(\bar{X}_{k})- \gr \eta_{k}(\bar{X}_{k})\right|\rt]\non\\
&\le&\delta \sigma (1-\alpha)\int\bigg| \int\left\{S^{M}(\bar{\eta}_{k}^{M})-\eta_{k} \rt\}(dx).\gr_{y}P(x,y) \bigg|\bar{\mu}_{k}^{N}(dy)\non\\&&+ \delta \sigma \alpha\int\bigg| \int\left\{\bar{\mu}_{k}^{N}-\mu_{k}\rt\}(dx).\gr_{y}P'(x,y) \bigg|\bar{\mu}_{k}^{N}(dy)\non\\
&\le& \delta \sigma (1-\alpha)l^{\gr}_{P}\mathcal{W}_{1}(S^{M}(\bar{\eta}_{k}^{M}),\eta_{k})+ \delta \sigma \alpha l^{\gr}_{P'}\mathcal{W}_{1}(\bar{\mu}_{k}^{N},\mu_{k}).\label{eq1.3.2}
\eeqn
(\ref{eq1.3.2}) follows by using Assumption \ref{As3}. About the first term in (\ref{eq1.3.2}) note that from triangular inequality,
\beqn
E\mathcal{W}_{1}(S^{M}(\bar{\eta}_{k}^{M}),\eta_{k})\le E\mathcal{W}_{1}(S^{M}(\bar{\eta}_{k}^{M}),\bar{\eta}_{k}^{M})+E\mathcal{W}_{1}(\bar{\eta}_{k}^{M},\eta_{k}).\label{eq1.3.3}
\eeqn
  The first term in (\ref{eq1.3.3}) can be written as 
\beqn
 E\mathcal{W}_{1}(S^{M}(\bar{\eta}_{k}^{M}),\bar{\eta}_{k}^{M}) &\le& E \sup_{f\in \Lip^{0}_1(\R^d)}|\left<f.1_{K_{k,\epsilon}}, S^{M}(\bar{\eta}_{k}^{M}) -\bar{\eta}_{k}^{M}\right>|+E\left<|x|.1_{K^c_{k,\epsilon}},S^{M}(\bar{\eta}_{k}^{M})\right>\nonumber\\&+&E\left<|x|.1_{K^c_{k,\epsilon}},\bar{\eta}_{k}^{M}\right>.\label{eq1.3.4}
\eeqn 
By Lemma \ref{lemth1}(b), for a specified $\eps>0,$ one can construct a compact set $K_{k,\epsilon}$ containing $0$ such that,
$$\sup_{M,N\ge 1} E \left<|x|.1_{K_{k,\epsilon}}, S^{M}(\bar{\eta}_{k}^{M}) + \bar{\eta}_{k}^{M}\rt> <\eps.$$
 Denote $m_{k,\epsilon}=\text{diam}(K_{k,\epsilon})$. Using Lemma \ref{app12} we have the L.H.S of (\ref{eq1.3.4})
\beqn
EE_{k}^{M,N}\bigg[ \sup_{\phi\in \text{Lip}^{0}_{1}(\R^d)}|\left<\phi.1_{K_{k,\epsilon}}, S^{M}(\bar{\eta}_{k}^{M}) -\bar{\eta}_{k}^{M}\right>|\bigg]+\eps\le  E E_{k}^{M,N}\bigg[\max_{\phi \in \mathcal{F}^{\epsilon}_{m_{k,\epsilon,1}}(K_{k,\epsilon})}|\left<\phi, S^{M}(\bar{\eta}_{k}^{M}) -\bar{\eta}_{k}^{M}\right>|\bigg]+2\eps\non
\eeqn
where (\ref{eq1.3.4}) follows from similar arguments used in (\ref{1.3.6}). Note that the Lemma \ref{lemth1} also suggests the compact set $K_{k,\eps}$ is non-random,  which only depends on $k$ and $\eps$ only. So from the display above we have
\beqn
E E^{M,N}_{k}\bigg[\sum_{\phi \in \mathcal{F}^{\epsilon}_{m_{k,\epsilon,1}}(K_{k,\epsilon})}|\left<\phi, S^{M}(\bar{\eta}_{k}^{M}) -\bar{\eta}_{k}^{M}\right>| \bigg]+2\eps\le \sum_{\phi \in \mathcal{F}^{\epsilon}_{m_{k,\epsilon,1}}(K_{k,\epsilon})}E|\left<\phi, S^{M}(\bar{\eta}_{k}^{M}) -\bar{\eta}_{k}^{M}\right>|+2\eps\,\,\,\,\,\label{eq1.3.5.5}
\eeqn
Using Lemma \ref{app22} we get the final bound of the first term in RHS of (\ref{eq1.3.5.5}) as $\frac{2m_{k,\epsilon}| \mathcal{F}^{\epsilon}_{m_{k,\epsilon,1}}(K_{k,\epsilon})|}{\sqrt{M}}$. Combining this estimate with (\ref{eq1.3.5.1}),(\ref{1.3.6}) and (\ref{eq1.3.2}) we now have
\beqn
E\mathcal{W}_{1}(\bar{\mu}_{k+1}^{N},\mu_{k+1})\le\left(\|A\|+\delta \sigma(2+ l_{PP'}^{\gr})+\delta \sigma \alpha l^{\gr}_{P'}\rt)E\mathcal{W}_{1}(\bar{\mu}_{k}^{N},\mu_{k})+ \delta \sigma (1-\alpha)l^{\gr}_{P}E\mathcal{W}_{1}(\bar{\eta}_{k}^{M},\eta_{k})\nonumber\\
+\frac{2 \delta \sigma (1-\alpha)l^{\gr}_{P}m_{k,\epsilon}| \mathcal{F}^{\epsilon}_{m_{k,\epsilon,1}}(K_{k,\epsilon})|}{\sqrt{M}}+\frac{2m_{\epsilon}| \mathcal{F}^{\epsilon}_{m_{\epsilon,1}}(K_{\epsilon})|}{\sqrt{N}}+\left(3+2 \delta \sigma (1-\alpha)l^{\gr}_{P}\rt)\eps.\s \label{finthm1}
\eeqn
For the term $E\mathcal{W}_{1}(\bar{\eta}_{k+1}^{M},\eta_{k+1}),$ we start with the following recursive form
\beqn
\bar{\eta}_{k+1}^{M}-\eta_{k+1} = (1-\alpha)\left[S^{M}(\bar{\eta}_{k}^{M})-\bar{\eta}_{k}^{M}\rt]P +(1-\alpha)\left[\bar{\eta}_{k}^{M}-\eta_{k}\rt]P + \alpha\left[\bar{\mu}_{k}^{N} - \mu_{k}\rt]P'
\eeqn
which leads to the following inequality
\beqn
\mathcal{W}_{1}(\bar{\eta}_{k+1}^{M},\eta_{k+1}) \le (1-\alpha)l(P)\mathcal{W}_{1}(S^{M}(\bar{\eta}_{k}^{M}),\bar{\eta}_{k}^{M}) +(1-\alpha)l(P)\mathcal{W}_{1}(\bar{\eta}_{k}^{M},\eta_{k}) + \alpha l(P')\mathcal{W}_{1}(\bar{\mu}_{k}^{N},\mu_{k}).\s
\eeqn
Using earlier estimates one has the final estimate for
\beqn
E\mathcal{W}_{1}(\bar{\eta}_{k+1}^{M},\eta_{k+1}) &\le&  2 (1-\alpha)l(P)\frac{m_{k,\epsilon}| \mathcal{F}^{\epsilon}_{m_{k,\epsilon,1}}(K_{k,\epsilon})|}{\sqrt{M}}+(1-\alpha)l(P)\mathcal{W}_{1}(\bar{\eta}_{k}^{M},\eta_{k}) + \alpha l(P')\mathcal{W}_{1}(\bar{\mu}_{k}^{N},\mu_{k})\nonumber\\&+&2 (1-\alpha)l(P)\eps.\label{finthm2}
\eeqn
Adding (\ref{finthm1}) and (\ref{finthm2}), using induction hypothesis and sending $M,N \to \infty$ we have 
$$E\mathcal{W}_{1}(\bar{\mu}_{k+1}^{N},\mu_{k+1}) +E\mathcal{W}_{1}(\bar{\eta}_{k+1}^{M},\eta_{k+1}) \le \left(3+2 \delta \sigma (1-\alpha)l^{\gr}_{P}+2 (1-\alpha)l(P)\rt)\eps.$$
Since $\eps >0$ arbitrary, the result follows.

Part (a) can be proved similarly. The change will come from the structural difference of $\bar{\eta}^{N}_{k}$ and $\eta^{N}_{k}$ because of the change in the updating kernel. So the term coming from the quantity $S^{M}(\bar{\eta}_{k}^{M})-\bar{\eta}_{k}^{M}$ won't appear here. Hence we get the following final estimate 
\beqn
E\left[ \clw_1(\mu_{k+1}^{N},\mu_{k+1})+ \clw_1(\eta_{k+1}^{N},\eta_{k+1})\right]\le\left[\|A\|+\delta \sigma(2+ l_{PP'}^{\gr})+\delta \sigma \alpha l^{\gr}_{P'}+\alpha l(P')\rt] E\mathcal{W}_{1}(\mu_{k}^{N},\mu_{k})\non\\+ \left[\delta \sigma (1-\alpha)l^{\gr}_{P}+(1-\alpha)l(P)\rt]E\mathcal{W}_{1}(\eta_{k}^{M},\eta_{k})+3\eps+\frac{2m_{\epsilon}| \mathcal{F}^{\epsilon}_{m_{\epsilon,1}}(K_{\epsilon})|}{\sqrt{N}}\non
\eeqn
from which the result follows by induction.

$\square$
\subsection{Proof of Proposition \ref{Thmstable}}

The techniques that we used is very similar with the contraction based method that was used in \cite{budhiraja2011discrete}. We will start with the following lemma and then prove the Proposition \ref{Thmstable} using it. Define the following distance on $\clp_{1}(\R^{d})\times\clp_{1}^{*}(\R^{d})$ for $(\mu_{n},\eta_{n}),(\mu'_{n},\eta'_{n})\in \clp_{1}(\R^{d})\times\clp_{1}^{*}(\R^{d})$ $$\mathcal{W}_{1}((\mu_{n},\eta_{n}),(\mu'_{n},\eta'_{n})):= \mathcal{W}_{1}(\mu_{n},\mu'_{n})+\mathcal{W}_{1}(\eta_{n},\eta'_{n}).$$ Note that it is a complete separable metric of the space $\clp_{1}(\R^{d})\times\clp_{1}^{*}(\R^{d}).$

\begin{lemma}\label{lem1}
Let $\mu_{0},\mu'_{0}\in \clp_{1}(\R^{d})$ and $\eta_{0},\eta'_{0}\in \clp_{1}^{*}(\R^{d})$. Suppose Assumptions \ref{As2},\ref{As2.5}, \ref{As3} and \ref{As7} hold.  Then the transformation $\Psi: \clp_{1}(\R^{d})\times\clp_{1}^{*}(\R^{d})\to \clp_{1}(\R^{d})\times\clp_{1}^{*}(\R^{d})$ is well defined if following hold 
\beqn
\delta  <a_{0}\s\s\text{and}\s\s\sum_{i=0}^{\infty}(1-\alpha)^{i}\int_{\R^d}|y|P'P^{i}(0,dy)<\infty.\label{cond0.5}
\eeqn
Moreover if Assumptions \ref{As3},\ref{As6} and \ref{As7} hold with the following condition:
\beqn
\max\Big\{\Big(\|A\|+\delta \sigma(2+ l_{PP'}^{\gr,\alpha}) + \alpha l(P')\Big ),(1-\alpha)l(P) \Big\}+\delta \sigma \max\big\{\alpha l_{P'}^{\gr},(1-\alpha) l_{P}^{\gr}\big\} < 1. \label{condition1}
\eeqn
Then there exist a $\theta\in(0,1)$ and a constant $a_1\in (0,\infty)$ such that for any $n\in \mathbb{N},$
$$\mathcal{W}_{1}\big(\Psi^{n}(\mu_{0},\eta_{0}),\Psi^{n}(\mu'_{0},\eta'_{0})\big) \leq a_{1}\theta^{n}.$$
\end{lemma}
\begin{remark}
The condition (\ref{condition1}) implies the first condition of (\ref{cond0.5}) while the second one is very general.
\end{remark}

\subsubsection{Proof of Lemma \ref{lem1}}
 For fixed $\mu_{0},\mu'_{0}\in \clp_{1}(\R^{d})$ and $\eta_{0},\eta'_{0}\in \clp_{1}^{*}(\R^{d})$ define the following quantities for $n\geq 1$
$$(\mu_{n},\eta_{n})=\Psi^{n}(\mu_{0},\eta_{0}),\quad\quad (\mu'_{n},\eta'_{n})=\Psi^{n}(\mu'_{0},\eta'_{0}) \quad \text{and} \quad \Psi^{0} = I.$$
First we will show that under transformation $\Psi$ the $(\mu_n,\nu_{n})\in\clp_1(\R^d)\times \clp^*_1(\R^d)$ for $(\mu_0,\nu_0)\in\clp_1(\R^d)\times \clp^{*}_1(\R^d),$ so that the quantity $\mathcal{W}_1(\mu_n,\mu'_{n})+\mathcal{W}_1(\nu_n,\nu'_{n})$ is well defined. Note that , if $\delta\in(0,a_0),$ then $\gamma =\|A\|+\delta \sigma\lt(2+  l^{\gr,\alpha}_{PP'}\rt) \in(0,1),$ implying $$\left<|x|,\mu_n\right> \le \gamma^{n}\left<|x|,\mu_0\right>+\frac{\delta [\sigma c^{\gr,\alpha}_{PP'}+\sigma_2]}{1-\gamma},$$ 
which follows similarly from the proof of  Lemma \ref{l0}(a). It means if $\delta\in (0,a_0)$ and $\left<|x|,\mu_0\right> <\infty$ hold, then $\mu_n \in \clp_1(\R^d)$ for all $n\ge 1$. Under conditions in (\ref{cond0.5}) one also has $\sup_{n>0}\left<|x|,\eta_{n}\right><\infty$ for all $n\in \mathbb{N}.$
One has $\gr \eta_{n+1}(y) =\int_{\R^{d}} \eta_{n}(x) [\gr_{y} R^{\alpha}_{\mu_{n}}(x,y)]dx$
 by Assumption \ref{As3} using DCT. From that condition it follows that for any $n\ge 1$, $\|\gr \eta_{n}(\cdot)\|_{1} < (1-\alpha)l_P^{\gr} +\alpha l^{\gr}_{P'}=l^{\gr,\alpha}_{PP'}<\infty$ showing $\eta_n \in \clp^*_1(\R^d)$ for all $n>0$ if $\eta_0 \in \clp^*_1(\R^d)$.

Now we will go back to the proof of the second part of the lemma regarding the contraction part. Assume $n\geq 2$. The first term of $\mathcal{W}_{1}((\mu_{n},\eta_{n}),(\mu'_{n},\eta'_{n}))$ can be expressed as
\beqn
\mathcal{W}_{1}(\mu_{n},\mu'_{n})= \mathcal{W}_{1}(\mu_{n-1}Q^{\eta_{n-1},\mu_{n-1}},\mu'_{n-1}Q^{\eta'_{n-1},\mu'_{n-1}}) &\leq& \mathcal{W}_{1}(\mu_{n-1}Q^{\eta_{n-1},\mu_{n-1}},\mu'_{n-1}Q^{\eta_{n-1},\mu'_{n-1}})  \nonumber\\ &+& \mathcal{W}_{1}(\mu'_{n-1}Q^{\eta_{n-1},\mu'_{n-1}},\mu'_{n-1}Q^{\eta'_{n-1},\mu'_{n-1}}) \nonumber\\ &=:& T_{1} +T_{2} .
\eeqn
\beqn
T_{1}&=&\mathcal{W}_{1}(\mu_{n-1}Q^{\eta_{n-1},\mu_{n-1}},\mu'_{n-1}Q^{\eta_{n-1},\mu'_{n-1}})  \leq \inf_{\{X,Y : \mathcal{L}(X,Y) = (\mu_{n-1},\mu'_{n-1}),  X,Y \perp \epsilon\}} E\big|A(X-Y) \nonumber\\&&+ \delta [f_{\delta}(\gr \eta_{n-1}(X),\mu_{n-1},X,\epsilon)-f_{\delta}(\gr \eta_{n-1}(Y),\mu'_{n-1},Y,\epsilon)]\big|\nonumber\\
&\leq& \inf_{\{X\sim \mu_{n-1},Y\sim \mu'_{n-1}\}} \big\{(\|A\|+\delta \sigma)E|X-Y| + \delta \sigma  E|\gr \eta_{n-1}(X) - \gr \eta_{n-1}(Y)|\big\}\non\\&&+\delta\sigma \clw_{1}(\mu_{n-1},\mu'_{n-1})\label{us}
\eeqn
The last inequality (\ref{us}) follows from Assumption \ref{As2}. As a consequence of Assumption \ref{As3} from (\ref{gr})  it follows that
\beqn
|\gr \eta_{n+1}(X) -\gr \eta_{n+1}(Y) | &\leq& \int_{\R^{d}}\eta_{n}(x)|\gr_{y} R^{\alpha}_{\mu_{n}}(x,X)-\gr_{y} R^{\alpha}_{\mu_{n}}(x,Y)|dx\nonumber\\
&\leq& (1-\alpha) \int_{\R^{d}}\eta_{n}(x)|\gr_{y}P(x,X)  - \gr_{y}P(x,Y) |dx\nonumber\\&&+\alpha |\gr_{y}\mu_{n}P'(X) - \gr_{y}\mu_{n}P'(Y)|\nonumber\\
&\leq& l_{PP'}^{\gr,\alpha}|X-Y|.\label{gradeq1}
\eeqn
With that estimate, taking infimum at R.H.S of (\ref{us}) with all possible couplings of $(X,Y)$ with marginals respectively $\mu_{n-1}$ and $\mu'_{n-1}$, one gets 
\beqn
T_{1}= \mathcal{W}_{1}\big(\mu_{n-1}Q^{\eta_{n-1},\mu_{n-1}},\mu'_{n-1}Q^{\eta_{n-1},\mu'_{n-1}}\big)  \leq \big(\|A\|+\delta \sigma(2+ l_{PP'}^{\gr,\alpha})\big)\mathcal{W}_{1}(\mu_{n-1},\mu'_{n-1}).\label{A1}
\eeqn
Let $X$ be a $\R^{d}$ valued random variable with law $\mu'_{n-1}$. Now about the term $T_{2}$,
\beqn
T_{2}&=&\mathcal{W}_{1}\big(\mu'_{n-1}Q^{\eta_{n-1},\mu'_{n-1}},\mu'_{n-1}Q^{\eta'_{n-1},\mu'_{n-1}}\big)\nonumber\\&\leq&\sup_{g\in \mbox{Lip}_{1}(\R^{d})} E\left|g(AX+\delta f_{\delta}(\gr \eta_{n-1}(X),\mu'_{n-1},X,\epsilon))-g(AX+\delta f_{\delta}(\gr \eta'_{n-1}(X),\mu'_{n-1},X,\epsilon))\right|\nonumber\\
&\leq& \delta \sigma E\left|\gr \eta_{n-1}(X) - \gr \eta'_{n-1}(X)\right|\nonumber\\
&\leq& \delta \sigma E\left|\int_{\R^{d}}\eta_{n-2}(x)(\gr_{y} R^{\alpha}_{\mu_{n-2}}(x,X))dx -\int_{\R^{d}}\eta'_{n-2}(x)(\gr_{y} R^{\alpha}_{\mu'_{n-2}}(x,X))dx\right| \non\\
&\leq& \alpha\delta \sigma \int_{\R^{d}}\eta_{n-2}(x)E\lt|\gr_{y} \mu_{n-2}P'(X) - \gr_{y} \mu'_{n-2}P'(X) \rt|dx  \non\\&&+(1-\alpha)\delta  \sigma E\lt| \int_{\R^{d}}\gr_{y} P(x,X)(\eta_{n-2}(x)-\eta'_{n-2}(x))dx \rt| \non\\
&=:& T^{(1)}_{2}+ T^{(2)}_{2}\label{th2eq1}
\eeqn
Note that 
\beqn
T^{(1)}_{2}:=\alpha\delta \sigma \int_{\R^{d}}\eta_{n-2}(x)\int_{\R^d}\mu'_{n-1}(dz)\lt|\int_{\R^{d}}\bigg( \mu_{n-2}(dy)\gr_{y}P'(y,z) -  \mu'_{n-2}(dy)\gr_{y}P'(y,z)\bigg) \rt|dx\label{B_{1}}
\eeqn
Since from Assumption \ref{As3} $\gr_{y}P'(\tilde{x},x)$ is a Lipschitz function with coefficient $l^{\gr}_{P'}$, the first integrand in (\ref{B_{1}}) will be bounded by $l^{\gr}_{P'}. \mathcal{W}_{1}(\mu_{n-2},\mu'_{n-2})$  which gives
\beqn
T^{(1)}_{2} \leq \alpha\delta \sigma l^{\gr}_{P'} \mathcal{W}_{1}(\mu_{n-2},\mu'_{n-2}). \label{fin B_{1}}
\eeqn
Now using Assumption \ref{As6} the second term $T^{(2)}_{2}$ gives similarly
\beqn
T^{(2)}_{2} &=&(1-\alpha)\delta  \sigma E\lt| \int_{\R^{d}}\gr_{y} P(x,X)(\eta_{n-2}(x)-\eta'_{n-2}(x))dx \rt| \non\\&\leq&   (1-\alpha)\delta\sigma \int_{\R^{d}} \lt|   \int_{\R^{d}}\gr_{y} P(x,y)\{\eta_{n-2}(x)-\eta'_{n-2}(x) \}dx         \rt| \mu'_{n-1}(dy)\non\\
&\leq&  (1-\alpha)\delta\sigma  l^{\gr}_{P} \mathcal{W}_{1}(\eta_{n-2},\eta'_{n-2}).\label{B_{2}}
\eeqn
Using the Assumption \ref{As7} we have 
\beqn
\mathcal{W}_{1}(\eta_{n},\eta'_{n}) \leq (1-\alpha) l(P) \mathcal{W}_{1}(\eta_{n-1},\eta'_{n-1})   + \alpha l(P') \mathcal{W}_{1}(\mu_{n-1},\mu'_{n-1})\label{eta}
\eeqn
Combining (\ref{fin B_{1}}),(\ref{B_{2}}) and (\ref{eta}) we have the following recursion for $n\geq 2,$
\beqn
\mathcal{W}_{1}(\mu_{n},\mu'_{n})\s+\mathcal{W}_{1}(\eta_{n},\eta'_{n}) \quad &\leq&\quad  \big(\|A\|+\delta \sigma(2+ l_{PP'}^{\gr,\alpha})\big)\mathcal{W}_{1}(\mu_{n-1},\mu'_{n-1})\s\non\\+   \alpha\delta \sigma l^{\gr}_{P'} \mathcal{W}_{1}(\mu_{n-2},\mu'_{n-2}) &+&\alpha l(P').\mathcal{W}_{1}(\mu_{n-1},\mu'_{n-1})\s+ (1-\alpha)\delta\sigma  l^{\gr}_{P} \mathcal{W}_{1}(\eta_{n-2},\eta'_{n-2})\non\\&+& (1-\alpha) l(P) \mathcal{W}_{1}(\eta_{n-1},\eta'_{n-1}).\label{rec}
\eeqn
Define a sequence $a_{n} := \mathcal{W}_{1}(\mu_{n},\mu'_{n})+\mathcal{W}_{1}(\eta_{n},\eta'_{n}),$ for $n\geq 2$ and and first two terms we set them to be
\beqn
a_{0}:= \mathcal{W}_{1}(\mu_{0},\mu'_{0})+\mathcal{W}_{1}(\eta_{0},\eta'_{0}), \s\s\s\s a_{1}:= \mathcal{W}_{1}(\mu_{1},\mu'_{1})+\mathcal{W}_{1}(\eta_{1},\eta'_{1})\non
\eeqn
which are well defined for  $\mu_{0},\mu'_{0}\in \clp_{1}(\R^{d})$ and $\eta_{0},\eta'_{0}\in \clp_{1}^{*}(\R^{d}).$ Then from (\ref{rec}) and denoting $c_{1}:=\max\left\{\left(\big(\|A\|+\delta \sigma(2+ l_{PP'}^{\gr,\alpha})\big) + \alpha l(P')\right),(1-\alpha)l(P)\right\}$, $c_{2}:=\delta \sigma \max\big\{\alpha l_{P'}^{\gr},(1-\alpha) l_{P}^{\gr}\big\}$ following holds
\beqn
a_{n} \leq c_{1} a_{n-1} +c_{2}  a_{n-2}\,\,\label{cond2}
\eeqn
for $n\geq 2.$
Given $(\omega,\delta,\alpha)$ if there exists a $\theta \in (0,1)$ for which the following inequality holds 
\beqn
\frac{c_{1}}{\theta} +\frac{ c_{2}}{\theta^{2}} \leq 1,\label{ineq}
\eeqn
then denoting $\lambda = \frac{c_2}{\theta},$ we have 
\beqn
a_{n} \leq \bigg[ \theta \bigg(1 - \frac{\lambda}{\theta}\bigg)\bigg]a_{n-1} + \theta \lambda a_{n-2}\quad \Leftrightarrow \s
a_{n}+\lambda a_{n-1} &\leq&  \theta (a_{n-1}+\lambda a_{n-2}).
\eeqn

Existence of a solution $\theta \in (0,1)$ satisfying (\ref{ineq}) is valid under $c_{1}+c_{2}<1$ which is equivalent to the condition
\beqn
\max\Big\{\Big(\big(\|A\|+\delta \sigma(2+ l_{PP'}^{\gr,\alpha})\big) + \alpha l(P')\Big ),(1-\alpha)l(P) \Big\}+\delta \sigma \max\big\{\alpha l_{P'}^{\gr},(1-\alpha) l_{P}^{\gr}\big\} < 1\label{cond}
\eeqn
in (\ref{condition1}) satisfied by $(\delta,\alpha,\|A\|)$. From (\ref{cond}) it follows   
$$a_{n}\leq a_{n}+\lambda a_{n-1} \leq \theta^{n-1} [a_{1}+\lambda a_{0}]$$
for $n\ge 2$. Since
\beqn
\mathcal{W}_{1}(\eta_{1},\eta'_{1}) &=& \mathcal{W}_{1}(\eta_{0}R^{\alpha}_{\mu_{0}},\eta'_{0}R^{\alpha}_{\mu'_{0}})\leq(1-\alpha) l(P) \mathcal{W}_{1}(\eta_{0},\eta'_{0})   + \alpha l(P') \mathcal{W}_{1}(\mu_{0},\mu'_{0}),\non\\
\mathcal{W}_{1}(\mu_{1},\mu'_{1}) &=& \mathcal{W}_{1}(\mu_{0}Q^{\eta_{0},\mu_{0}}, \mu'_{0}Q^{\eta'_{0},\mu'_{0}})\leq  \mathcal{W}_{1}(\mu_{0}Q^{\eta_{0},\mu_{0}}, \mu'_{0}Q^{\eta_{0},\mu'_{0}})+\mathcal{W}_{1}(\mu'_{0}Q^{\eta_{0},\mu'_{0}}, \mu'_{0}Q^{\eta'_{0},\mu'_{0}}) \non\\
&\leq&\Big (\|A\|+\delta\sigma(2 +l_{PP'}^{\gr})\Big)\mathcal{W}_{1}(\mu_{0},\mu'_{0})+\delta \sigma E\lt| \gr \eta_{0}(X) - \gr\eta'_{0}(X)\rt|\non
\eeqn
where $X \sim \mu'_{0}.$
Final estimate for $a_{n}$ is
$$a_{n}\leq \theta^{n-1}\lt[\Big (\max\Big\{\Big(\|A\|+\delta \sigma(2+ l_{PP'}^{\gr,\alpha}) + \alpha l(P')\Big ),(1-\alpha)l(P) \Big\}+\lambda \Big)a_{0}+\delta \sigma E\lt| \gr \eta_{0}(X) - \gr\eta'_{0}(X)\rt|\rt].$$
 Since $X \sim \mu'_{0}\in \clp_{1}(\R^d)$ and $\gr\eta_0,\gr\eta'_{0}$  have linear growth (since $\eta_0,\eta'_{0}\in \clp^{*}_{1}(\R^{d})$), the second term inside the bracket is finite.   
A general formula can be observed for $a_n$
\beqn
\mathcal{W}_{1}(\Psi^{n}(\mu_{0},\eta_{0}),\Psi^{n}(\mu'_{0},\eta'_{0}))\le \theta^{n}\left[ a\mathcal{W}_{1}((\mu_{0},\eta_{0}),(\mu'_{0},\eta'_{0}))+b \mathcal{W}_{1}(\mu'_{0}Q^{\eta_{0},\mu'_{0}}, \mu'_{0}Q^{\eta'_{0},\mu'_{0}}) \rt]\label{use}
\eeqn
where $$a=\frac{\max\Big\{\Big(\|A\|+\delta \sigma(2+ l_{PP'}^{\gr,\alpha}) + \alpha l(P')\Big ),(1-\alpha)l(P) \Big\}+\lambda}{\theta} ,\s\s b=\frac{1}{\theta}.$$ Observe that the quantity inside the bracket of RHS of (\ref{use}) is finite for  $\mu_{0},\mu'_{0}\in \clp_{1}(\R^{d})$ and $\eta_{0},\eta'_{0}\in \clp_{1}^{*}(\R^{d})$. Hence proved the lemma.

$\square$

	We now complete the proof of the theorem. Given $l(PP')<1$ from Assumption (\ref{As7}), one can always find $(\omega_{0},\alpha_{0},\delta_{0})\in (0,1)\times(0,1)\times(0,1)$ for which (\ref{cond}) holds under $$\|A\| <\omega_{0},\quad \alpha<\alpha_{0},\quad\delta<\delta_{0}.$$For existence we need to show that under $\mathcal{W}_{1}\lt((\cdot,\cdot),(\cdot,\cdot)\rt)$ distance $\clp_{1}(\R^{d})\times \clp_{1}^{*}(\R^{d})$ is complete. From Lemma \ref{lem1} one can choose $(\omega,\alpha,\delta)$ such that (\ref{condition1}) holds. 
 It follows that using the $\theta$ from that lemma the sequence $\{\Psi^{n}(\mu_{0},\eta_{0})\}_{n\geq 1}^{\infty}$ is a cauchy sequence in $\clp_{1}(\R^{d})\times \clp_{1}(\R^{d})$ which is a complete metric space under $\mathcal{W}_{1}\lt((\cdot,\cdot),(\cdot,\cdot)\rt).$ So there exists a $(\mu_{\infty},\eta_{\infty}) \in \clp_{1}(\R^{d})\times \clp_{1}^{}(\R^{d})$ such that $\Psi^{n}(\mu_{0},\eta_{0}) \to (\mu_{\infty},\eta_{\infty})$ as $n\to\infty$. Our assertion for existence will be proved if we prove  $\eta_{\infty} \in \mathcal{P}_{1}^*(\R^{d}).$
 Given the initial conditon $\|\gr\eta_{0}(x)\|_{1}<\infty,$ we will always have from (\ref{gr}) $\|\gr\eta_{k}(x)\|_{1}<\infty \quad \forall\quad k>1.$
Note that for $\eta_{0}\in \clp^{*}_{1}(\R^{d})$, one has $\eta_{k} \in  \clp^{*}_{1}(\R^{d})$ for all $k.$ This implies $\eta_{\infty} \in  \clp^{*}_{1}(\R^{d}).$ So
$$(\mu_{\infty},\eta_{\infty}) \in \clp_{1}(\R^{d})\times \clp_{1}^{*}(\R^{d}) .$$ 
Observe further for $\theta\in (0,1)$ in (\ref{use}) of Lemma \ref{lem1} 
\beqn
\clw_{1}\big(\Psi^{n}(\mu_0,\eta_0),(\mu_{\infty},\eta_{\infty})\big)&=&\clw_{1}\big(\Psi^{n}(\mu_0,\eta_0),\Psi^{n}(\mu_{\infty},\eta_{\infty})\big)\non\\&\le& \theta^{n}\big[a\mathcal{W}_{1}\big((\mu_{0},\eta_{0}),(\mu_{\infty},\eta_{\infty})\big)+b\mathcal{W}_{1}(\mu_{\infty}Q^{\eta_{0},\mu_{\infty}}, \mu_{\infty}Q^{\eta_{\infty},\mu_{\infty}})\big].\label{uniq}
\eeqn
Uniqueness of fixed points follows immediately from (\ref{uniq}).

$\square$

\subsection{Proof of Theorem \ref{unifintime}}
We will prove part (b) of the theorem. Part (a) will follow similarly. We need to prove the following Lemma first.

\begin{lemma}\label{lemthm3.1}
Consider the second particle system $\mathbb{IPS}_{2}.$ Suppose that Assumptions \ref{unifas1},\ref{unifas2} hold. Denote $N_{1}=\min{\{N,M\}}.$
Then there exist a constant $C\in(0,\infty)$ such that the upper-bound $b(\tau,d)$ of the quantity $\sup_{k\ge 1}E\mathcal{W}_1\big((\bar{\mu}_{k}^{N},\bar{\eta}_{k}^{M}),\Psi(\bar{\mu}_{k-1}^{N},\bar{\eta}_{k-1}^{M})\big)$ can be given as $b(N_1,\tau,d)$ as defined in Theorem \ref{unifintime}. The constant $C$ will vary for dfferent cases.
\end{lemma}

\subsubsection{Proof of Lemma \ref{lemthm3.1}}
 We start with the fact that 
\beqn
&&E\mathcal{W}_1\big((\bar{\mu}_{k}^{N},\bar{\eta}_{k}^{M}),\Psi(\bar{\mu}_{k-1}^{N},\bar{\eta}_{k-1}^{M})\big)=E\mathcal{W}_{1}(\bar{\mu}_{k}^{N},\bar{\mu}_{k-1}^{N}Q^{\bar{\eta}_{k-1}^{M},\bar{\mu}_{k-1}^{N}})+E\mathcal{W}_{1}(\bar{\eta}_{k}^{M},\bar{\eta}_{k-1}^{M}R_{\bar{\mu}_{k-1}^{N}}^{\alpha})\nonumber\\
&&\le E\mathcal{W}_{1}(\bar{\mu}_{k}^{N},\bar{\mu}_{k-1}^{N}Q^{\bar{\eta}_{k-1}^{M},\bar{\mu}_{k-1}^{N}})+(1-\alpha)E\mathcal{W}_{1}(S^M(\bar{\eta}_{k-1}^{M}),\bar{\eta}_{k-1}^{M})\s\s\s\non\\&&=E\mathcal{W}_{1}(\bar{\mu}_{k}^{N},\bar{\mu}_{k-1}^{N}Q^{\bar{\eta}_{k-1}^{M},\bar{\mu}_{k-1}^{N}})+(1-\alpha)E\lt[E\mathcal{W}_{1}(S^M(\bar{\eta}_{k-1}^{M}),\bar{\eta}_{k-1}^{M})\big| \mathcal{F}_{k-1}^{M,N}\rt].\label{lemthm3e0}
\eeqn
In order to bound both terms in (\ref{lemthm3e0}) we borrow the following formulation from \cite{fournier2013rate} about the convergence rate of empirical distribution of iid random variables to its common distribution, where the key idea of bounding Wasserstein distance came from the constructive quantization context \cite{dereich2013constructive}. A similar idea was also developed in \cite{boissard2014mean}. We will maintain the same notation used in \cite{fournier2013rate}. Let $\clp_{l}$ be the natural partition of $(-1,1]^{d}$ into $2^{dl}$ translations of $(-2^{-l},2^{-l}]^d.$ Define a sequence of sets $\{B_{n}\}_{n\ge 0}$ such that $B_0:=(-1,1]^{d}$ and, for $n\ge 1$, $B_n:=(-2^n,2^n]^{d}\setminus (-2^{n-1},2^{n-1}]^{d}.$ For a set $F\subset\R^d$ denote the set $2^{n}F$ as $\{2^{n}x: x\in F\}.$ For any two probability measures $\mu$ and $\nu$, combining Lemma $5$ and $6$ of \cite{fournier2013rate} one has the following inequality for the Wasserstein-$1$ distance,
\beqn
\clw_{1}(\mu,\nu)\le 3C.2^{(1+\frac{d}{2})}\sum_{n\ge 0}2^{n}\sum_{l\ge 0}2^{-l}\sum_{F\in \clp_{l}}\left[\mu(2^{n}F\cap B_n)- \nu(2^{n}F\cap B_n)\rt],\label{thm3lem1e1}
\eeqn 
where $C$ is a constant depends only on $d.$ We denote $a_{k}^{i,M,N}:=\delta_{\bar{X}_{k}^{i}} - \delta_{\bar{X}_{k-1}^{i}}Q^{\bar{\eta}_{k-1}^{M},\bar{\mu}_{k-1}^{N}}.$ It follows that $\bar{\mu}_{k}^{N}-\bar{\mu}_{k-1}^{N}Q^{\bar{\eta}_{k-1}^{M},\bar{\mu}_{k-1}^{N}}=\frac{1}{N}\sum_{i=1}^{N}a_{k}^{i,M,N}.$ Note that on conditioned upon $\mathcal{F}_{k-1}^{M,N},$ the family of signed measures $\{a_{k}^{i,M,N}\}_{i=1,\ldots,M}$ is an independent class of measures while unconditionally they are just identical. Using the fact that for any set $A\in \mathcal{B}(\R^d),\s \delta_{\bar{X}_{k}^{i}}(A)\bigg| \mathcal{F}_{k-1}^{M,N} \sim \text{Bernoulli}(\delta_{\bar{X}_{k-1}^{i}}Q^{\bar{\eta}_{k-1}^{M},\bar{\mu}_{k-1}^{N}}(A)),$ we have
\beqn
E\left[\left(a_{k}^{i,M,N}(A)\rt)^2\bigg|\mathcal{F}_{k-1}^{M,N}\rt]=\delta_{\bar{X}_{k-1}^{i}}Q^{\bar{\eta}_{k-1}^{M},\bar{\mu}_{k-1}^{N}}(A)\lt[1-\delta_{\bar{X}_{k-1}^{i}}Q^{\bar{\eta}_{k-1}^{M},\bar{\mu}_{k-1}^{N}}(A)\rt]\le \delta_{\bar{X}_{k-1}^{i}}Q^{\bar{\eta}_{k-1}^{M},\bar{\mu}_{k-1}^{N}}(A)\,\label{lem3.1main}
\eeqn
which implies the unconditional expectation $E\big[\big(a_{k}^{i,M,N}(A)\big)^2\big]\le P\big[\bar{X}^i_{k-1}+\delta f_{\delta}(\gr \bar{\eta}^{M}_{k-1},\bar{\mu}_{k-1}^{N},\bar{X}^i_{k-1},\eps^{N}_{k})\in A\big].$ Using all these we have
\beqn
E\lt|\bar{\mu}_{k}^{N}(A)-\bar{\mu}_{k-1}^{N}Q^{\bar{\eta}_{k-1}^{M},\bar{\mu}_{k-1}^{N}}(A)\rt|^2=E\bigg|\frac{1}{N}\sum_{i=1}^{N}a_{k}^{i,M,N}(A)\bigg|^{2}\le\frac{E\lt[a_{k}^{i,M,N}(A)\rt]^{2}}{N}\nonumber\\
\le \frac{P\big[\bar{X}^i_{k-1}+\delta f_{\delta}(\gr \bar{\eta}^{M}_{k-1},\bar{\mu}_{k-1}^{N},\bar{X}^i_{k-1},\eps^{N}_{k})\in A\big]}{N}= \frac{ E\lt[ \delta_{\bar{X}_{k-1}^{i}}Q^{\bar{\eta}_{k-1}^{M},\bar{\mu}_{k-1}^{N}}(A)\rt]}{N}.\,\,\,\,\,\,\,\,\,\,\,\,\non
\eeqn
Using these with Cauchy-Schwarz inequality one gets following bound 
\beqn
E\lt|\bar{\mu}_{k}^{N}(A)-\bar{\mu}_{k-1}^{N}Q^{\bar{\eta}_{k-1}^{M},\bar{\mu}_{k-1}^{N}}(A)\rt|\le \min\bigg\{\sqrt{\frac{ E\big[ \delta_{\bar{X}_{k-1}^{i}}Q^{\bar{\eta}_{k-1}^{M},\bar{\mu}_{k-1}^{N}}(A)\big]}{N}},\s 2 E\lt[ \delta_{\bar{X}_{k-1}^{i}}Q^{\bar{\eta}_{k-1}^{M},\bar{\mu}_{k-1}^{N}}(A)\rt]\bigg\}\,\label{thm3ineq}
\eeqn
where second term inside the bracket of RHS of (\ref{thm3ineq}) follows trivially.
Denoting the whole constant in R.H.S of (\ref{thm3lem1e1}) as $C_d,$ we have
\beqn
E\clw_{1}(\bar{\mu}_{k}^{N},\bar{\mu}_{k-1}^{N}Q^{\bar{\eta}_{k-1}^{M},\bar{\mu}_{k-1}^{N}})\le C_d\sum_{n\ge 0}2^{n}\sum_{l\ge 0}2^{-l}E\sum_{F\in \clp_{l}}\big[\bar{\mu}_{k}^{N}(2^{n}F\cap B_n)- \bar{\mu}_{k-1}^{N}Q^{\bar{\eta}_{k-1}^{M},\bar{\mu}_{k-1}^{N}}(2^{n}F\cap B_n)\big]\label{thm3lem1e2}
\eeqn 
Note that  $\#\clp_{l}=2^{dl}$. Using Cauchy-Schwarz inequality with (\ref{thm3ineq}) and Jensen's inequality $E\sqrt{X}\le \sqrt{EX}$ for non-negative random variable $X$ , the last sum $E\sum_{F\in \clp_{l}}\big[\bar{\mu}_{k}^{N}(2^{n}F\cap B_n)- \bar{\mu}_{k-1}^{N}Q^{\bar{\eta}_{k-1}^{M},\bar{\mu}_{k-1}^{N}}(2^{n}F\cap B_n)\big]$ in the R.H.S of (\ref{thm3lem1e2}) can be bounded by
\beqn
\le \min\bigg\{2^{\frac{dl}{2}}\bigg[\frac{ E\big[ \delta_{\bar{X}_{k-1}^{i}}Q^{\bar{\eta}_{k-1}^{M},\bar{\mu}_{k-1}^{N}}(B_n)\big]}{N}\bigg]^{\frac{1}{2}},  2 E\lt[ \delta_{\bar{X}_{k-1}^{i}}Q^{\bar{\eta}_{k-1}^{M},\bar{\mu}_{k-1}^{N}}(B_n)\rt]\bigg\}.\s 
\eeqn
Now using Remark \ref{rem4.1} along with Lemma \ref{l0}, if $\delta\in(0,a(\tau))$ the quantity $\sup_{n\ge 0}\sup_{M,N \ge 1}E|\bar{X}_{n}^{i}|^{1+\tau}:=b(\tau)<\infty,$ one has by Chebyshev inequality for $n\ge 1,$
$$ \sup_{k\ge 1}E\lt[ \delta_{\bar{X}_{k-1}^{i}}Q^{\bar{\eta}_{k-1}^{M},\bar{\mu}_{k-1}^{N}}(B_{n})\rt]\le \sup_{k\ge 1}P[|\bar{X}_{k}^{i}|>2^{(n-1)}]\le \frac{b(\tau)}{2^{(1+\tau)(n-1)}} = b(\tau)2^{-(1+\tau)(n-1)}.$$
Note that $a(\tau)^{\frac{1}{1+\tau}}\to a_0$ as $\tau \to 0$ and $\delta \in(0,a_0), $ we can find $\tau_0 \in (0,a(\tau))$ such that $\delta \in (0,a(\tau_{0})^{\frac{1}{1+\tau_{0}}}).$ So the bound in (\ref{thm3lem1e2}) can be restated as 
\beqn
\sup_{k \ge 1}E\clw_{1}(\bar{\mu}_{k}^{N},\bar{\mu}_{k-1}^{N}Q^{\bar{\eta}_{k-1}^{M},\bar{\mu}_{k-1}^{N}})&\le& C_d\sum_{n\ge 0}2^{n}\sum_{l\ge 0}2^{-l}\min\bigg\{2^{\frac{dl}{2}}\sqrt{\frac{b(\tau)2^{-(1+\tau)(n-1)}}{N}},2b(\tau)2^{-(1+\tau)(n-1)}\bigg\}\non\\
&\le& C'_d\sum_{n\ge 0}2^{n}\sum_{l\ge 0}2^{-l}\min\bigg\{2^{\frac{dl}{2}}\frac{2^{-\frac{(1+\tau)n}{2}}}{\sqrt{N}},2^{-(1+\tau)n}\bigg\}. \s\label{thm3lem1e3}
\eeqn
where $b(\tau)$ is just a constant and the last display is obtained by accumulating upper bounds of all the constants to $C'_{d}$. Now proceeding exactly like step 1 to step 4 of the proof of Theorem 1 (for $p=1,q=1+\tau$) in \cite{fournier2013rate} one gets the following bounds
\begin{eqnarray}\label{total}
\sup_{k\ge 1}E\clw_{1}(\bar{\mu}_{k}^{N},\bar{\mu}_{k-1}^{N}Q^{\bar{\eta}_{k-1}^{M},\bar{\mu}_{k-1}^{N}})=C
\begin{cases}
 N^{- \max\{\frac{1}{2},\frac{\tau}{1+\tau}\}}\s\s\s\s\s\s\text{if}\s\s  d=1,\tau \neq 1,\s \\
N^{-\frac{1}{2}}\log(1+N)+N^{-\frac{\tau}{1+\tau}}\s\text{if}\s\s d=2,\tau\neq 1,\s\nonumber \\
N^{-\max\{\frac{1}{d},\frac{\tau}{1+\tau}\}}\s\s\s\s\s\s\text{if}\s\s d>2,\tau \neq \frac{1}{d-1}. 
\end{cases}
\end{eqnarray}

Now we will fill the gaps for each of the three special cases $\tau=1,\tau=1$ and $\tau=\frac{1}{d-1}$ of three regimes  respectively $d=1,d=2$ and $d>2$. We note that one can generalize the choice of $\l_{N,\ep}$ done in step 1 of Theorem 1 of \cite{fournier2013rate} where $l_{N,\ep}$ could be taken as $\frac{\frac{1}{2}\log(\ep N)}{d\log 2}\vee 0$ instead of $\frac{\log(2+\ep N)}{d\log 2}$ though it doesn't change the conclusion of the main theorem. After step $1$ with $p=1, q=1+\tau, \ep = 2^{-(1+\tau)n}$ one will get
\begin{eqnarray}\label{total}
\sum_{l\ge 0}2^{-l}\min\bigg\{2^{\frac{dl}{2}}\sqrt{\frac{\ep}{N}},\ep\bigg\}=C
\begin{cases}
\min\{\ep,\lt(\frac{\ep}{N}\rt)^{\frac{1}{2}}\} \s\s\s\s\s\s\s\,\,\,\, \text{if}\s d=1, \\
\min\{\ep,\lt(\frac{\ep}{N}\rt)^{\frac{1}{2}}[\log(\ep N)\vee 0]\}\s\s \text{if}\s d=2,\nonumber \\
\min\{\ep,\ep\lt(\ep N\rt)^{-\frac{1}{d}}\}\,\,\s\s\s\s\s\s \text{if}\s d>2,
\end{cases}
\end{eqnarray}
where the constant $C$ will vary from case to cases. Suppose $d=1.$ From (\ref{thm3lem1e3}) for general $\tau>0$ one has
\beqn
\sup_{k\ge 1}E\clw_{1}(\bar{\mu}_{k}^{N},\bar{\mu}_{k-1}^{N}Q^{\bar{\eta}_{k-1}^{M},\bar{\mu}_{k-1}^{N}})&\le&C'_{d}\sum_{n\ge 0}2^{n}\min\bigg\{2^{-(1+\tau)n},\lt(\frac{2^{-(1+\tau)n}}{N}\rt)^{\frac{1}{2}}\bigg\}.
\eeqn
Note that for $n\ge n_{N,\tau}:=\frac{\log N}{(1+\tau)\log2},$  one has $2^{-(1+\tau)n}\le\lt(\frac{2^{-(1+\tau)n}}{N}\rt)^{\frac{1}{2}}.$ So for $\tau =1,$ 
\beqn
\sum_{n\ge 0}2^{n}\min\bigg\{2^{-2n},\lt(\frac{2^{-2n}}{N}\rt)^{\frac{1}{2}}\bigg\}&\le& \sum_{n< n_{N,1}}2^{n}\lt(\frac{2^{-2n}}{N}\rt)^{\frac{1}{2}} + \sum_{n\ge n_{N,1}}2^{-n}\nonumber\\
=n_{N,1}N^{-\frac{1}{2}}+C2^{-n_{N,1}}&=& N^{-\frac{1}{2}}\frac{\log N}{2\log 2}+CN^{-\frac{1}{2}}.\label{lem3.1e}
\eeqn
For $d=2,$ from  (\ref{thm3lem1e3}) for general $\tau>0$ one has
\beqn
\sup_{k\ge 1}E\clw_{1}(\bar{\mu}_{k}^{N},\bar{\mu}_{k-1}^{N}Q^{\bar{\eta}_{k-1}^{M},\bar{\mu}_{k-1}^{N}})&\le&C'_{d}\sum_{n\ge 0}2^{n}\min\bigg\{2^{-(1+\tau)n},\lt(\frac{2^{-(1+\tau)n}}{N}\rt)^{\frac{1}{2}}\big[\log\big(2^{-(1+\tau)n}N\big)\vee 0\big]\bigg\}.\non
\eeqn
For $\tau = 1,$ $\ep=2^{-2n}.$ Note that if $n<n^{(2)}_{N}:=\log_{4}N - \log_{2}\lt(\log N\rt),$ then one has $$\ep=2^{-2n}>\lt(\frac{2^{-2n}}{N}\rt)^{\frac{1}{2}}\lt[\log\lt(2^{-2n}N\rt)\vee 0\rt].$$ 
\beqn
&& \sum_{n\ge 0}2^{n}\min\lt\{2^{-2n},\lt(\frac{2^{-2n}}{N}\lt[\log\lt(2^{-2n}N\rt)\vee 0\rt]\rt)^{\frac{1}{2}}\rt\}\non\\ &\le& \sum_{n< n^{(2)}_{N}}2^{n}\lt(\frac{2^{-2n}}{N}\rt)^{\frac{1}{2}}\lt[\log\lt(2^{-2n}N\rt)\vee 0\rt]+ \sum_{n\ge n^{(2)}_{N}}2^{-n}\le n^{(2)}_{N}\frac{\lt[\log\lt(N\rt)\vee 0\rt]}{N^{\frac{1}{2}}}+C2^{-n^{(2)}_{N}}\non\\&\le& C_{1}N^{-\frac{1}{2}}\lt[(\log N)^{2}-\log N\log_{2}(\log N) \rt]+C_{2}\frac{\log N}{\sqrt{N}}.\label{lem3.1e}
\eeqn

By proceeding similarly, for all non regular cases we will end up getting the following results (the constant $C$ will vary from case to cases):
\begin{eqnarray}\label{total}
\sup_{k\ge 1}E\clw_{1}(\bar{\mu}_{k}^{N},\bar{\mu}_{k-1}^{N}Q^{\bar{\eta}_{k-1}^{M},\bar{\mu}_{k-1}^{N}})=C
\begin{cases}
 N^{-\frac{1}{2}}\log N+N^{-\frac{1}{2}}\s\s\s\s\s\s\s\s\s\s\s\,\,\,\,\,\text{if}\s\s  d=1,\tau = 1,\s \\
N^{-\frac{1}{2}}\lt[(\log N)^{2}-\log N.\log_{2}(\log N) \rt]+\frac{\log N}{\sqrt{N}}\s\text{if}\s\s d=2,\tau = 1,\s\nonumber \\
\frac{\log_{2}N}{N^{\frac{1}{d}}}+N^{-\frac{1}{d}}\s\s\s\s\s\s\s\s\s\s\s\s\s\s\,\text{if}\s\s d>2,\tau = \frac{1}{d-1}. 
\end{cases}
\end{eqnarray}
Now about the second term of (\ref{lemthm3e0}) using (\ref{thm3lem1e1}), the upperbound of $E\mathcal{W}_{1}(S^M(\bar{\eta}_{k-1}^{M}) \bar{\eta}_{k-1}^{M})$ is
\beqn
3C2^{(1+\frac{d}{2})}\sum_{n\ge 0}2^{n}\sum_{l\ge 0}2^{-l}E\sum_{F\in \clp_{l}}\left[S^M(\bar{\eta}_{k-1}^{M})(2^{n}F\cap B_n)- \bar{\eta}_{k-1}^{M}(2^{n}F\cap B_n)\rt].\s
\eeqn
By Cauchy Schwarz inequality and using Jensen inequality $E\sqrt{X}\le \sqrt{EX}$ for a nonnegative random variable $X,$ one gets the upperbound of 
\beqn
&& E\lt[\sum_{F\in \clp_{l}}\left[S^M(\bar{\eta}_{k-1}^{M})(2^{n}F\cap B_n)- \bar{\eta}_{k-1}^{M}(2^{n}F\cap B_n)\rt]\bigg|\clf_{k-1}^{M,N}\rt]\non\\
&\le& 2^{\frac{dl}{2}}\bigg[\sum_{F\in\clp_{l}} E\big[\big(\frac{1}{M}\sum_{i=1}^{M}\delta_{Y^{i,M}_{k-1}}(2^{n}F\cap B_n)-\bar{\eta}_{k-1}^{M}(2^{n}F\cap B_n)\big)^{2}\big| \mathcal{F}_{k-1}^{M,N}\big]\bigg]^{\frac{1}{2}}.\label{lem3.1f}
\eeqn
Using similar argument used in (\ref{lem3.1main}) the R.H.S of (\ref{lem3.1f}) will be less than
\beqn
&&2^{\frac{dl}{2}}\bigg[\frac{\sum_{F\in\clp_{l}}       \bar{\eta}_{k-1}^{M}(2^{n}F\cap B_n)\lt(1-\bar{\eta}_{k-1}^{M}(2^{n}F\cap B_n)\rt)}{M}\bigg]^{\frac{1}{2}}\le 2^{\frac{dl}{2}}\bigg[\frac{\bar{\eta}_{k-1}^{M}( B_n)}{M}\bigg]^{\frac{1}{2}}\non\\ &\le& 2^{\frac{dl}{2}}\bigg[\frac{\bar{\eta}_{k-1}^{M}(x:|x|>2^{n-1})}{M}\bigg]^{\frac{1}{2}}
\le 2^{\frac{dl}{2}}\bigg[\frac{\lt<|x|^{1+\tau},\bar{\eta}_{k-1}^{M}\rt>2^{-(n-1)(1+\tau)}}{M}\bigg]^{\frac{1}{2}}.
\eeqn 
Finally using Jensen inequality $E\sqrt{X}\le \sqrt{EX},$ and from Corollary \ref{cor} followed by Lemma \ref{l1}(b) denoting $c(\tau):= \sup_{k\ge 1}\sup_{M\ge 1} E\lt<|x|^{1+\tau},\bar{\eta}_{k-1}^{M}\rt>$ one gets 
\beqn
\sup_{k\ge 1}E\sum_{F\in \clp_{l}}\left[S^M(\bar{\eta}_{k-1}^{M})(2^{n}F\cap B_n)- \bar{\eta}_{k-1}^{M}(2^{n}F\cap B_n)\rt]&\le& 2^{\frac{dl}{2}} \sup_{k\ge 1} E\bigg[\frac{\lt<|x|^{1+\tau},\bar{\eta}_{k-1}^{M}\rt>2^{-(n-1)(1+\tau)}}{M}\bigg]^{\frac{1}{2}}\non\\
\le  2^{\frac{dl}{2}} \sup_{k\ge 1} \bigg[\frac{E\lt<|x|^{1+\tau},\bar{\eta}_{k-1}^{M}\rt>2^{-(n-1)(1+\tau)}}{M}\bigg]^{\frac{1}{2}}&\le& 2^{\frac{dl}{2}} \bigg[\frac{\sup_{k\ge 1} E\lt<|x|^{1+\tau},\bar{\eta}_{k-1}^{M}\rt>2^{-(n-1)(1+\tau)}}{M}\bigg]^{\frac{1}{2}}\non\\
&\le& 2^{\frac{dl}{2}}  \bigg[\frac{ c(\tau)2^{-(n-1)(1+\tau)}}{M}\bigg]^{\frac{1}{2}}.
\eeqn
 Hence the conclusion about the upper bound of $E\mathcal{W}_{1}(S^M(\bar{\eta}_{k-1}^{M}),\bar{\eta}_{k-1}^{M})$ will be similar to the first term of (\ref{lemthm3e0}). It will be a function of the sample size of the concentration gradient $M$ in place of $N$ in the bound of $E\clw_{1}(\bar{\mu}_{k}^{N},\bar{\mu}_{k-1}^{N}Q^{\bar{\eta}_{k-1}^{M},\bar{\mu}_{k-1}^{N}})$. Combining this with the conclusion about the first term of (\ref{lemthm3e0}) we can state the bound in terms of $N_{1}=\min\{M,N\}$ and the result of Lemma \ref{lemthm3.1} will follow.

$\square$

Now we will complete the theorem. Observe the following identity
\beqn
(\bar{\mu}_{n}^{N},\bar{\eta}_{n}^{M})-(\mu_{n},\eta_{n}) = \sum_{i=1}^{n} \bigg[\Psi^{(n-i)}(\bar{\mu}_{i}^{N},\bar{\eta}_{i}^{M})-\Psi^{(n-i)}\circ \Psi(\bar{\mu}_{i-1}^{N},\bar{\eta}_{i-1}^{M})\bigg]+\big[\Psi^{n}(\bar{\mu}_{0}^{N},\bar{\eta}_{0}^{M})-\Psi^{n}(\mu_{0},\eta_{0})\big].\non
\eeqn
Using Triangular inequality and Lemma \ref{lem1} following holds
\beqn
&&\clw_1\lt((\bar{\mu}_{n}^{N},\bar{\eta}_{n}^{M}),(\bar{\mu}_{n},\bar{\eta}_{n}) \rt)\non\\
&\le& \sum_{i=1}^{n}\clw_1\lt(\Psi^{(n-i)}(\bar{\mu}_{i}^{N},\bar{\eta}_{i}^{M}),\Psi^{(n-i)}\circ \Psi(\bar{\mu}_{i-1}^{N},\bar{\eta}_{i-1}^{M})\rt)+\clw_1\lt(\Psi^{n}(\bar{\mu}_{0}^{N},\bar{\eta}_{0}^{M}),\Psi^{n}(\mu_{0},\eta_{0})\rt)\non\\
&\le&\sum_{i=1}^{n}\theta^{n-i}\bigg[a\clw_1\lt((\bar{\mu}_{i}^{N},\bar{\eta}_{i}^{M}),\Psi(\bar{\mu}_{i-1}^{N},\bar{\eta}_{i-1}^{M})\rt)+ b\mathcal{W}_{1}\lt(\bar{\mu}^{(i-1)}_{M,N}Q^{\bar{\eta}_{i}^{M},\bar{\mu}^{(i-1)}_{M,N}}, \bar{\mu}^{(i-1)}_{M,N}Q^{\bar{\eta}_{i-1}^{M}R^{\alpha}_{\bar{\mu}_{i-1}^{N}}, \bar{\mu}^{(i-1)}_{M,N}}\rt)\bigg] \non\\ &&+\theta^{n}\lt[a\clw_1\lt((\bar{\mu}_{0}^{N},\bar{\eta}_{0}^{M}),(\mu_{0},\eta_{0})\rt)+b\mathcal{W}_{1}(\mu_{0}Q^{\bar{\eta}_{0}^{M},\mu_{0}}, \mu_{0}Q^{\eta_{0},\mu_{0}})\rt]\s\s\s\s\s\s\s\s\s\label{thm3.1}
\eeqn
where (\ref{thm3.1}) follows from (\ref{use}) with specified constants $a$ and $b$ and $\bar{\mu}^{(i-1)}_{M,N}:=\bar{\mu}^{N}_{i-1}Q^{\bar{\eta}_{i-1}^{M},\bar{\mu}^{N}_{i-1}}$. Let $X_{i}^{M,N}$ be a random variable, conditioned on $\mathcal{F}_{i-1}^{M,N},$ sampled from $\bar{\mu}^{(i-1)}_{M,N}.$ We have
\beqn
&&\mathcal{W}_{1}\lt(\bar{\mu}^{(i-1)}_{M,N}Q^{\bar{\eta}_{i}^{M},\bar{\mu}^{(i-1)}_{M,N}}, \bar{\mu}^{(i-1)}_{M,N}Q^{\bar{\eta}_{i-1}^{M}R^{\alpha}_{\bar{\mu}_{i-1}^{N}}, \bar{\mu}^{(i-1)}_{M,N}}\rt)\non\\&\le& \sup_{g\in \mbox{Lip}_{1}(\R^{d})} E\bigg|g(AX_{i}^{M,N}+\delta f_{\delta}(\gr \bar{\eta}_{i}^{M},\bar{\mu}^{(i-1)}_{M,N},X_{i}^{M,N},\epsilon)) -g(AX_{i}^{M,N}\non\\&&+\delta f_{\delta}(\gr(\bar{\eta}_{i-1}^{M}R^{\alpha}_{\bar{\mu}_{i-1}^{N}}),\bar{\mu}^{(i-1)}_{M,N},X_{i}^{M,N},\epsilon))\bigg| \leq \delta \sigma E\lt[\left|\gr \bar{\eta}_{i}^{M}(X_{i}^{M,N}) - \gr\bar{\eta}_{i-1}^{M}R^{\alpha}_{\bar{\mu}_{i-1}^{N}}(X_{i}^{M,N})\right|\big| \mathcal{F}_{i-1}^{M,N}\rt]\non\\&=&(1-\alpha)\int\bigg|\int\lt[S^{M}(\bar{\eta}^{M}_{i-1})-\bar{\eta}_{i-1}^{M} \rt](dx)\gr_{y}P(x,y)\bigg|(\bar{\mu}^{N}_{i-1}Q^{\bar{\eta}_{i-1}^{M}})(dy)\non\\&\le& l_{P}^{\gr}(1-\alpha)\mathcal{W}_{1}\lt(S^{M}(\bar{\eta}^{M}_{i-1}),\bar{\eta}_{i-1}^{M}\rt).\s\s\s\label{thm3.1e2}
\eeqn
Last display follows from Assumption \ref{As3}. Since $\bar{\eta}_{0}^{M}=\eta_0,$ one has
\beqn
\mathcal{W}_{1}(\mu_{0}Q^{\bar{\eta}_{0}^{M},\mu_{0}}, \mu_{0}Q^{\eta_{0},\mu_{0}})=0.\label{thm3.1e3}
\eeqn
Combining the results (\ref{thm3.1e2}),(\ref{thm3.1e3}), with (\ref{thm3.1}) we get for each $n,$
\beqn
E\clw_1\lt((\bar{\mu}_{n}^{N},\bar{\eta}_{n}^{M}),(\mu_n,\eta_n) \rt)&\le&\frac{a}{1-\theta}\sup_{k\ge 1}E\mathcal{W}_1\big((\bar{\mu}_{k}^{N},\bar{\eta}_{k}^{M}),\Psi(\bar{\mu}_{k-1}^{N},\bar{\eta}_{k-1}^{M})\big)\non\\+\frac{b l_{P}^{\gr}(1-\alpha)}{1-\theta}\sup_{k \ge 1}E\mathcal{W}_{1}\lt(S^{M}(\bar{\eta}^{M}_{k-1}),\bar{\eta}_{k-1}^{M}\rt)&+&a\theta^{n}E\clw_1\lt((\bar{\mu}_{0}^{N},\bar{\eta}_{0}^{M}),(\mu_{0},\eta_{0})\rt).\label{final}
\eeqn
 Using Lemma \ref{lemthm3.1} the result follows.

$\square$
\subsection{Proof of Corollary \ref{cor2}:}
Using triangular inequality and from (\ref{use}) one gets
\beqn
&&E \mathcal{W}_{1}\lt((\bar{\mu}_{n}^{N},\bar{\eta}_{n}^{M}),(\mu_{\infty},\eta_{\infty})\rt)  \le\mathcal{W}_{1}\big((\mu_{n},\eta_{n}),(\mu_{\infty},\eta_{\infty})\big)+ E \mathcal{W}_{1}\lt((\bar{\mu}_{n}^{N},\bar{\eta}_{n}^{M}),(\mu_{n},\eta_{n})\rt)\non\\
&\le&  \theta^{n}\bigg[a\mathcal{W}_{1}((\mu_{0},\eta_{0}),(\mu_{\infty},\eta_{\infty}))+b \mathcal{W}_{1}(\mu_{0}Q^{\eta_{0}}, \mu_{0}Q^{\eta_{\infty}})\bigg]+E \mathcal{W}_{1}\lt((\bar{\mu}_{n}^{N},\bar{\eta}_{n}^{M}),(\mu_{n},\eta_{n})\rt).
\eeqn
Combining this with (\ref{final}) we get
\beqn
&& E \mathcal{W}_{1}\lt((\bar{\mu}_{n}^{N},\bar{\eta}_{n}^{M}),(\mu_{\infty},\eta_{\infty})\rt)  \le \theta^{n}\bigg[a\mathcal{W}_{1}((\mu_{0},\eta_{0}),(\mu_{\infty},\eta_{\infty}))+b \mathcal{W}_{1}(\mu_{0}Q^{\eta_{0}}, \mu_{0}Q^{\eta_{\infty}})\bigg]\non\\&+&
\frac{a}{1-\theta}\sup_{k\ge 1}E\mathcal{W}_1\big((\bar{\mu}_{k}^{N},\bar{\eta}_{k}^{M}),\Psi(\bar{\mu}_{k-1}^{N},\bar{\eta}_{k-1}^{M})\big)+\frac{b l_{P}^{\gr}(1-\alpha)}{1-\theta}\sup_{k \ge 1}E\mathcal{W}_{1}\lt(S^{M}(\bar{\eta}^{M}_{k-1}),\bar{\eta}_{k-1}^{M}\rt).\non
\eeqn
The result is obvious after using Lemma \ref{lemthm3.1}. 

$\square$

\subsection{Proof of Theorem \ref{steadystate}:}
Fix $N$ and $M$. Define $\Theta_{n}^{N,M} \in \mathcal{P}((\mathbb{R}^{d})^{N}\times \clp_{1}^{*}(\R^{d})\times \clp_{}^{}(\R^{d}))$ as 
\begin{eqnarray}
\langle \phi,\Theta_{n}^{N,M}  \rangle &=& \frac{1}{n} \sum_{j=1}^{n} E \phi\big(\bar{X}_{j}(N), \eta_{j}^{M},S^{M}(\eta_{j}^{M})\big), \quad \phi\in BM\big((\mathbb{R}^{d})^{N}\times \clp_{1}^{*}(\R^{d}\big)\times \clp_{}^{}(\R^{d}))\,\,
\end{eqnarray}
for $N\geq 1, M\geq 1$ and $n\in \mathbb{N}_{0}$ where $\{(\bar{X}_{j}(N),\bar{\eta}_{j}^{M},S^{M}(\bar{\eta}_{j}^{M})),\,\,  j\in \mathbb{N}_{0},i=1,\ldots,N\}$ are as defined in the context of $\mathbb{IPS}_{2}$.
Note that $(\mathbb{R}^{d})^{N}\times \clp_{1}^{*}(\R^{d})\times \clp_{}^{}(\R^{d})$ is a complete separable metric space with metric $d((x,\mu_{1},\mu_{3}),(y,\mu_{2},\mu_{4})):=\|x-y\|+\frac{1}{2}\clw_{1}(\mu_1,\mu_2)+\frac{1}{2}\clw_{1}\big(\mu_{3},\mu_{4}\big)$ where $\|x\|:=\frac{1}{N}\sum_{i=1}^{N}|x_{i}|$ for $x=(x_1,\ldots,x_{N})\in(\mathbb{R}^{d})^{N}$. From Lemma \ref{l0} and \ref{l1} it follows that, for each $N,M\geq 1,$ the sequence $\{\Theta_{n}^{N,M},n\geq 1\}$ is relatively compact (By Prohorov's Theorem) and using Assumption \ref{As2} it is easy to see that any limit point $\Theta_{\infty}^{N,M}$ of $\Theta_{n}^{N,M}$ (as $n\to \infty$) is an invariant measure of the Markov chain $\{X_n(N),\bar{\eta}_{n}^{M},S^{M}(\bar{\eta}_{n}^{M})\}_{n\ge 0}$ and from Lemma \ref{l0} it satisfies $\int_{(\mathbb{R}^d)^N\times  \clp_{1}^{*}(\R^{d})\times \clp(\R^d)} |x|\; \Theta_{\infty}^{N,M}(dx) < \infty$ (Taking the norm of the product space as $|(x,y,z)|=\|x\|+\frac{1}{2}\|y\|_{1}+\frac{1}{2}\|z\|_{1}$ where $(x,y,z)\in (\mathbb{R}^{d})^{N}\times \clp_{1}^{*}(\R^{d})\times\clp_{}^{}(\R^{d})$ ). Uniqueness of invariant measure can be proved by the following simple coupling argument (see for example \cite{budhiraja2014long}): Suppose
$ \Theta_{\infty}^{N,M}$, $\tilde{\Theta}_{\infty}^{N,M}$ are two invariant measures that satisfy
$\int_{(\mathbb{R}^d)^N\times  \clp_{1}^{*}(\R^{d})\times \clp(\R^d)} |x|\; \Theta_{\infty}^{N,M}(dx) < \infty$, $\int_{(\mathbb{R}^d)^N\times  \clp_{1}^{*}(\R^{d})\times \clp(\R^d)} |x|  \tilde{\Theta}_{\infty}^{N,M}(dx) < \infty$.

Let $\big(X_{0}(N), \eta_{0}^{M},S^{M}(\eta_{0}^{M})\big)$ and $\big(\tilde{X}_0(N),\tilde{\eta}_{0}^{M},S^{M}(\tilde{\eta}_{0}^{M})\big)$ with
probability laws $\Theta_{\infty}^{N,M}$ and $\tilde{\Theta}_{\infty}^{N,M}$ respectively be given on a common probability space under same noise sequence (i.e in which an i.i.d. array of $\mathbb{R}^m$ valued random variables $\{\epsilon^i_n,  i= 1, \ldots, N, n \ge 1\}$ are defined 
that is independent of $(X_0(N),\eta_{0}^{M}, \tilde X_0(N),\tilde{\eta}_{0}^{M})$ with common probability law $\theta$) and the evolution equations are following.
\begin{align*}
X_{n+1}^{i} &= AX_{n}^{i} + \delta f_{\delta}(X_{n}^{i}, \gr \eta_{n}^{M}(X_{n}^{i}),\mu_{n}^{N},\epsilon_{n+1}^{i}),\; \mu_{n}^{N} = \frac{1}{N}\sum_{i=1}^N \delta_{X_{n}^{i}}, \,\,\eta_{k}^{M}= (1-\alpha)(S^{M}(\eta_{k-1}^{M})P) +\alpha \mu^{N}_{k-1}P', \\
\tilde X_{n+1}^{i} &= A\tilde X_{n}^{i} +  \delta f_{\delta}(\tilde X_{n}^{i},\gr \tilde{\eta}_{n}^{M}(\tilde X_{n}^{i}), \tilde\mu_{n}^{N},\epsilon_{n+1}^{i}),\; \tilde\mu_{n}^{N} = \frac{1}{N}\sum_{i=1}^N \delta_{\tilde X_{n}^{i}}, \,\, \tilde{\eta}_{k}^{M}= (1-\alpha)(S^{M}(\tilde{\eta}_{k-1}^{M})P) +\alpha \tilde{\mu}^{N}_{k-1}P',
\end{align*}
where recall $f_{\delta}(\cdot,\cdot,\cdot,x)=f(\cdot,\cdot,\cdot,x)+\frac{B(x)}{\delta}.$ Note that  
\beqn
\clw_{1}\big(\frac{1}{N}\sum_{i=1}^{N}\delta_{X_{i}},\frac{1}{N}\sum_{i=1}^{N}\delta_{Y_{i}}\big)\le \frac{1}{N}\sum_{i=1}^{N}|X_{i}-Y_{i}|\label{wass}
\eeqn
for any two arrays $\{X_i\}_{i=1}^{N}$ and $\{Y_i\}_{i=1}^{N}$. Using the independence of the noise sequence along with (\ref{wass}) and Assumption  \ref{As2} we have
\beqn
E|X_{n+1}^{i}- \tilde X_{n+1}^{i}| &\le& (\|A\|+ \delta \sigma) E|X_{n}^{i}- \tilde X_{n}^{i}|
+ \delta \sigma \frac{1}{N} \sum_{j=1}^N E|X_{n}^{j}- \tilde X_{n}^{j}|\non\\&&+\delta \sigma E|\gr\eta_{n}^{M}(X_{n}^{i}) - \gr\tilde{\eta}_{n}^{M}(\tilde{X}_{n}^{i})|.\label{recc1}
\eeqn
Now applying Assumption \ref{As3} (doing similar calculations as in \eqref{th2eq1},\eqref{fin B_{1}},\eqref{B_{2}}) following inequality holds
\beqn
E|\gr\eta_{n}^{M}(X_{n}^{i}) - \gr\tilde{\eta}_{n}^{M}(\tilde{X}_{n}^{i})|&\le& E|\gr\eta_{n}^{M}(X_{n}^{i}) - \gr\eta_{n}^{M}(\tilde{X}_{n}^{i})|+E|\gr\eta_{n}^{M}(\tilde{X}_{n}^{i})- \gr\tilde{\eta}_{n}^{M}(\tilde{X}_{n}^{i})|\label{stead1}\\
\le l_{pp'}^{\gr}E|X_{n}^{i}-\tilde{X}_{n}^{i}|&+&\alpha l_{P'}^{\gr}E\clw_{1}\big(\mu_{n-1}^{N},\tilde{\mu}_{n-1}^{N}\big)+(1-\alpha) l_{P}^{\gr}E\clw_{1}\big(S^{M}(\eta_{n-1}^{M}),S^{M}(\tilde{\eta}_{n-1}^{M})\big).\non
\eeqn

Note that (\ref{wass}) implies 
\beqn
E\big[\clw_{1}(S^{M}(\eta^{M}_{k-1}),S^{M}(\tilde{\eta}^{M}_{k-1}))\big|\mathcal{F}_{k-1}^{M,N}\big]\le\clw_{1}(\eta^{M}_{k-1},\tilde{\eta}^{M}_{k-1})\label{wass2}
\eeqn
from which following holds from \eqref{stead1}
\beqn
E|\gr\eta_{n}^{M}(X_{n}^{i}) - \gr\tilde{\eta}_{n}^{M}(\tilde{X}_{n}^{i})|&\le& l_{pp'}^{\gr}E|X_{n}^{i}-\tilde{X}_{n}^{i}|+\alpha l_{P'}^{\gr}E|X_{n-1}^{i}-\tilde{X}_{n-1}^{i}|\non\\&&+(1-\alpha) l_{P}^{\gr}E\clw_{1}\big(\eta_{n-1}^{M},\tilde{\eta}_{n-1}^{M}\big).\label{recc2}
\eeqn
We also have
\beqn
\clw_{1}\big(\eta^{M}_{n+1},\tilde{\eta}^{M}_{n+1}\big)\le (1-\alpha)l(P)\clw_{1}\big(S^{M}(\eta^{M}_{n}),S^{M}(\tilde{\eta}^{M}_{n})\big)+\alpha l(P')\clw_{1}\big(\mu_{n}^{N},\tilde{\mu}_{n}^{N}\big)
\eeqn
and after taking expectation
\beqn
E\clw_{1}\big(\eta^{M}_{n+1},\tilde{\eta}^{M}_{n+1}\big)&\le& (1-\alpha)l(P)E\clw_{1}\big(\eta^{M}_{n},\tilde{\eta}^{M}_{n}\big)+\alpha l(P')E|X_{n}^{i}-\tilde{X}_{n}^{i}|.\label{recc3}
\eeqn
Letting
$A^{(M,N)}_{n+1}:= \frac{1}{N}\sum_{i=1}^N |X_{n+1}^{i}- \tilde X_{n+1}^{i}|+\clw_{1}\big(\eta^{M}_{n+1},\tilde{\eta}^{M}_{n+1}\big)$, we have the following recursion relation combining (\ref{recc1}),(\ref{recc2}) and (\ref{recc3})
\beqn
EA^{(M,N)}_{n+1}&\le&\max\Big\{\Big(\|A\|+\delta \sigma(2+ l_{PP'}^{\gr,\alpha}) + \alpha l(P')\Big ),(1-\alpha)l(P) \Big\} EA^{(M,N)}_{n}\non\\&&+\delta\sigma\max\{(1-\alpha)l_{P}^{\gr},\alpha l_{P'}^{\gr}\}EA^{(M,N)}_{n-1}
\eeqn
which is the same recursion as in (\ref{cond2}). Now for the chosen $\delta,\alpha$ satisfying (\ref{cond}) there exists a $\theta\in (0,1)$ such that
\beqn
EA^{(M,N)}_{n}\le\theta^{n-1}[EA^{(M,N)}_{0}+EA^{(M,N)}_{1}].\label{recmn}
\eeqn
 Also, since $\Theta_{\infty}^{N,M}$ and $\tilde \Theta_{\infty}^{N,M}$ are invariant distributions, for every $n \in \mathbb{N}_0$, $\big(X_{n+1}(N),\eta_{n+1}^{M},S^{M}(\eta_{n+1}^{M})\big)$ is distributed as $\Theta_{\infty}^{N,M}$ and
$\big(\tilde{X}_{n+1}(N),\tilde{\eta}_{n+1}^{M},S^{M}(\tilde{\eta}_{n+1}^{M})\big)$ is distributed as $\tilde{\Theta}_{\infty}^{M,N}$.
Thus \newline $(X_{n+1}(N),\eta_{n+1}^{M},S^{M}(\eta_{n+1}^{M}))$ and $\big(\tilde{X}_{n+1}(N),\tilde{\eta}_{n+1}^{M},S^{M}(\tilde{\eta}_{n+1}^{M})\big)$ define a coupling of random variables with laws $\Theta_{\infty}^{N,M}$ and
$\tilde \Theta_{\infty}^{N,M}$ respectively. From \eqref{recmn} we then have
$$\mathcal{W}_1(\Theta_{\infty}^{N,M}, \tilde \Theta_{\infty}^{M,N}) \le Ed\big((X_{n}(N),\eta_{n}^{M},S^{M}(\eta_{n}^{M})),(\tilde{X}_{n}(N),\tilde{\eta}_{n}^{M},S^{M}(\tilde{\eta}_{n}^{M}))\big) \le EA_{n}^{M,N}\to 0,$$ as $n\to\infty.$ So there exists a unique invariant measure $\Theta_{\infty}^{N,M} \in \mathcal{P}_{1}\big((\mathbb{R}^{d})^{N}\times \clp_{1}^{*}(\R^{d})\times\clp(\R^d)\big)$ for this Markov chain and, as $n \to \infty$,
\begin{eqnarray}
\Theta_{n}^{N,M} \rightarrow \Theta_{\infty}^{N,M}.\label{4.5.2}
\end{eqnarray}
This proves the first part of the theorem.  
Denote $\Theta^{N,M}_{\infty}\lt(\cdot,\clp^{*}_1(\R^d),\clp(\R^d)\rt)$ by $\Theta^{1,N,M}_{\infty}$ and \newline $\Theta^{N,M}_{n}\lt(\cdot,\clp^{*}_1(\R^d),\clp(\R^d)\rt)$ by $\Theta^{1,N,M}_{n}$.

Define $r_{N}:(\mathbb{R}^{d})^{N} \to \mathcal{P}(\mathbb{R}^{d})$  as 
$$ r_{N}(x_{1},\ldots, x_{N}) = \frac{1}{N}\sum_{i=1}^{N}\delta_{x_{i}},\s\s (x_1, \ldots x_N) \in (\R^d)^N.$$
Let $\nu_{n}^{N,M} = \Theta_{n}^{1,N,M}  \circ r_{N}^{-1}$ and $\nu_{\infty}^{N,M} = \Theta_{\infty}^{1,N,M}  \circ r_{N}^{-1}$. In order to prove that $\Theta_{\infty}^{1,N,M}$ is $\mu_{\infty}$-chaotic, it suffices to argue that (cf. \cite{sznitman1991topics})
\begin{equation}
	\nu^{N,M}_{\infty} \to \delta_{\mu_{\infty}} \mbox { in } \mathcal{P}(\mathcal{P}(\mathbb{R}^{d})), \mbox{ as } N,M \to \infty.
\label{eq:eq727}	
\end{equation}
We first argue that as $n\to \infty$ 
\begin{equation}\nu_{n}^{N,M} \rightarrow \nu_{\infty}^{N,M}\quad\quad \text{in $\mathcal{P}(\mathcal{P}(\mathbb{R}^{d}))$}.\label{eq:eq728}\end{equation}
It suffices to show that $\langle F,\nu_{n}^{N,M}\rangle \to \langle F,\nu_{\infty}^{N,M}\rangle$ for any continuous and bounded function $F:\mathcal{P}(\mathbb{R}^{d}) \to \mathbb{R}$.
But this is immediate on observing that 
$$\langle F, \nu_n^{N,M} \rangle = \langle F\circ r_N , \Theta_{n}^{1,N,M} \rangle, \; \langle F, \nu_{\infty}^{N,M} \rangle = \langle F\circ r_N , \Theta_{\infty}^{1,N,M}\rangle,$$
the continuity of the map $r_N$ and the weak convergence of $\Theta_{n}^{N,M}$ to $\Theta_{\infty}^{N,M}$.
Next, for any 
$f \in BL_{1}(\mathcal{P}(\mathbb{R}^{d}))$
$$
\Big|\langle f,\nu_{n}^{N,M} \rangle - \langle f,\delta_{\mu_{\infty}} \rangle\Big| = \Big| \frac{1}{n} \sum_{j=1}^n Ef(\bar{\mu}_j^{N}) - f(\mu_{\infty})\Big| 
\le \frac{1}{n} \sum_{j=1}^n E \mathcal{W}_{1}(\bar{\mu}_{j}^{N},\mu_{\infty}).$$
Fix $\epsilon > 0$.  For every $N,M \in \mathbb{N}$ there exists $n_0(N,M) \in \N$ such that for all $n \ge n_0(N,M)$
$$
E \mathcal{W}_{1}(\bar{\mu}_{n}^{N},\mu_{\infty}) \le \limsup_{n \to \infty} E \mathcal{W}_{1}(\bar{\mu}_{n}^{N},\mu_{\infty}) + \epsilon.$$
Thus for all $n, N,M \in \mathbb{N}$
\begin{equation}
	|\langle f,\nu_{n}^{N,M} \rangle - \langle f,\delta_{\mu_{\infty}} \rangle| \le \frac{n_0(N,M)}{n} \max_{1 \le j \le n_0(N,M)} E \mathcal{W}_{1}(\bar{\mu}_{j}^{N},\mu_{\infty})
	+ \limsup_{n \to \infty} E \mathcal{W}_{1}(\bar{\mu}_{n}^{N},\mu_{\infty}) + \epsilon.
	\label{eq:eq740}
	\end{equation}
Finally
\begin{align*}
	\limsup_{N,M \to \infty} |\langle f,\nu_{\infty}^{N,M} \rangle - \langle f,\delta_{\mu_{\infty}} \rangle|
	= & 	\limsup_{\min{\{N,M\}} \to \infty} \lim_{n \to \infty}|\langle f,\nu_{n}^{N,M} \rangle - \langle f,\delta_{\mu_{\infty}} \rangle|  \\
	\le & \limsup_{\min{\{N,M\}} \to \infty} \limsup_{n \to \infty} E \mathcal{W}_{1}(\bar{\mu}_{n}^{N},\mu_{\infty}) + \epsilon\\
	\le & \epsilon,
\end{align*}
where the first equality is from \eqref{eq:eq728}, the second uses \eqref{eq:eq740} and the third is a consequence of Corollary \ref{cor2}.  Since $\epsilon > 0$ is arbitrary, we have
\eqref{eq:eq727}	and the result follows.

$\square$

\subsection{Proof of Concentration bounds:}
\subsubsection{Proof of Theorem \ref{conv} \eqref{polyconv}:}
We start with the following lemma where we  establish a concentration bound for $\mathcal{W}_1\big((\bar{\mu}_{n}^{N},\bar{\eta}_{n}^{M}),\Psi(\bar{\mu}_{n-1}^{N},\bar{\eta}_{n-1}^{M})\big)$ for each fixed time $n\in\mathbb{N}$ and then combine it with the estimate in \eqref{thm3.1} in order to get the desired result.  
\begin{lemma}
	\label{lempoly}
Let $N_{1}=\min\{M,N\}.$ Assumptions (\ref{As2}-\ref{As3}) and Assumptions (\ref{unifas1}),(\ref{unifas2}) hold for some $\tau>0$. Suppose that
 $\delta \in(0,a(\tau)^{\frac{1}{1+\tau}})$, and $(1-\alpha)l_{\tau}(P)<1.$ Then there exist \newline $a_1, a_2, a_3,a'_1, a'_2, a'_3  \in (0, \infty)$ such that for all $\eps, R > 0, n \in \mathbb{N}$, and $N_{1} \ge \max\{1, a_1 (\frac{R}{\eps})^{d+2}\}.$ 
	\begin{eqnarray}
	P[\mathcal{W}_1\big((\bar{\mu}_{n}^{N},\bar{\eta}_{n}^{M}),\Psi(\bar{\mu}_{n-1}^{N},\bar{\eta}_{n-1}^{M})\big) >\varepsilon] &\leq& a_3 \left(e^{-a_2 \frac{N_{1} \varepsilon^{2}}{R^{2}}} + \frac{R^{-\tau}}{\varepsilon}\right),\label{lempoly1}\\
P[\mathcal{W}_{1}\lt(S^{M}(\bar{\eta}^{M}_{n-1}),\bar{\eta}_{n-1}^{M}\rt) >\varepsilon] &\leq& a'_3 \left(e^{-a'_2 \frac{N_{1} \varepsilon^{2}}{R^{2}}} + \frac{R^{-\tau}}{\varepsilon}\right).\label{lempoly2}
	\end{eqnarray}
	 
\end{lemma}
\subsubsection{Proof of Lemma \ref{lempoly}}   
Second concentration bound will follow by proceeding as Lemma 4.5 of \cite{budhiraja2014long}. The proof relies on an idea of restricting measures to a
compact set and estimates on metric entropy \cite{bolley2007quantitative} (see also \cite{villani2003topics}). The basic idea is to first obtain a concentration bound for the $\mathcal{W}_1$ distance between the truncated law and its corresponding empirical law in a compact ball of radius $R$ and getting a tail estimate from Lemma \ref{l1} and Corollary \ref{cor} after conditioning by $\mathcal{F}_{n-1}^{M,N}$. With the notations (for example $\mu_{R}$ is the truncated measure of $\mu$ restricted on a ball $B_{\R}(0)$ of $R$ radius) introduced in  Lemma 4.5 of \cite{budhiraja2014long} we sketch the proof of the second bound. With that notation the truncated version of $\bar{\eta}_{n-1}^{M}$  is denoted by $\bar{\eta}_{n-1,R}^{M}$. Suppoe $\{Y_{n-1}^{i,M}:i=1,\ldots,M\}$ are iid from $\bar{\eta}_{n-1}^{M}$ conditioned on $\mathcal{F}_{n-1}^{M,N}.$ where $\{Z_{i}^{M,R}:i=1,\ldots,M\}$ are iid from $\bar{\eta}_{n-1,R}^{M}$ conditioned under $\mathcal{F}_{n-1}^{M,N}.$ Define
\[
X_{n-1}^{i,M} = \left \{ 
\begin{array}{cc}
Y_{n-1}^{i,M} & \text{ when }|Y_{n-1}^{i,M}| \leq  R, \\ \ \\
Z_{n-1}^{i,M} & \text{ otherwise }.
\end{array} \right.
\]
Note that $P(X_{n-1}^{i,M} \in A \mid \mathcal{F}_{n-1}^{M,N}) = P(Z_{n-1}^{i,M} \in A \mid \mathcal{F}_{n-1}^{M,N}).$ Denote $S_{}^{M}(\bar{\eta}^{M}_{n-1,R}):=\frac{1}{M}\sum_{i=1}^{M}\delta_{X_{n-1}^{i,M}}.$
Now denoting $ a(1+\tau):=\sup_{n\ge 0} \sup_{M,N}E\lt<|x|^{1+\tau},\bar{\eta}^{M}_{n}\rt>$, from \eqref{wass} we have
\beqn
P\big[\mathcal{W}_{1}(S_{}^{M}(\bar{\eta}^{M}_{n-1,R}),S^{M}(\bar{\eta}^{M}_{n-1})) >\frac{\e}{3}\big]&\le& 3\frac{E[\mathcal{W}_{1}(S_{}^{M}(\bar{\eta}^{M}_{n-1,R}),S^{M}(\bar{\eta}^{M}_{n-1}))]}{\e}\non\\
\le\frac{3}{\e}E E\big[|X_{n-1}^{i,M}-Y_{n-1}^{i,M}|\big|\mathcal{F}_{n-1}^{M,N}\big]&=&\frac{3}{\e}E E\big[|Z_{n-1}^{i,M}-Y_{n-1}^{i,M}|1_{|Y_{n-1}^{i,M}|>R}\big|\mathcal{F}_{n-1}^{M,N}\big]\non\\
\le \frac{6}{\e}EE\big[|Y_{n-1}^{i,M}|1_{|Y_{n-1}^{i,M}|>R}\big|\mathcal{F}_{n-1}^{M,N}\big]&\le& 6 a(1+\tau)\frac{R^{-\tau}}{\e}.\label{lempolye1}
\eeqn

Now using Azuma Hoeffding inequality as done in display (4.35) of Lemma 4.5 in \cite{budhiraja2014long} one has 
\beqn
P\big[\mathcal{W}_{1}(S_{}^{M}(\bar{\eta}^{M}_{n-1,R}),\bar{\eta}_{n-1,R}^{M})>\frac{\e}{3}\big]\le \max\left \{2, \frac{16R}{\varepsilon}(2\sqrt{d}+1)3^{[\frac{8R}{\varepsilon}(\sqrt{d}+1)]^{d}  }\right\} e^{ -\frac{M \varepsilon^{2}}{288 R^{2}}}.\label{lempolye2}
\eeqn

From the definition of $\bar{\eta}_{n-1,R}^{M}$
\beqn
P\big[\mathcal{W}_{1}(\bar{\eta}_{n-1}^{M},\bar{\eta}_{n-1,R}^{M})\ge\frac{\e}{3}\big]\le\frac{6}{\e} E\big[|Y_{n-1}^{i,M}|1_{|Y_{n-1}^{i,M}|>R}\big]\le3 a(1+\tau)\frac{R^{-\tau}}{\e}.\label{lempolye3}
\eeqn
Using triangular inequality 
\begin{eqnarray}
\mathcal{W}_{1}(S^{M}(\bar{\eta}^{M}_{n-1}),\bar{\eta}_{n-1}^{M}) \leq \mathcal{W}_{1}(S_{}^{M}(\bar{\eta}^{M}_{n-1,R}),S^{M}(\bar{\eta}^{M}_{n-1})) + \mathcal{W}_{1}(S_{}^{M}(\bar{\eta}^{M}_{n-1,R}),\bar{\eta}_{n-1,R}^{M}) + \mathcal{W}_{1}(\bar{\eta}_{n-1}^{M},\bar{\eta}_{n-1,R}^{M})\;\;\;\non
\end{eqnarray}
combining \eqref{lempolye1},\eqref{lempolye2} and \eqref{lempolye3} the result \eqref{lempoly2} will follow.

 The first one \eqref{lempoly1} follows by noting that 
\beqn
P\big[\mathcal{W}_1\big((\bar{\mu}_{n}^{N},\bar{\eta}_{n}^{M}),\Psi(\bar{\mu}_{n-1}^{N},\bar{\eta}_{n-1}^{M})\big)>\varepsilon\big]\le  P\big[\mathcal{W}_1\big(\bar{\mu}^{N}_{n},\bar{\mu}_{n-1}^{N}Q^{\bar{\mu}_{n-1}^{N},\bar{\eta}_{n-1}^{M}}\big)> \frac{\epsilon}{2}\big]\label{ftermlempoly}\\+P\bigg[\mathcal{W}_1\big(S^{M}(\bar{\eta}_{n-1}^{M}),\bar{\eta}_{n-1}^{M}\big)>\frac{\epsilon}{2(1-\alpha)l(P)}\bigg].\non
\eeqn
Proceeding like Lemma 4.5 of \cite{budhiraja2014long} the bound for the first term in RHS of \eqref{ftermlempoly} can be established.

\qed

\subsubsection{Proof of Theorem \ref{conv}(a)}
Combining \eqref{thm3.1},\eqref{thm3.1e2} and \eqref{thm3.1e3} it follows that
\beqn
\mathcal{W}_{1}\lt((\bar{\mu}_{n}^{N},\bar{\eta}_{n}^{M}),(\mu_{n},\eta_{n})\rt)\s\le\s\sum_{i=1}^{n}\theta^{n-i}\bigg[a\clw_1\lt((\bar{\mu}_{i}^{N},\bar{\eta}_{i}^{M}),\Psi(\bar{\mu}_{i-1}^{N},\bar{\eta}_{i-1}^{M})\rt)\non\\ + b l_{P}^{\gr}(1-\alpha)\mathcal{W}_{1}\lt(S^{M}(\bar{\eta}^{M}_{i-1}),\bar{\eta}_{i-1}^{M}\rt)\bigg] +a \theta^{n}\clw_1\lt((\bar{\mu}_{0}^{N},\bar{\eta}_{0}^{M}),(\mu_{0},\eta_{0})\rt).\label{mart}
\eeqn

 Denoting $c_{1}:=\max\left\{\left(\big(\|A\|+\delta \sigma(2+ l_{PP'}^{\gr,\alpha})\big) + \alpha l(P')\right),(1-\alpha)l(P)\right\}$, $c_{2}:=\delta \sigma \max\big\{\alpha l_{P'}^{\gr},(1-\alpha) l_{P}^{\gr}\big\}$ define the function $g_{0}(\cdot)$ as
$$g_{0}(\gamma):= c_{2}+(1-\gamma)c_{1}-(1-\gamma)^{2}.$$
Since $g_{0}(0)=c_{2}+c_{1}-1<0$ (from the assumption), $g_{0}(1)=c_{2}>0$ and  $g(\cdot)$ is continuous. So there exists a $\gamma>0$ such that $g_{0}(\gamma)<0$ or equivalently 
$$\frac{c_{1}}{1-\gamma}+\frac{c_{2}}{(1-\gamma)^2}<1.$$

So there exists a $\theta\in(0,1-\gamma)$ such that statement of Lemma \ref{lem1} holds. Now using that $\gamma$ from \eqref{mart} one has
\begin{eqnarray}
&&P\Big[\mathcal{W}_{1}\lt((\bar{\mu}_{n}^{N},\bar{\eta}_{n}^{M}),(\mu_{n},\eta_{n})\rt)>\varepsilon\Big] \leq P\Big[\bigcup_{i=1}^{n}\big\{ a\theta^{n- i}\clw_1\lt((\bar{\mu}_{i}^{N},\bar{\eta}_{i}^{M}),\Psi(\bar{\mu}_{i-1}^{N},\bar{\eta}_{i-1}^{M})\rt)>\frac{\gamma}{2} (1-\gamma)^{n-i}\varepsilon \big\}\nonumber\\ &&\bigcup_{i=1}^{n}\big\{b l_{P}^{\gr}(1-\alpha)\theta^{n-i}\mathcal{W}_{1}\lt(S^{M}(\bar{\eta}^{M}_{i-1}),\bar{\eta}_{i-1}^{M}\rt)>\frac{\gamma}{2} (1-\gamma)^{n-i}\varepsilon\big\}\bigcup_{i=1}^{n}\big\{\theta^{n}\clw_1\lt((\bar{\mu}_{0}^{N},\bar{\eta}_{0}^{M}),(\mu_{0},\eta_{0})\rt) \non\\&&> \gamma (1-\gamma)^{n}\varepsilon\big\}\Big] 
\s\le\s \sum_{i=1}^{n}P\big[\clw_1\lt((\bar{\mu}_{i}^{N},\bar{\eta}_{i}^{M}),\Psi(\bar{\mu}_{i-1}^{N},\bar{\eta}_{i-1}^{M})\rt) >\frac{\gamma \varepsilon}{2a} \big(\frac{1-\gamma}{\theta}\big)^{n-i}\big]+\non\\ &&\sum_{i=1}^{n}P\big[\mathcal{W}_{1}\lt(S^{M}(\bar{\eta}^{M}_{i-1}),\bar{\eta}_{i-1}^{M}\rt)>\frac{\gamma\varepsilon}{2b l_{P}^{\gr}(1-\alpha)}\big(\frac{1-\gamma}{\theta}\big)^{n-i}\big]+P\big[\clw_1\lt((\bar{\mu}_{0}^{N},\bar{\eta}_{0}^{M}),(\mu_{0},\eta_{0})\rt)  >\gamma\varepsilon\Big(\frac{1-\gamma}{\theta}\Big)^{n}\big].\non
\label{unif}
\end{eqnarray}
Let $\beta_{1}=\frac{\gamma \varepsilon}{2a},\,\,\beta_{2}= \frac{\gamma\varepsilon}{2b l_{P}^{\gr}(1-\alpha)}\,\,\beta_{3}=\gamma\varepsilon .$ Note that $\nu:=\big(\frac{1-\gamma}{\theta}\big) >1,$ from our choice of $\gamma$.  Therefore denoting $\beta:=\min\{\beta_{1},\beta_{2}\},\s$
$N_{1}\ge a_{1}\Big(\frac{R}{\beta}\Big)^{d+2}\vee 1$ implies $N_{1}\ge a_{1}\Big(\frac{R}{\beta \nu^{n}}\Big)^{d+2}\vee 1$ for all $n\in\mathbb{N}_{0}$ and a consequence of Lemma \ref{lempoly} gives 
\beqn
&&P\Big[\mathcal{W}_{1}\lt((\bar{\mu}_{n}^{N},\bar{\eta}_{n}^{M}),(\mu_{n},\eta_{n})\rt)>\varepsilon\Big] \leq \sum_{i=1}^{n}P\big[\clw_1\lt((\bar{\mu}_{i}^{N},\bar{\eta}_{i}^{M}),\Psi(\bar{\mu}_{i-1}^{N},\bar{\eta}_{i-1}^{M})\rt) >\beta_{1} \nu^{n-i}\big]\non\\
&+&\sum_{i=1}^{n}P\big[\mathcal{W}_{1}\lt(S^{M}(\bar{\eta}^{M}_{i-1}),\bar{\eta}_{i-1}^{M}\rt)>\beta_{2}\nu^{n-i}\big]+P\big[\clw_1\lt((\bar{\mu}_{0}^{N},\bar{\eta}_{0}^{M}),(\mu_{0},\eta_{0})\rt)  >\beta_{3}\nu^{n}\big]\label{finex}\\
&\le& a_3 \sum_{i=1}^{n}\left(e^{-a_2 \frac{N_{1} \beta^{2}_{} \nu^{2i}}{R^{2}}} + \frac{R^{-\tau}}{\beta \nu^{i}}\right)+a'_3 \sum_{i=1}^{n}\left(e^{-a'_2 \frac{N_{1} \beta^{2}_{} \nu^{2i}}{R^{2}}} + \frac{R^{-\tau}}{\beta \nu^{i}}\right)\non\\&&+P\big[\clw_1\lt((\bar{\mu}_{0}^{N},\bar{\eta}_{0}^{M}),(\mu_{0},\eta_{0})\rt)  >\beta_{3}\nu^{n}\big].\non
\eeqn
Now proceeding similarly like the proof of Theorem 3.7 of \cite{budhiraja2014long} through  optimizing the value of $R$ the conclusion will follow.

\subsubsection{Proof of Theorem \ref{conv}(b)}
Second part regarding the exponential concentration bound will follow similarly (like Theorem 3.8 of \cite{budhiraja2014long}) under the following lemmas on uniform exponential integrability.

\begin{lemma}
	\label{lemexpconv}
	Suppose Assumptions \ref{exponas1} and \ref{exponas2} hold. Suppose there exists $\alpha^{*}>0$ such that $$\alpha^{*}|h_{1}(0)|+h_{2}(\alpha^*)=-\log(1-\alpha).$$ Then for all $\alpha_1 \in [0, \min\big\{\alpha^{*},\frac{\alpha(\delta)}{\delta}\big\})$ and 
	$\delta \in \Big[0,\frac{1- \|A\|}{(2+l^{\gr,\alpha}_{PP'})K}\Big)$, 	
\beqn
\sup_{n\ge 0} \sup_{M,N \ge 1} E e^{\alpha_1|X_n^{1,N}|} < \infty, \s\s \sup_{n\ge 0}\sup_{M,N \ge 1}E\lt<e^{\alpha_{1}|x|},\bar{\eta}^{M}_{n}\rt><\infty. \label{exponlem}
\eeqn
\end{lemma}
{\bf Proof.}
We will start by proving the second inequality. Note that from Corollary \ref{cor} the conditions for $``\sup_{n\ge 0}\sup_{M,N \ge 1}E\lt<e^{\alpha_{1}|x|},\bar{\eta}^{M}_{n}\rt><\infty$" are same as the conditions for $\sup_{n\ge 0}\sup_{N \ge 1}E\lt<e^{\alpha_{1}|x|},\eta^{N}_{n}\rt><\infty$ in $\mathbb{IPS}_{1}$ and from Lemma \ref{l1} they are again same as the conditions for finiteness of $\sup_{n\ge 0}\lt<e^{\alpha_{1}|x|},\eta_{n}\rt>.$ Note that
\beqn\label{etaltau11}
\left<\eta_{k+1},e^{\alpha_{1}|x|}\right> = \sum_{i=0}^{k}\left[\alpha(1-\alpha)^{i}\left<\mu_{k-i}P'P^{i},e^{\alpha_{1}|x|}\right>\right]+(1-\alpha)^{k+1}\left<\eta_{0}P^{k+1},e^{\alpha_{1}|x|}\right>.
\eeqn
Now from Assumption \ref{exponas2}, using lipshitz property $|h_{1}(x)|\le l_{h_{1}}|x|+|h_{1}(0)|$ one has $\left<\mu_{k}P'P^{i},e^{\alpha_{1}|x|}\right>\le e^{\alpha_{1}|h_{1}(0)|+h_{2}(\alpha)} \left<\mu_{k}P'P^{i-1},e^{\alpha_{1}l_{h_{1}}|x|}\right>+e^{h_{3}(\alpha_1)+h_{2}(\alpha_1)}.$ So we have an upperbound of $\big<\mu P' P^{i},e^{\alpha_{1}|x|}\big>$  that is
\beqn
&\le& e^{\sum_{k=0}^{i-1}\big[h_{2}(\alpha_{1}l^{k}_{h_{1}})+\alpha_{1}l^{k}_{h_{1}}|h_{1}(0)|\big]}\big<\mu P',e^{\alpha_{1}l^{i}_{h_{1}}|x|}\big>+\sum_{k=0}^{i-1}e^{h_{3}(\alpha_{1}l^{k}_{h_{1}})+h_{2}(\alpha_{1}l^{k}_{h_{1}})+\sum_{m=1}^{k\vee 1}\big(h_{2}(\alpha_{1}l^{m-1}_{h_{1}})+\alpha l^{m-1}_{h_{1}}|h_{1}(0)|\big)}\non\\
&\le& \big<\mu P',e^{\alpha_{1}|x|}\big>e^{i\big(h_{2}(\alpha_{1})+\alpha_{1}|h_{1}(0)|\big)}+ \frac{e^{(i-1)\big(h_{2}(\alpha_{1})+\alpha_{1}|h_{1}(0)|\big)}-1}{e^{h_{2}(\alpha_{1})+\alpha_{1}|h_{1}(0)|} -1}e^{h_{3}(\alpha_{1})+h_{2}(\alpha_{1})}.\non
\eeqn
Last inequality follows since $h_{2}(\cdot), h_{3}(\cdot)$ are non-decreasing and $l_{h_{1}}\le 1.$
Using \eqref{etaltau11} under the condition $\sup_{n\ge 0}  \lt< e^{\alpha_1|x|},\mu_{n}\rt> < \infty$ (which we prove shortly) we conclude that $\sup_{k\ge 0}\left<\eta_{k+1},e^{\alpha_{1}|x|}\right><\infty$ or equivalently $\sum_{i=0}^{\infty}(1-\alpha)^{i}e^{i\big[h_{2}(\alpha)+\alpha_{1}|h_{1}(0)|\big]}<\infty$ if there exists an $\alpha_{1}$ such that $\alpha_{1}|h_{1}(0)|+h_{2}(\alpha_{1})+\log(1-\alpha)<0.$ Since $g(\alpha_{1}):=\alpha_{1}|h_{1}(0)|+h_{2}(\alpha_{1}),$ is an increasing function of $\alpha_{1}$ and $g(0)=0$. From the definition of $\alpha^*$ we can always find $0<\alpha_{1}<\alpha^*$ such that $\sup_{n\ge 0}\lt<e^{\alpha_{1}|x|},\eta_{n}\rt><\infty.$

 Now we prove $\sup_{n\ge 0}  \lt< e^{\alpha_1|x|},\mu_{n}\rt> < \infty$ or equivalently the first term in \eqref{exponlem}. Note that from \eqref{l1eq1} for $n \ge 1$
\beqn
|X_{n+1}^{i}| &\le& \|A\| |X_{n}^{i}| + \delta A_1(\epsilon_{n+1}^{i})[|\gr \eta_{n}^{N}(X_{n}^{i})|+\|\mu_{n}^{N}\|_{1}+| X_{n}^{i}|]+\delta A_2(\epsilon_{n+1}^{i})+|B(\epsilon_{n+1}^{i})|\non\\
&\le& |X_{n}^{i}|\Big[\|A\|+\delta K\Big( 1+l^{\gr,\alpha}_{PP'}\Big)\Big]+\delta K \|\mu_{n}^{N}\|_{1}+\delta\Big(A_{2}(\epsilon_{n+1}^{i})+Kc_{PP'}^{\alpha}+\frac{B(\epsilon_{n+1}^{i})}{\delta}\Big).\non
\eeqn

Now from the choice $\alpha_{1}\le\frac{\alpha(\delta)}{\delta},$ taking expectation after having exponential 
\beqn
Ee^{\alpha_{1}|X_{n+1}^{i}|}\le  Ee^{\alpha_{1}|X_{n}^{i}|\big[\|A\|+\delta K\big( 1+l^{\gr,\alpha}_{PP'}\big)\big]+\alpha_{1}\delta K \|\mu_{n}^{N}\|_{1}} \cle_1(\alpha_{1})\label{exponma}
\eeqn
where  $\cle_1(\alpha_1 )=e^{\alpha_{1}\delta K c_{PP'}^{\alpha}}\int e^{\alpha_{1}\delta \big(A_{2}(z) + \frac{|B(z)|}{\delta}\big)} \theta(dz).$ We note that from Assumption \ref{exponas2} there always exist $\alpha^{**}<\frac{\alpha(\delta)}{\delta},\s  c_{3}$ such that for all $\alpha_{1}\in(0,\alpha^{**})$
\beqn
\cle_1(\alpha_1 )\le e^{c_{3}\alpha_{1}}.\label{exponmain}
\eeqn

Using conditioning argument we have
\beqn
Ee^{\alpha_{1}|X_{n}^{i}|\big[\|A\|+\delta K\big( 1+l^{\gr,\alpha}_{PP'}\big)\big]+\alpha_{1}\delta K \|\mu_{n}^{N}\|_{1}}&=&EE\bigg[e^{\alpha_{1}|X_{n}^{i}|\big[\|A\|+\delta K\big( 1+l^{\gr,\alpha}_{PP'}\big)\big]+\alpha_{1}\delta K \|\mu_{n}^{N}\|_{1}}\bigg|\sigma\Big(\frac{1}{N}\sum_{i=1}^{N}\delta_{X^{i,N}_{n}}\Big)\bigg]\non\\
&=&E\bigg[e^{\alpha_{1}\delta K \|\mu_{n}^{N}\|_{1}}E\Big[e^{\alpha_{1}|X_{n}^{i}|\big[\|A\|+\delta K\big( 1+l^{\gr,\alpha}_{PP'}\big)\big]}\bigg|\sigma\Big(\frac{1}{N}\sum_{i=1}^{N}\delta_{X^{i,N}_{n}}\Big)\Big]\bigg]\non\\
&=&E\bigg[e^{\alpha_{1}\delta K \|\mu_{n}^{N}\|_{1}}\frac{1}{N}\sum_{i=1}^{N}e^{\alpha_{1}|X^{i,N}_{n}|\big[\|A\|+\delta K\big( 1+l^{\gr,\alpha}_{PP'}\big)\big]}\bigg]\label{exponlem3}
\eeqn
where  \eqref{exponlem3} follows from exchangeability of $\{X_{n}^{i,N}\}_{i=1,\ldots,N}$. Observing $\|\mu_{n}^{N}\|_{1}=\int|x|\mu_{n}^{N}(dx)$ and using Jensen's inequality applied to the function $x\to e^{\alpha_{1}\delta Kx},$ we have after taking expectation
\beqn
E\bigg[e^{\alpha_{1}\delta K \|\mu_{n}^{N}\|_{1}}\frac{1}{N}\sum_{i=1}^{N}e^{\alpha_{1}|X^{i,N}_{n}|\big[\|A\|+\delta K\big( 1+l^{\gr,\alpha}_{PP'}\big)\big]}\bigg]\le E\bigg[\frac{1}{N}\sum_{i=1}^{N}e^{\alpha_{1}\delta K |X^{i,N}_{n}|}\frac{1}{N}\sum_{i=1}^{N}e^{\alpha_{1}|X^{i,N}_{n}|\big[\|A\|+\delta K\big( 1+l^{\gr,\alpha}_{PP'}\big)\big]}\bigg].\non
\eeqn
Since $f_{1}(x):=e^{\alpha_{1}\delta Kx}$ and $f_{2}(x):=e^{\alpha_{1}x\big[\|A\|+\delta K\big( 1+l^{\gr,\alpha}_{PP'}\big)\big]}$ are both non-decreasing, so putting $\mu=\mu_{n}^{N}$ almost surely in the following inequality $\int f_{1}(x)f_{2}(x) \mu(dx) \geq \int f_{1}(x) \mu(dx) \int f_{2}(y) \mu(dy)$  and taking expectation we have
\beqn
 E\bigg[\frac{1}{N}\sum_{i=1}^{N}e^{\alpha_{1}\delta K |X^{i,N}_{n}|}\frac{1}{N}\sum_{i=1}^{N}e^{\alpha_{1}|X^{i,N}_{n}|\big[\|A\|+\delta K\big( 1+l^{\gr,\alpha}_{PP'}\big)\big]}\bigg] &\le& E \frac{1}{N}\sum^{N}_{i=1}e^{\alpha_{1}|X_{n}^{i,N}|\big[\|A\|+\delta K\big(2+l^{\gr,\alpha}_{PP'}\big)\big]}\non\\
&\le& Ee^{\alpha_{1}|X_{n}^{i}|\big[\|A\|+\delta K\big(2+l^{\gr,\alpha}_{PP'}\big)\big]}.\non
\eeqn
From our choice of $\delta,$ $\kappa:=\|A\|+\delta K\big(2+l^{\gr,\alpha}_{PP'}\big)\in (0,1).$ Denoting $F_{n+1}(\alpha_{1}):=Ee^{\alpha_{1}|X_{n+1}^{i}|}$ from \eqref{exponma} we have the following recursive inequality: 
\beqn
F_{n+1}(\alpha_{1})\le F_{n}(\alpha_{1}\kappa)\cle_1(\alpha_1 ).
\eeqn
 Iterating the above inequality we have for all $n \ge 1$
$$F_n(\alpha_1) \le F_0(\alpha_1) \prod_{j=0}^{n-1} \cle_1(\alpha_1  \kappa_1^j)
\le F_0(\alpha_1) e^{c_3\alpha_1  \sum_{j=0}^{n-1}\kappa_1^j} \le F_0(\alpha_1) e^{c_3\alpha_1 /(1-\kappa_1)}$$
where the second inequality is a consequence of \eqref{exponmain}.

Note further for the system in \eqref{extra} let $\{X_{n}\}_{n\in\mathbb{N}_{0}}$ be defined as the random variables with laws $\mathcal{L}(X_{n}):=\mu_{n}$ for $n\in\mathbb{N}_{0}.$ Then starting similarly from 
$$|X_{n+1}^{}| \le |X_{n}^{}|\Big[\|A\|+\delta K\Big( 1+l^{\gr,\alpha}_{PP'}\Big)\Big]+\delta K \|\mu_{n}^{}\|_{1}+\delta\Big(A_{2}(\epsilon_{n+1}^{})+Kc_{PP'}^{\alpha}+\frac{B(\epsilon_{n+1}^{})}{\delta}\Big)$$
using the inequality $\int f_{1}(x)f_{2}(x) \mu(dx) \geq \int f_{1}(x) \mu(dx) \int f_{2}(y) \mu(dy)$ (similar to Lemma 4.11 of \cite{budhiraja2014long}) one can prove 
\beqn
\sup_{n\ge 0}\lt<e^{\alpha_{1}|x|},\mu_{n}\rt>\le \lt<e^{\alpha_{1}|x|},\mu_{0}\rt> e^{\frac{c_3\alpha_1}{1-\kappa_1}}.\label{expononlin}
\eeqn
under same conditions on $\delta, \alpha_{1}.$  This is needed for proving $\sup_{n\ge 0}\lt<e^{\alpha_{1}|x|},\eta_{n}\rt><\infty.$ The result follows.

\qed

\begin{lemma}\label{expolem2}
	Then there exist $a_1, a_2, a_3,a'_1, a'_2, a'_3  \in (0, \infty)$ such that for all $\eps, R > 0$ and $n \in \mathbb{N}$, and $N_{1} \ge \max\{1, \ti a_1 (\frac{R}{\eps})^{d+2}\}$
	\begin{eqnarray}
	P[\mathcal{W}_1\big((\bar{\mu}_{n}^{N},\bar{\eta}_{n}^{M}),\Psi(\bar{\mu}_{n-1}^{N},\bar{\eta}_{n-1}^{M})\big) >\varepsilon] &\leq& a_3 \left(e^{-a_2 \frac{N_{1} \varepsilon^{2}}{R^{2}}} +  \ti B_{1}(\alpha_1)\frac{e^{-\alpha_1 R}}{\varepsilon}\right),\label{lempoly1}\\
P[\mathcal{W}_{1}\lt(S^{M}(\bar{\eta}^{M}_{n-1}),\bar{\eta}_{n-1}^{M}\rt) >\varepsilon] &\leq& a'_3 \left(e^{-a'_2 \frac{N_{1} \varepsilon^{2}}{R^{2}}} +  \ti B_{2}(\alpha_1)\frac{e^{-\alpha_1 R}}{\varepsilon}\right).\label{lempoly2}
	\end{eqnarray}
	\end{lemma}
\subsubsection{Proof of Lemma \ref{expolem2}:}
Follows from similar decompositions given in Lemma \ref{lempoly} and Lemma 4.7 of \cite{budhiraja2014long}.
\qed
\subsubsection{Proof of Theorem \ref{conv}(b):}
Starting from \eqref{mart}, the conclusion will follow by applying Lemma \ref{expolem2} in \eqref{finex}.

\qed
\vspace{1cm}

\subsection{Proof of Theorem \ref{conviid}}
We will start by introducing a coupling. Consider a system of $\R^d$ valued auxiliary random variables $\{Y_{n}^{i,N},i=1,\ldots,N\}_{n\ge 0}$ defined as follows. 
\beqn
Y_{n+1}^{i,N} &=& AY_{n}^{i,N} + \delta f(\gr \eta_{n}^{}(Y_{n}^{i,N}), \mu_{n},Y_{n}^{i,N},\epsilon_{n+1}^{i})+B(\epsilon_{n+1}^{i}),\text{\hspace{10 mm}} \quad i=1,\ldots,N, \quad n\in\mathbb{N}_{0}.\non\\
\eta_{n+1}&=& \eta_{n}R^{\alpha}_{\mu_{n}},\non\\
Y_{0}^{i,N}&=&X_{0}^{i,N}.\label{aux}
\eeqn
Now for each $n\in\mathbb{N},\s\{Y_{n}^{i,N},i=1,\ldots,N\}$ is a set of $\R^d$ valued iid random variables under initial assumption $\mathcal{L}(\{X^{i,N}_{0}\}_{i=1,\ldots,N})=\mu^{\otimes N}_{0}.$ Suppose $\zeta_{n}^{N}:=\frac{1}{N}\sum_{i=1}^{N}\delta_{Y^{i,N}_{n}}$. The following Lemma will make a connection between $\zeta_{n}^{N}$ and $\mu_{n}^{N}$.

\begin{lemma}(Coupling with the auxiliary system)\label{couplingin}
Suppose Assumptions  \ref{As2},\ref{As3},\ref{As7} and \ref{exponas1} hold. Then for every $n \ge 0$ and $N \ge 1,$ with the $C_{1},$ and $\chi_{1}$ defined in \eqref{cou1},\eqref{cou2}
	\begin{eqnarray}
	\mathcal{W}_{1}(\mu_{n+1}^{N},\mu_{n+1}) &\leq& \mathcal{W}_{1}(\zeta_{n+1}^{N},\mu_{n+1}) +C_{1}\sum_{k=0}^{n} \chi_{1}^{n-k}\mathcal{W}_{1}(\zeta_{k}^{N},\mu_{k}).\label{Gron}
	\end{eqnarray}
\end{lemma}

\textbf{Proof.}
Since by Assumption \ref{As2} and  $A_{1}(\epsilon)\leq K$, we have for each $j=1,\ldots,N$
\begin{eqnarray*}
|X_{n+1}^{j} - Y_{n+1}^{j,N}| &\leq& \|A\| |X_{n}^{j} - Y_{n}^{j,N}| + \delta K \big\{|\gr\eta_{n}^{N}(X_{n}^{j,N}) - \gr\eta_{n}(Y_{n}^{j,N})|+|X_{n}^{j,N} - Y_{n}^{j,N}| +\mathcal{W}_{1}(\mu_{n}^{N},\mu_{n})\big\}\non\\
\end{eqnarray*}
Using the calculations in \eqref{gradeq1},\eqref{th2eq1},\eqref{B_{1}} and \eqref{B_{2}} 
\beqn
|\gr\eta_{n}^{N}(X_{n}^{j,N}) - \gr\eta_{n}(Y_{n}^{j,N})|&\le&|\gr\eta_{n}^{N}(X_{n}^{j,N}) - \gr\eta_{n}^{N}(Y_{n}^{j,N})|+|\gr\eta_{n}^{N}(Y_{n}^{j,N}) - \gr\eta_{n}(Y_{n}^{j,N})|\non\\
&\le& l_{PP'}^{\gr,\alpha}|X_{n}^{j} - Y_{n}^{j,N}|+(1-\alpha)l_{P}^{\gr}\clw_{1}(\eta_{n-1}^{N},\eta_{n-1})+\alpha l^{\gr}_{P'}\clw_{1}(\mu_{n-1}^{N},\mu_{n-1})\non
\eeqn
Thus 
\beqn
|X_{n+1}^{j,N} - Y_{n+1}^{j,N}| &\le& \big[\|A\|+\delta K(1+l^{\gr,\alpha}_{PP'}) \big]|X_{n}^{j} - Y_{n}^{j,N}|+\delta K\bigg[\clw_{1}(\mu_{n}^{N},\mu_{n})\non\\&&+(1-\alpha)l^{\gr}_{P}\clw_{1}(\eta_{n-1}^{N},\eta_{n-1})+\alpha l^{\gr}_{P'}\clw_{1}(\mu_{n-1}^{N},\mu_{n-1})\bigg]\label{recfin}
\eeqn
Using \eqref{recfin} as the recursion on $a^{j}_{n+1}:=|X_{n+1}^{j,N} - Y_{n+1}^{j,N}|$ with $a^{j}_{0}=0,$ we get 
\beqn
a^{j}_{n+1}\le \delta K \sum_{k=1}^{n}\big[\|A\|+\delta K(1+l^{\gr,\alpha}_{PP'}) \big]^{n-k}\bigg[\clw_{1}(\mu_{k}^{N},\mu_{k})+(1-\alpha)l^{\gr}_{P}\clw_{1}(\eta_{k-1}^{N},\eta_{k-1})\non\\+\alpha l^{\gr}_{P'}\clw_{1}(\mu_{k-1}^{N},\mu_{k-1})\bigg].\label{recfin2}
\eeqn
Denote $\|A\|+\delta K(1+l^{\gr,\alpha}_{PP'})$ by $\chi$.
Observe that 
\beqn
\clw_{1}(\eta_{n-1}^{N},\eta_{n-1})=(1-\alpha)l(P)\clw_{1}(\eta_{n-2}^{N},\eta_{n-2})+\alpha l(P')\clw_{1}(\mu_{n-2}^{N},\mu_{n-2}).\label{recfin3}
\eeqn 
Denote the quantity in the third bracket of RHS of \eqref{recfin2} by $b_{k}.$ Using \eqref{recfin3} and $\eta_{0}^{N}=\eta_{0}$ we have
\beqn
b_{k}&=&\clw_{1}(\mu_{k}^{N},\mu_{k})+(1-\alpha)l^{\gr}_{P}\clw_{1}(\eta_{k-1}^{N},\eta_{k-1})+\alpha l^{\gr}_{P'}\clw_{1}(\mu_{k-1}^{N},\mu_{k-1})\non\\
&=&\clw_{1}(\mu_{k}^{N},\mu_{k})+(1-\alpha)l^{\gr}_{P}\alpha l(P')\sum_{i=0}^{k-2}[(1-\alpha)l(P)]^{k-2-i}\clw_{1}(\mu_{i}^{N},\mu_{i})+\alpha l^{\gr}_{P'}\clw_{1}(\mu_{k-1}^{N},\mu_{k-1})\non\\
&\le& c_{4} \sum_{i=0}^{k}c^{k-i}_{5}\clw_{1}(\mu_{i}^{N},\mu_{i}).
\eeqn
 where $c_{4}:=\max\{1,(1-\alpha)l^{\gr}_{P}\alpha l(P')\}$ and $c_{5}:=\max\{\alpha l^{\gr}_{P'},(1-\alpha)l(P)\}.$
Thus from \eqref{recfin2} we have
\beqn
a^{j}_{n+1}\le \delta K c_{4} \sum_{k=0}^{n} \chi^{n-k} \sum_{i=0}^{k} c^{k-i}_{5}\clw_{1}(\mu_{i}^{N},\mu_{i}).
\eeqn
 Now applying Lemma \ref{app42} we have 
\beqn
a^{j}_{n+1}&\le& \delta K c_{4} \sum_{i=0}^{n}  \clw_{1}(\mu_{i}^{N},\mu_{i})\bigg[\frac{\chi^{n+1-i} - c^{n+1-i}_{5}}{\chi - c_{5}}\bigg]\non\\
&\le& \delta K c_{7}\sum_{i=0}^{n}\chi_{2}^{n+1-i}  \clw_{1}(\mu_{i}^{N},\mu_{i}) \label{recfin4}
\eeqn
where $\chi_{2}:=\max\{\chi,c_{5}\}$ and $c_{7}:= \frac{c_{4}}{|\chi-c_{5}|}.$ Note that from \eqref{wass}
we have for all $n\ge 0,$ $$\clw_{1}(\zeta_{n}^{N},\mu_{n}^{N})\le \frac{1}{N}\sum_{j=1}^{N}a^{j}_{n}.$$
Combining the result above and using triangle inequality in \eqref{recfin4} 
\beqn
\clw_{1}(\zeta_{n+1}^{N},\mu_{n+1}^{N})&\leq& \delta K c_{7} \sum_{k=0}^{n}\chi_{2}^{n+1-k}\clw_{1}(\zeta_{k}^{N},\mu_{k}^{N})+\delta K c_{7} \sum_{k=0}^{n}\chi_{2}^{n+1-k}\mathcal{W}_{1}(\zeta_{k}^{N},\mu_{k}).\nonumber 
\eeqn
Applying Lemma \ref{app42} with
$$a_n = \chi_{2}^{-n}\mathcal{W}_{1}(\zeta_{n}^{N},\mu_{n}^{N}), \; b_n = \delta K c_{7} \sum_{k=0}^{n-1} \chi_{2}^{-k} \mathcal{W}_{1}(\eta_{k}^{N},\mu_{k}), \;
p_n = \delta K c_{7}, \; n \ge 0.$$
We have 
\begin{eqnarray}
\chi_{2}^{-(n+1)}\mathcal{W}_{1}(\zeta_{n+1}^{N},\mu_{n+1}^{N}) &\leq& b_{n+1}+  \sum_{k=0}^{n}(\delta K c_{7})^{2} \sum_{i=0}^{k-1} \chi_{2}^{-i}\mathcal{W}_{1}(\zeta_{i}^{N},\mu_{i}) \left(1+\delta K c_{7}\right)^{n-k}\nonumber\\
&=& b_{n+1}+ \sum_{i=0}^{n}\sum_{k=i+1}^{n}(\delta K c_{7})^{2}(1+\delta K c_{7})^{n-k}\chi_{2}^{-i}\mathcal{W}_{1}(\zeta_{i}^{N},\mu_{i}) \nonumber\\
&=& b_{n+1}+\sum_{i=0}^{n} (\delta K c_{7})^{2}\chi_{2}^{-i}.\mathcal{W}_{1}(\zeta_{i}^{N},\mu_{i})\sum_{m=0}^{n-i-1} (1+\delta K c_{7})^{m}\nonumber\\
&=& b_{n+1}+\sum_{i=0}^{n} (\delta K c_{7})\chi_{2}^{-i}\mathcal{W}_{1}(\zeta_{i}^{N},\mu_{i})[(1+\delta K c_{7})^{n-i} -1]\label{Gron1}. 
\end{eqnarray}
Simplifying (\ref{Gron1}) one gets 
\begin{eqnarray}
\mathcal{W}_{1}(\zeta_{n+1}^{N},\mu_{n+1}^{N}) &\leq&\delta K c_{7} \sum_{k=0}^{n} \chi_{2}^{n+1-k} \mathcal{W}_{1}(\zeta_{k}^{N},\mu_{k}) +\sum_{k=0}^{n} (\delta K c_{7})\chi_{2}^{n+1-k}\mathcal{W}_{1}(\zeta_{k}^{N},\mu_{k})[(1+\delta K c_{7})^{n-k} -1]\nonumber\\
&=&\sum_{k=0}^{n} (\delta K c_{7})\chi_{2}^{n+1-k}\mathcal{W}_{1}(\zeta_{k}^{N},\mu_{k})(1+\delta K c_{7})^{n-k} .\nonumber\\
&=&\delta K c_{7}\chi_{2}\sum_{k=0}^{n} (\chi_{2}+\delta K c_{7}\chi_{2})^{n-k}\mathcal{W}_{1}(\zeta_{k}^{N},\mu_{k}) .\nonumber
\end{eqnarray}
Note that $\delta K c_{7}\chi_{2} = C_{1}$ and $\chi_{2}+C_{1}=\chi_{1}$ as defined in \eqref{cou1} \and \eqref{cou2} respectively. Thus we have
$$\mathcal{W}_{1}(\zeta_{n+1}^{N},\mu_{n+1}^{N}) \le C_{1}\sum_{k=0}^{n} \chi_{1}^{n-k}\mathcal{W}_{1}(\zeta_{k}^{N},\mu_{k}).$$
The result now follows by an application of triangle inequality.
\qed

\subsubsection{Proof of Theorem \ref{conviid}}
Since $\chi_{1}<1.$ So we can find $\gamma>0$ such that $\chi_{1}<1-\gamma.$ Taking that $\gamma,$ we have $\nu_{1}:=\frac{1-\gamma}{\chi_{1}}>1$.
For any $\e > 0$, From Lemma \ref{conviid}
\begin{eqnarray}
&&P[\mathcal{W}_{1}(\mu_{n}^{N},\mu_{n}) >\varepsilon] \leq P[\mathcal{W}_{1}(\zeta_{n}^{N},\mu_{n})> \gamma \varepsilon] +\sum_{i=0}^{n-1} P[C_{1}\chi_{1}^{n-1-i}\mathcal{W}_{1}(\zeta_{i}^{N},\mu_{i}) \ge \gamma\varepsilon (1-\gamma)^{n-i}]\quad\quad\nonumber\\
&=& P[\mathcal{W}_{1}(\zeta_{n}^{N},\mu_{n})> \gamma \varepsilon] +\sum_{i=1}^{n} P[\mathcal{W}_{1}(\zeta_{n-i}^{N},\mu_{n-i}) \ge \frac{\gamma\varepsilon\chi_{1}}{C_{1}} \nu^{i}]\label{unifexp1}\\
&=&P[\mathcal{W}_{1}(\zeta_{n}^{N},\mu_{n})> \gamma \varepsilon] +\sum_{i=1}^{i_{\varepsilon}} P[\mathcal{W}_{1}(\zeta_{n-i}^{N},\mu_{n-i}) \ge \frac{\gamma\chi_{1}\varepsilon}{C_{1}} \nu^{i}]+\sum_{i=i_{\varepsilon}+1}^{n} P[\mathcal{W}_{1}(\zeta_{n-i}^{N},\mu_{n-i}) \ge \frac{\gamma\chi_1\varepsilon}{C_{1}} \nu^{i}] ,\;\;\;\non
\end{eqnarray}
where $i_{\varepsilon}:=\max\{i\ge 0:  \frac{\gamma\chi_{1}\varepsilon}{C_{1}} \nu^{i}<1\}.$
Note that for $\delta \in \Big[0,\frac{1- \|A\|}{(2+l^{\gr,\alpha}_{PP'})K}\Big),$ and $\alpha_{1}\in (0,\frac{\alpha(\delta)}{\delta})$ from  \eqref{expononlin} we have $\sup_{n\ge 0}\lt<e^{\alpha_{1}|x|,\mu_{n}}\rt><\infty.$ That implies from the statement of Theorem 2 of \cite{fournier2013rate} that for all $N>0,$

$$P[\clw_{1}(\zeta_{n}^{N},\mu_{n})\ge \varepsilon]\le a(N,\varepsilon)1_{\{\varepsilon \le 1\}}+b(N,\varepsilon).$$
where $a(N,\varepsilon)=e^{-cN\varepsilon^{2}}1_{\{d=1\}}+e^{-cN\big(\frac{\varepsilon}{\log(2+\frac{1}{\varepsilon})}\big)^2}1_{\{d=2\}}+e^{-cN\varepsilon^{d}}1_{\{d>2\}}$
and $b(N,\varepsilon)=e^{-cN\varepsilon}.$ In order to prove \eqref{totalconc} we will prove only for one case $d>2.$ Rest will follow similarly. There exists $C'_{1},C'_{2},C'_{3}$

\beqn
\sum_{i=i_{\varepsilon}+1}^{n} P[\mathcal{W}_{1}(\zeta_{n-i}^{N},\mu_{n-i}) \ge \frac{\gamma\chi_1\varepsilon}{C_{1}} \nu^{i}]&\le& \sum_{i=i_{\varepsilon}+1}^{n} b(N,\frac{\gamma\chi_1\varepsilon}{C_{1}} \nu^{i})\le \sum_{i=i_{\varepsilon}+1}^{n} e^{-C'_{1}\e N\nu^{i}}\label{leq21}\\
\sum_{i=1}^{i_{\varepsilon}} P[\mathcal{W}_{1}(\zeta_{n-i}^{N},\mu_{n-i}) \ge \frac{\gamma\chi_1\varepsilon}{C_{1}} \nu^{i}]&\le&\sum_{i=1}^{i_{\varepsilon}} a(N,\frac{\gamma\chi_1\varepsilon}{C_{1}} \nu^{i})\le  \sum_{i=1}^{i_{\varepsilon}} e^{-C'_{2} N(\e\nu^{i})^{d}}\le \sum_{i=1}^{i_{\varepsilon}} e^{-C'_{2} N\e^{d} \nu^{i}}\;\;\;\;\;\;\label{leq22}\\
P[\mathcal{W}_{1}(\zeta_{n}^{N},\mu_{n})> \gamma \varepsilon]&\le&e^{-C'_{3}N \e^{d}\wedge\e}\label{leq23}
\eeqn

Suppose $k_{0}$ such that $\nu^{i}\ge k_{0} i$ for all $i\ge 1.$ Combining \eqref{leq21},\eqref{leq22},\eqref{leq23} we have for all $N>1$ and $a''_{2}=k_{0}\min\{C'_{1},C'_{2},C'_{3}\}.$
\beqn
\sup_{n}P[\mathcal{W}_{1}(\mu_{n}^{N},\mu_{n}) >\varepsilon]&\leq& \sum_{i=0}^{\infty}e^{-a''_{2} N i \e^{d}\wedge\e }\le\frac{e^{-a''_{2}N \e^{d}\wedge\e}}{1-e^{-a''_{2} N  \e^{d}\wedge\e}}.
\eeqn
Now there exists $N_{3}:=-\frac{1}{a''_{2}}\log(1-\frac{1}{a''_{1}})$ such that $N\ge N_{3}\max\{\frac{1}{\e},\frac{1}{\e^{d}}\}$ we have
$$\sup_{n}P[\mathcal{W}_{1}(\mu_{n}^{N},\mu_{n}) >\varepsilon]\leq  a''_{1}e^{-a''_{2} N \e^{d}\wedge\e}.$$

\qed

\section{Acknowledgements}
A part this article was part of author's Phd thesis. The author is thankful to Prof. Amarjit Budhiraja for his comments on an earlier version of the manuscript.






 \appendix
  \section*{Appendix}
  \renewcommand{\thesection}{A} 


The first part of the following lemma is an immediate consequence of Ascoli-Arzela theorem where as the second follows from Lemma 5 in \cite{como2009scaling}.
%
%
%
\begin{lemma}\label{app12}
(a) For a compact set $K$ in $\mathbb{R}^{d}$ let $\mathcal{F}_{a,b}(K)$ be the space of functions $f: K \to \mathbb{R}$ such that $\sup_{x\in K}|f(x)|\leq a$  and $|f(x) - f(y)|\leq b|x-y|$ for all $x,y \in K$. Then for any $\epsilon > 0$ there is a finite subset $\mathcal{F}_{a,b}^{\epsilon}(K)$ of $\mathcal{F}_{a,b}(K)$ such that for any signed measure $\mu$
$$\sup_{f \in \mathcal{F}_{a,b}(K)} |\langle f,\mu\rangle | \leq \max_{g \in \mathcal{F}_{a,b}^{\epsilon}(K)} |\langle g,\mu \rangle| + \epsilon |\mu|_{TV}.$$
\end{lemma}
The next lemma is straightforward.
\begin{lemma}\label{app22}
Let $P:\mathbb{R}^{d}\times \mathcal{B}(\mathbb{R}^d) \to [0,1]$ be a transition probability kernel. Fix $N\geq 1$ and let $y_{1},y_{2},...,y_{N} \in \mathbb{R}^{d}$. Let $X_{1},X_{2},...,X_{N}$ be independent random variables such that $\mathcal{L}(X_{i})=\delta_{y_{i}}P.$ Let $f\in BM(\mathbb{R}^{d})$ and let $m_{0}^{N}= \frac{1}{N}\sum_{i=1}^{N}\delta_{y_{i}}$, $m_{1}^{N}= \frac{1}{N}\sum_{i=1}^{N}\delta_{X_{i}}$. Then 
$$E|\langle f, m_{1}^{N}- m_{0}^{N}P  \rangle| \leq \frac{2\|f\|_{\infty}}{\sqrt{N}}.$$
\end{lemma}
The following is a discrete version of Gronwall's lemma.
\begin{lemma}\label{app42}
\begin{enumerate}[(a)]
\item Let $\{ a_{i}\}_{i=0}^{\infty},\{ b_{i}\}_{i=0}^{\infty},\{ p_{i}\}_{i=0}^{\infty}$ be  non-negative sequences. Suppose that
$$  a_{n}  \leq b_{n} + \sum_{k=0}^{n-1} p_{k} a_{k} \;   \mbox{ for all } n \ge 0.$$
Then 
$$ a_{n}  \leq b_{n} + \sum_{k=0}^{n-1}\left [ p_{k} b_{k}\left(\prod_{j=k+1}^{n-1}(1+ p_{j})\right) \right ]\;   \mbox{ for all } n \ge 0.$$
\item For any $a,b>0$ and $\{C_{i}\}_{i\ge 0}$ be a nonnegative sequence of elements, then for all $n\ge 0$
$$\sum_{k=0}^{n}a^{n-k} \sum_{i=0}^{k}b^{k-i}C_{i} = \sum_{i=0}^{n} C_{i}\Bigg[\frac{a^{n+1-i}-b^{n+1-i}}{a-b}\Bigg].$$
\end{enumerate}
\end{lemma}

\bibliographystyle{plain}
\bibliography{ei}

\skp

{\sc

\skp

\noi
Abhishek Pal Majumder\\
Department of Mathematical Science\\
University of Copenhagen\\
email: abhishek@math.ku.dk

}

\end{document}